\documentclass[10pt,oneside]{jm}
\usepackage{amsfonts,amssymb,amsmath,rotating} 

\newcommand{\abs}{\vskip 0.5em\noindent\rm}

\newcommand{\AbsT}[1]{\paragraph{\hspace{-1em} #1}\rm}

\newcommand{\Pf}{\paragraph*{Proof.}\rm}

\newcommand{\Rem}{\AbsT{Remark.}}

\newcommand{\PropT}[1]{\AbsT{Proposition: #1.}}

\newcommand{\Thm}{\AbsT{Theorem.}}
\newcommand{\ThmT}[1]{\AbsT{Theorem: #1.}}

\newcommand{\bfi}{\noindent {\bf i)} }

\newcommand{\bfii}{\noindent {\bf ii)} }

\newcommand{\bfa}{\noindent {\bf a)} }
\newcommand{\bfb}{\noindent {\bf b)} }

\newcommand{\absb}{\abs$\bu$ }
\newcommand{\absc}{\abs$\circ$ }

\newcommand{\QED}{\hfill $\sharp$}

\newcommand{\N}{\mathbb N}
\newcommand{\Z}{\mathbb Z}
\newcommand{\Q}{\mathbb Q}
\newcommand{\R}{\mathbb R}
\newcommand{\C}{\mathbb C}
\newcommand{\F}{\mathbb F}

\newcommand{\cA}{\mathcal A}

\newcommand{\cD}{\mathcal D}

\newcommand{\cF}{\mathcal F}

\newcommand{\cH}{\mathcal H}

\newcommand{\cL}{\mathcal L}

\newcommand{\cM}{\mathcal M}

\newcommand{\cP}{\mathcal P}

\newcommand{\cS}{\mathcal S}

\newcommand{\cU}{\mathcal U}

\newcommand{\al}{\alpha}

\newcommand{\bt}{\beta}

\newcommand{\dt}{\delta}
\newcommand{\Dt}{\Delta}
\newcommand{\eps}{\epsilon}

\newcommand{\lb}{\lambda}

\newcommand{\zt}{\zeta}

\newcommand{\Hom}{\text{Hom}}

\newcommand{\id}{\text{id}}

\newcommand{\modcat}[2]{\textbf{mod}_{#1}\text{-}#2}

\newcommand{\spc}{\rule{1.5em}{0em}}
\newcommand{\vspc}{\rule{0em}{1.2em}}
\newcommand{\baraut}{\ov{\text{\hspace{0.5em}\rule{0em}{0.55em}}}}

\newcommand{\ld}{,\ldots\hskip0em ,}

\newcommand{\ov}[1]{\overline{#1}}

\newcommand{\fl}[1]{\lfloor #1\rfloor}

\newcommand{\mt}{\mapsto}
\newcommand{\lmt}{\longmapsto}

\newcommand{\da}{\downarrow}

\newcommand{\ra}{\rightarrow}

\newcommand{\cn}{\colon}
\newcommand{\dcup}{\stackrel{.}{\cup}}

\newcommand{\bu}{\bullet}

\newcommand{\sseq}{\subseteq}
\newcommand{\smin}{\setminus}
\newcommand{\emp}{\emptyset}

\newcommand{\tm}{\times}
\newcommand{\otm}{\otimes}

\newcommand{\GAP}{{\sf GAP}}

\renewcommand{\PropT}[1]{\AbsT{Proposition: #1.}\it}
\renewcommand{\Thm}{\AbsT{Theorem.}\it}
\renewcommand{\ThmT}[1]{\AbsT{Theorem: #1.}\it}

\begin{document}
\raggedbottom
\pagestyle{myheadings}
\markboth{}{}
\thispagestyle{empty}
\setcounter{page}{1}
\begin{center} \Large\bf
On low-degree representations \vspace*{0.2em} \\
of the symmetric group \vspace*{1em} \\
\large\rm
J\"urgen M\"uller \vspace*{1em} \\
\end{center}

\begin{abstract} \noindent
The aim of the present paper is to obtain a classification of all the
irreducible modular representations of the symmetric group on $n$ letters
of dimension at most $n^3$, including dimension formulae. This is achieved
by improving an idea, originally due to G. James, to get hands on dimension
bounds, by building on the current knowledge about decomposition numbers of
symmetric groups and their associated Iwahori-Hecke algebras, and by
employing a mixture of theory and computation. 

\abs 
\end{abstract}


\section{Introduction}\label{intro}

\abs
It is a major open problem in representation theory of finite groups
to understand the irreducible modular representations of the symmetric
group $\cS_n$ on $n$ letters. Although considerable progress has been made
in recent decades, not even their dimensions are in general known today.
The situation becomes somewhat better if we restrict ourselves to the
representations whose dimensions are bounded above by a polynomial in $n$:
The aim of the present paper is to obtain an explicit classification 
of all the irreducible modular representations of $\cS_n$ of degree
at most $n^3$, together with explicit degree formulae, the final result
being described in \ref{strategy}, and explicit lists in characteristic
$p\leq 7$ being given in the subsequent sections from \ref{low2} on.

\abs
Our original motivation to investigate into this was a particular question
asked a while ago by E.~O'Brien, in view of the, at that time forthcoming,
paper \cite{FawBreSax}. Accordingly, we are proud to be able to say that
the present results are now used significantly in \cite[Sect.4]{FawBreSax}. 
From a broader perspective, the present paper lies well within the
recent philosophy of collecting results on representations of low degree
of finite simple groups and their close relatives, and aims at contributing
to this programme;
see for example \cite{HisMal,Lub}, where both viewpoints of either
fixing a degree bound, or, for the case of groups of Lie type, allowing
for a polynomial degree bound in terms of the Lie rank, are pursued. 

\abs
In order to achieve the goal specified above, it has turned out that
quite a few additional pieces of information have to be collected or
newly derived, in particular by improving an idea, originally due to
G.~James, to get hands on dimension bounds, by building on the current
knowledge about decomposition numbers of symmetric groups and their
associated Iwahori-Hecke algebras, and by applying a mixture of theoretical
reasoning and computational techniques. 

\absb
The starting point of our considerations is an observation due to G.~James
\cite{JamMinimal}, describing the growth behavior of irreducible
modular representations of the symmetric group $\cS_n$ on $n$ letters,
when $n$ tends to infinity. More precisely:

\abs  
Given a rational prime $p$, let $d^\mu:=\dim_{\F_p}(D^\mu)$
be the dimension of the irreducible $\cS_n$-module $D^\mu$
parameterized by the $p$-regular partition $\mu=[\mu_1,\mu_2,\ldots]$ of $n$.
Fixing a non-negative integer $m$, and assuming that the largest part of $\mu$ 
equals $\mu_1=n-m$, it is shown in \cite[Thm.1]{JamMinimal} that
$d^\mu\sim\frac{n^m}{n!}\cdot d^{\ov{\mu}}$, for $n>\!\!>0$,
where $\ov{\mu}=[\mu_2,\mu_3,\ldots]$, a partition of $n-m$.
In particular, by \cite[Cor.2]{JamMinimal},
$d^\mu$ always is bounded above by $d^\mu\leq n^m$, while there only 
is an asymptotical lower bound $n^{m-1}<d^\mu$ for $n>\!\!>0$. 
Since for small $n$ there are cases such that 
the latter inequality fails, the question arises whether there is
an explicit bound $n_0$ such that $n^{m-1}<d^\mu$ for all $n\geq n_0$.

\abs
A general strategy to find effective lower bound 
functions $f(n)$, fulfilling $f(n)\leq d^\mu$ for all $n\geq n_0$
and some explicitly given $n_0$, is already introduced and used
in the proofs of \cite[Thm.5]{JamMinimal} and its key Lemma 
\cite[La.4]{JamMinimal}.
The major ingredients to these proofs are the branching rule for
irreducible ordinary representations of $\cS_n$, see for example 
\cite[Thm.9.2]{JamLN}, the hook length formula \cite{FraRobThr} 
for the degree of irreducible ordinary representations,
see also \cite[Thm.20.1]{JamLN},
and the unitriangularity of the decomposition matrix of $\cS_n$.
As main theoretical results of the present paper we are going to 
derive two improvements of \cite[La.4]{JamMinimal}:

\absb
The first main result, being valid for arbitrary rational primes $p$, 
is given in \ref{jamesthm}, whose improvement consists of weaker conditions
imposed on candidate lower bound functions $f(n)$;
it is detailed in \ref{jamesrem} how this compares to \cite[La.4]{JamMinimal}.
In practice, it turns out that the weaker checks necessary save
quite a bit of explicit computation, and lends itself to the
treatment independent of $p$ carried out in \ref{degbnd}. 
The new ingredient to the proof is to replace the ordinary branching rule
by the modular branching rule, which is recalled in \ref{modbranch};
of course, the latter had not been available at the time of writing
of \cite{JamMinimal}.

\abs
Although the above theorem also holds for $p=2$,
in this case it would not be strong enough for our purposes, 
the reason being the existence of `very small' representations escaping
the desired growth behavior in low degrees.
The representations to be treated separately are the first and second
basic spin representations, whose branching behavior is described 
\ref{basicspin}, where we also give degree formulae.
This leads to our second main result, given in \ref{james2thm},
which holds only for $p=2$ but takes care of these
exceptional cases. The strategy of proof is reminiscent of 
\cite[Thm.4.3]{Dan}, but again we get away with weaker conditions
to be imposed on candidate lower bound functions $f$,
as is detailed in \ref{danzrem}, which again 
in practice saves quite a bit of explicit computation.

\abs
Actually, we manage to prove both \ref{jamesthm} and \ref{james2thm} 
in a conceptual and unified manner.
In particular, for our major case of interest $m=3$, 
in the cases $p=3$ and $p=2$ the improvements lead to shorter 
proofs of \cite[Prop.3.1]{Dan} and \cite[Prop.4.4]{Dan},
together with smaller lower bounds $n_0$, where regular
behavior sets in. Moreover, in the proofs of \ref{jamesthm} and 
\ref{james2thm} Jantzen-Seitz partitions \cite{JanSei} feature prominently,
where an answer to the following problem would yield a small further
improvement to  \ref{jamesthm}, as is detailed in \ref{seitzrem}: 

\abs
{\bf Open problem:} \it 
Classify the Jantzen-Seitz partitions whose largest part occurs at
least twice; in particular show that there are none for $p=3$.

\absb
The above discussion leads to the question of how to find candidate 
lower degree functions $f(n)$ in the first place. In our main case
of interest $m=3$, due to G.~James's observations \cite{JamMinimal}, 
it seems reasonable to take as a candidate function $f(n)$ a lower bound
function for the degrees $d^\mu$, where $\mu$ runs through all $p$-regular
partitions of $n$ with largest part $\mu_1=n-(m+1)=n-4$. 

\abs
Best possible polynomial lower bound functions for the degrees $d^\mu$, 
where $\mu$ runs through all $p$-regular partitions of $n$ with
largest part at least $n-4$, are collected in \ref{degformulae}. 
For the cases $p\geq 5$ these follow immediately from the degree
formulae given in \cite[La.1.21]{BruKle}, while the cases $p=3$ and $p=2$, 
which are of most interest for us, are spared there.
Actually, lower bound functions, for all $p\geq 2$, can also be derived 
from a close inspection of the results given without proof in 
\cite[App.]{JamMinimal}, but the degree formulae given there are
far from being complete.

\abs
Since, apart from the absence of explicit proofs, for our 
ultimate classification aim we need degree formulae anyway, 
we are deriving these for the missing cases $p=3$ and $p=2$,
and all the partitions $\mu$ mentioned above, in \ref{deg3formulae} and
\ref{deg2formulae}, respectively. As it turns out, to collect complete and 
precise data, but not just estimates, it is necessary to determine most
(if not all) of the rows of the decomposition matrix of $\cS_n$ being
parameterized by partitions $\lb$ of $n$ with largest part $\lb_1\geq n-4$.
To our knowledge, sufficiently far reaching complete results of this 
generic kind, where the partitions considered are restricted to
interesting classes, but $n$ is treated as a parameter, have not been 
determined so far. Prototypical computations in this direction 
are done in \cite{JamWil}, from which we in particular use the
observation \cite[Prop.2.1]{JamWil}, being recalled in \ref{jamestrick},
several times, in order to find upper bounds for 
certain decomposition numbers.

\abs
To make our approach to finding the relevant decomposition numbers of $\cS_n$
as conceptual as possible, we first determine the corresponding part of
the crystallized decomposition matrix of the generic
Iwahori-Hecke algebra $\cH_n(u)$ associated with $\cS_n$, 
with respect to specializing the indeterminate parameter $u$ to
a $p$-th complex root of unity; the necessary background on the relationship
between these decomposition matrices is recalled in \ref{hecke}.
Crystallized decomposition matrices can be computed recursively by 
applying the LLT algorithm \cite{LasLecThi}, which is recalled
in \ref{llt}. We use a modification of it, as detailed in \ref{llttrunc}:

\abs
Let the Fock space $\cF$ be the free abelian group having a basis
consisting of all partitions of all non-negative integers $n$, 
and let $\cF^{>m}$ be the subgroup spanned by the partitions whose 
largest part is smaller than $n-m$, where $m$ is a fixed non-negative integer.
Then the truncated Fock space $\cF(\leq m):=\cF/\cF^{>m}$ has a basis 
consisting of the partitions with largest part at least $n-m$, and thus 
captures precisely the part of the decomposition matrix we are interested in. 
As it turns out, instead of running the LLT algorithm in $\cF$
and subsequently projecting onto $\cF(\leq m)$, a modified version of it
can just be applied directly to $\cF(\leq m)$.
Moreover, $\cF(\leq m)$ has components being parameterized by $n$, 
which for $n$ large enough are all naturally parameterized by the
partitions of the numbers $1\ld m$. Hence it is conceivable that 
computations can be performed generically, with $n$ as a parameter,
in terms of the latter partitions.

\abs
Indeed, this is successfully carried out in \ref{crystdec3} and
\ref{crystdec2} for $m=4$, and the cases $p=3$ and $p=2$, respectively.
Then, in \ref{decmat3} and 
\ref{decmat2}, respectively, the part of the crystallized decomposition
matrix found is used, together with {\it ad hoc} arguments from the 
toolbox of modular representation theory of symmetric groups,
to determine the corresponding part of the decomposition matrix of $\cS_n$,
Finally, in \ref{deg3formulae} and \ref{deg2formulae}, respectively,
degree formulae are collected.
Moreover, the observations made for the cases treated here lead us
to ask the following question, which is discussed in some more
detail in \ref{lltrem}:

\abs
{\bf Open problem:} \it
Show that the crystallized decomposition matrix corresponding to the $n$-th
component of $\cF(\leq m)$ depends only on the congruence class of 
$n$ modulo $p$, whenever $n\geq 2m+1$.

\absb 
Having this theoretical machinery in place, we are finally prepared to 
tackle our original problem of classifying the $p$-regular partitions
$\mu$ of $n$ such that $d^\mu\leq n^3$:
As already indicated above, as a test function to employ 
\ref{jamesthm} and \ref{james2thm}, we essentially take the lower 
bound function for the degrees $d^\mu$ where $\mu$ has largest part
$n-4$ given in \ref{degformulae}. Actually, as detailed in \ref{degbnd},
in order to obtain a description independent of $p$, and a strong bound
for $n_0$, in low degrees we are using a suitable modification of the
natural lower bound function.

\abs
This covers the generic behavior, and leaves finitely many $n$
to be considered separately. The computational strategy used to
do so, as automatic as possible for fixed $p$, is described in
\ref{strategy}. We are employing the computer algebra system
\GAP{} \cite{GAP} as our major computational environment. Here,
we also use explicit data on decomposition numbers of symmetric
groups accessible in the character table database 
{\sf CTblLib} \cite{CTblLib} of \GAP{}, and on the web page 
\cite{ModularAtlasProject}
of the {\sf ModularAtlas} project \cite{ModularAtlas}.
An implementation of the LLT algorithm is
provided in the {\sf SPECHT} package \cite{SPECHT}, which is 
available through the computer algebra system {\sf CHEVIE} \cite{CHEVIE}.
The explicit classification thus obtained for the cases $p=3$ and 
$p=2$ is given in \ref{low3} and \ref{low2}, respectively.

\abs\abs
{\bf Acknowledgement:}
The author is grateful for financial support by the 
German Science Foundation (DFG) Scientific Priority Programmes
SPP-1388 `Representation Theory' and SPP-1478 `Algorithmic and Experimental
Methods in Algebra, Geometry, and Number Theory', to which this paper
is a contribution. Moreover, he thanks the anonymous referee for their
comments, which helped to improve the exposition.

\section{Background}\label{back}

\abs
We assume the reader familiar with ordinary and modular representation
theory of symmetric groups, as well as with the representation theory
of Iwahori-Hecke algebras; as general references see \cite{JamLN,JamKer}
and \cite{Mat}, respectively. We recall the necessary facts, where
in particular we take this opportunity to fix notation.

\AbsT{Partitions.}\label{partitions}
For $n\in\N_0$ let $\cP_n$ denote the set of all partitions 
$\lb=[\lb_1,\lb_2\ld\lb_l]$ of $n$, where we assume 
$\lb_1\geq\cdots\geq\lb_l>0$ and $l=l(\lb)\in\N_0$ is the length of $\lb$;
for convenience we let $\lb_{l+1}:=0$.
Moreover, given $\lb\in\cP_n$, whenever $n\geq 1$ let 
$\ov{\lb}:=[\lb_2,\lb_3\ld\lb_l]\in\cP_{n-\lb_1}$
be the partition obtained from $\lb$ by removing its first part;
for $n=0$ we let $\ov{[]}=[]\in\cP_0$. Since a partition $\lb$
is uniquely determined, as soon as $n$ is understood, by its tail part 
$\ov{\lb}$, we may also write it as $(\ov{\lb})$;
here we are using round brackets to distinguish this 
from our standard square bracket notation for partitions.
This convention will prove convenient in view of the following
key combinatorial definition of this paper:

\abs
For $m\in\N_0$ let 
$$ \cP_n(m):=\{\lb\in\cP_n;\lb_1=n-m\}=\{\lb\in\cP_n;\ov{\lb}\in\cP_m\}
   \sseq\cP_n$$
and $\cP_n(\leq m):=\coprod_{j=0}^m \cP_n(j)\sseq\cP_n$.
Hence we get an injective map $\cP_n(m)\ra\cP_m\cn\lb\mt\ov{\lb}$,
where $\ov{\lb}\in\cP_m$ actually is a tail part of some partition in $\cP_n$
if and only if $n-m\geq\ov{\lb}_1$. In particular,
we have $\cP_n(m)=\emp$ whenever $n<m$ or $n=m>0$,
and the map $\cP_n(m)\ra\cP_m$ is surjective if and only if $n\geq 2m$.

\abs
In the sequel, a particular role will be played by the sets $\cP_n(m)$ 
for $m\leq 4$, the round bracket notation introduced above yielding
a compact description independent of $n$: we have 
$\cP_n(0)=\{()\}$ for $n\geq 0$,
$\cP_n(1)=\{(1)\}$ for $n\geq 2$, 
$\cP_n(2)=\{(2),(1^2)\}$ for $n\geq 4$,
$\cP_n(3)=\{(3),(2,1),(1^3)\}$ for $n\geq 6$, and
$\cP_n(4)=\{(4),(3,1),(2^2),(2,1^2),(1^4)\}$ for $n\geq 8$.

\absb
Given a rational prime $p$, let $\cP_n^{p\text{-reg}}\sseq\cP_n$ 
be the set of all $p$-regular partitions, that is those for
which any part occurs less than $p$-times, and let
$\cP_n^{p\text{-reg}}(m):=\cP_n^{p\text{-reg}}\cap\cP_n(m)\sseq\cP_n$ and 
$\cP_n^{p\text{-reg}}(\leq m):=\cP_n^{p\text{-reg}}\cap\cP_n(\leq m)
 \sseq\cP_n$. 
As usual, we assume both $\cP_n$ and $\cP_n^{p\text{-reg}}$ to be ordered
reversed lexicographically. Then the sets $\cP_n(\leq m)$ and 
$\cP_n^{p\text{-reg}}(\leq m)$ are order ideals in $\cP_n$ and 
$\cP_n^{p\text{-reg}}$, respectively.

\abs
Moreover, let $\cP_n\ra\cP_n^{p\text{-reg}}\cn\lb\mt\lb^R$
be the $p$-regularisation map, see \cite[6.3.48]{JamKer};
of course, this map restricts to the identity map on $\cP_n^{p\text{-reg}}$.
Finally, let $\cP_n^{p\text{-reg}}\ra\cP_n^{p\text{-reg}}\cn\mu\mt\mu^M$
be the Mullineux map, see \cite{BesOls,Mul};
recall that this is the identity map in the case $p=2$,
and by \cite{ForKle} records the effect of tensoring 
irreducible $p$-modular representations with the sign representation.

\AbsT{Young diagrams.}
As usual, we identify a partition $\lb=[\lb_1\ld\lb_l]\in\cP_n$ 
with its Young diagram, that is an array of $l$ rows of `nodes',
where row $a$, counting from top to bottom, contains $\lb_a$ nodes,
counting from left to right. The node being in row $a$ and column $b$ 
of a Young diagram is called its $(a,b)$-node. Given a rational prime $p$, 
the $(a,b)$-node is called an $i$-node, or is said to have
$p$-residue $i\in\{0\ld p-1\}$, if $i\equiv b-a\pmod{p}$.

\abs
The $(a,\lb_a)$-node $x$, for some $a\in\{1\ld l\}$, is called removable
if $\lb_a>\lb_{a+1}$; in this case, we let $\lb\smin\{x\}\in\cP_{n-1}$
be the partition obtained by removing $x$. 
Similarly, the $(a,\lb_a+1)$-node $y$, for some $a\in\{1\ld l+1\}$, 
is called addable if either $a=1$ or $\lb_a<\lb_{a-1}$;
in this case, we let $\lb\cup\{y\}\in\cP_{n+1}$
be the partition obtained by adding $y$.
Let $R_i(\lb)$ and $A_i(\lb)$ be the sets of removable and addable 
$i$-nodes of $\lb$, respectively.  
If $x\in R_i(\lb)$ or $x\in A_i(\lb)$, then let $r_i(\lb,x)$
and $a_i(\lb,x)$ be the number of removable and addable $i$-nodes 
strictly to the right of $x$, respectively.

\abs
Let $\mu\in\cP_n^{p\text{-reg}}$, and let $i\in\{0\ld p-1\}$. 
Labelling each removable and addable $i$-node 
of $\mu$ by `$-$' and `$+$', respectively, and reading off these labels
from left to right yields a sequence called the $i$-signature of $\mu$,
and repeatedly cancelling `$-+$' pairs in the sequence we end up with
the reduced $i$-signature of $\mu$.
The removable $i$-nodes surviving in the reduced $i$-signature 
are called $i$-normal, and if there is an $i$-normal node
then the leftmost one is called $i$-good.
Let $N_i(\mu)$ be the set of $i$-normal nodes of $\mu$;
note that $1+r_i(\mu,x)-a_i(\mu,x)$ is the number of $i$-normal nodes
to the right of $x\in N_i(\mu)$, including $x$. 
Similarly, the addable $i$-nodes surviving in the reduced $i$-signature
are called $i$-conormal, and if there is an $i$-conormal node
then the rightmost one is called $i$-cogood.

\AbsT{Grothendieck groups.}
For $n\in\N_0$ let $G(\modcat{}{K\cS_n})$, where $K$ is a field, be 
the Grothendieck group of the category of finitely generated $K\cS_n$-modules,
the equivalence class of a module $M$ being denoted by $[M]$.

\abs
$G(\modcat{}{\Q\cS_n})$ has the classes of the
(absolutely irreducible) Specht modules $\{[S^\lb];\lb\in\cP_n\}$
as its standard $\Z$-basis, and thus can be naturally identified
with the free abelian group $\Z[\cP_n]$.
Similarly, $G(\modcat{}{\F_p\cS_n})$, where $p$ is a rational prime,
has the classes of the (absolutely) irreducible modules
$\{[D^\lb];\lb\in\cP_n^{p\text{-reg}}\}$ as its standard $\Z$-basis,
and thus can be identified with the free abelian group
$\Z[\cP_n^{p\text{-reg}}]$.
As it turns out, this identification is natural as well:

\abs
Associated to this setting there is the decomposition map
$$ \Dt_n^p\cn\Z[\cP_n]\ra\Z[\cP_n^{p\text{-reg}}]\cn
   \lb\mt\sum_{\mu\in\cP_n^{p\text{-reg}}}[S^\lb\cn D^\mu]\cdot\mu ,$$
where the decomposition number $[S^\lb\cn D^\mu]\in\N_0$ is the 
multiplicity of the constituent $D^\mu$ in a composition series
of a $p$-modular reduction of $S^\lb$.
The map $\Dt_n^p$ has the following `unitriangularity' properties:
We have $[S^\lb\cn D^\mu]>0$ only if $\lb\unlhd\mu$, 
where `$\unlhd$' denotes the dominance partial order of $\cP_n$,
and $[S^\mu\cn D^\mu]=1$ for all $\mu\in\cP_n^{p\text{-reg}}$.
Hence this gives rise to maps 
$\Dt_n^p(\leq m)\cn\Z[\cP_n(\leq m)]\ra\Z[\cP_n^{p\text{-reg}}(\leq m)]$.
The matrices $\cD_n$ and $\cD_n(\leq m)$ of $\Dt_n^p$ and $\Dt_n^p(\leq m)$,
respectively, with respect to reversed lexicographical ordering of 
partitions, are called the associated decomposition matrices.

\abs
Moreover, for $n\geq 1$, the ordinary branching rule for Specht
modules, see for example \cite[Thm.9.2]{JamLN},
gives rise to $i$-restriction maps, where $i\in\{0\ld p-1\}$,
$$ \da_{i}\cn\Z[\cP_n]\ra\Z[\cP_{n-1}]\cn
   \lb\mt\sum_{x\in R_i(\lb)}\lb\smin\{x\} .$$
(Similarly, there are $i$-induction maps, but these we will not need.)
By way of $p$-modular reduction this induces maps
$$ \da_{i}\cn\Z[\cP_n^{p\text{-reg}}]\ra\Z[\cP_{n-1}^{p\text{-reg}}]\cn
   \mu\mt\sum_{\tau\in\cP_{n-1}^{p\text{-reg}}}
   [D^\mu\da_i\cn D^\tau]\cdot\tau, $$
where $[D^\mu\da_i\cn D^\tau]\in\N_0$ is the multiplicity of the 
constituent $D^\tau$ in a composition series of $D^\mu\da_i$.
These maps are not at all well-understood, but are partly described by 
the following theorem, where we only state the facts needed later:

\ThmT{Modular branching rule; \cite{Kle,KleII}, 
      see also \cite{BruKle2,KleBook}}\label{modbranch}\mbox{}\\ 
Let $\mu\in\cP_n^{p\text{-reg}}$, where $n\geq 1$, 
and let $i\in\{0\ld p-1\}$.

\absb\it
We have $D^\mu\da_{i}=\{0\}$ if and only if $N_i(\mu)=\emp$.
Moreover, $D^\mu\da_{i}$ is irreducible if and only if
$N_i(\mu)$ is a singleton set.

\absb\it
Let $x\in N_i(\mu)$ be an $i$-normal node such that 
$\mu\smin\{x\}\in\cP_{n-1}^{p\text{-reg}}$;
note that this in particular holds if $x$ is the $i$-good node. Then we have  
$$ [D^\mu\da_{i}\cn D^{\mu\smin\{x\}}]= 1+r_i(\mu,x)-a_i(\mu,x) .$$

\vspace*{-2.2em}\QED

\abs\abs\abs
We recall a straightforward observation, relating decomposition numbers
of $\cS_n$ and $\cS_{n-1}$ with constituent multiplicities of restrictions 
of irreducible $\cS_n$-modules to $\cS_{n-1}$. This actually proves to be
powerful to find upper bounds for certain decomposition numbers by induction:

\PropT{\cite[Prop.2.1]{JamWil}}\label{jamestrick}\mbox{}\\
Let $\lb\in\cP_n$, where $n\geq 1$, and $\tau\in\cP_{n-1}^{p\text{-reg}}$.
Then for $i\in\{0\ld p-1\}$ we have
$$ \sum_{x\in R_i(\lb)}[S^{\lb\smin\{x\}}\cn D^\tau]
  =\sum_{\mu\in\cP_n^{p\text{-reg}}}[S^\lb\cn D^\mu]
   \cdot [D^\mu\da_i\cn D^\tau] .$$
In particular, for any $\mu\in\cP_n^{p\text{-reg}}$ we have 
$$ [S^\lb\cn D^\mu]\cdot[D^\mu\da_i\cn D^\tau] 
   \leq\sum_{x\in R_i(\lb)}[S^{\lb\smin\{x\}}\cn D^\tau], $$
where equality holds if and only if $D^\mu$ is the only constituent
of $S^\lb$ such that $[D^\mu\da_i\cn D^\tau]>0$.
\QED

\AbsT{Iwahori-Hecke algebras.}\label{hecke}
In order to find lower bounds for decomposition numbers of $\cS_n$
we make use of the generic Iwahori-Hecke algebra $\cH_n(u)$ associated 
with $\cS_n$. Here, the parameter $u$ is an indeterminate, and $\cH_n(u)$ 
is an algebra over the field $\Q(u)$.
By Tits' Deformation Theorem, the (split semi-simple) module categories 
$\modcat{}{\Q(u)\cS_n}$ and $\modcat{}{\cH_n(u)}$ can be identified,
hence their Grothendieck groups can be identified as well.
More precisely, there again is a natural choice of Specht modules, 
see \cite{DipJam}, which we again denote by $S^\lb$, for $\lb\in\cP_n$.
Next, we recall the necessary facts about the decomposition theory of
$\cH_n(u)$, for more details see \cite[Ch.6.2]{Mat}:

\abs
Specializing the parameter $u$ to become a primitive $p$-th root of unity
$\zt:=\zt_p\in\C$ yields the $\Q(\zt)$-algebra $\cH_n(\zt)$. 
Then $G(\modcat{}{\cH_n(\zt)})$ has the classes of the (absolutely)
irreducible modules $\{[D_\zt^\lb];\lb\in\cP_n^{p\text{-reg}}\}$ as
its standard $\Z$-basis, and can also be identified naturally with
the free abelian group $\Z[\cP_n^{p\text{-reg}}]$. Again 
there is an associated decomposition map
$$ \Dt_n^\zt\cn\Z[\cP_n]\ra\Z[\cP_n^{p\text{-reg}}]\cn
   \lb\mt\sum_{\mu\in\cP_n^{p\text{-reg}}}[S^\lb\cn D_\zt^\mu]\cdot\mu ,$$
where the decomposition number $[S^\lb\cn D_\zt^\mu]\in\N_0$ is the
multiplicity of the constituent $D_\zt^\mu$ in a composition series
of a $\zt$-modular reduction of $S^\lb$. The map $\Dt_n^\zt$
enjoys the same unitriangularity properties as $\Dt_n^p$, and thus similarly 
gives rise to decomposition matrices $\cD_n^\zt$ and $\cD_n^\zt(\leq m)$.

\absb
Considering $p$-modular reduction again, which specializes the 
$p$-th root of unity $\zt\in\C$ to $1\in\F_p$, it turns out, 
see \cite{GeckTrees}, that there is a decomposition map
$$ A_n^\zt\cn\Z[\cP_n^{p\text{-reg}}]\ra\Z[\cP_n^{p\text{-reg}}]\cn
   \tau\mt\sum_{\mu\in\cP_n^{p\text{-reg}}}[D_\zt^\tau\cn D^\mu]\cdot\mu ,$$
where the decomposition number $[D_\zt^\tau\cn D^\mu]\in\N_0$ is the
multiplicity of the constituent $D^\mu$ in a composition series
of a $p$-modular reduction of $D_\zt^\tau$.
Then we have $\Dt_n^p=\Dt_n^\zt\cdot A_n^\zt$, which implies that 
$A_n^\zt$ also has the celebrated unitriangularity properties, 
giving rise to maps 
$A_n^\zt(\leq m)\cn\Z[\cP_n^{p\text{-reg}}(\leq m)]
                \ra\Z[\cP_n^{p\text{-reg}}(\leq m)]$.
The matrices $\cA_n$ and $\cA_n(\leq m)$ associated with $A_n^\zt$
and $A_n^\zt(\leq m)$, respectively, being square lower unitriangular
matrices, are called adjustment matrices.

\abs
In particular, from this we indeed get lower bounds for 
decomposition numbers: we have $[S^\lb\cn D^\mu]\geq [S^\lb\cn D_\zt^\mu]$,
for all $\lb\in\cP_n$ and $\mu\in\cP_n^{p\text{-reg}}$;
or phrased in terms of dimensions we get 
$\dim_{\F_p}(D^\mu)\leq\dim_{\Q(\zt)}(D_\zt^\mu)$.
In the cases of small $m$ we are interested in, it eventually turns
out that the adjustment matrices $\cA_n(\leq m)$ actually are sparse 
with very small entries, so that the above inequalities tend to provide 
good lower bounds. Hence to make this approach effective, we need to 
get hands on the matrices $\cD_n^\zt$. This is actually achieved as follows:  

\AbsT{The LLT algorithm.}\label{llt}
We proceed towards a recursive description of $\cD_n^\zt$, where 
we only describe the facts necessary; for more details see \cite[Ch.6.1]{Mat}:

\abs
We consider the Fock space $\cF:=\bigoplus_{n\in\N_0}\Z[\cP_n]$.
Let $\cF_q:=\Z[q,q^{-1}]\otm_\Z\cF$, where $q$ is an indeterminate.
For $i\in\{0\ld p-1\}$ running through the residue classes modulo $p$
and $k\in\N_0$, in view of \cite[La.6.16]{Mat} we define 
$\Z[q,q^{-1}]$-linear maps by
$$ F_i^{(k)}\cn\cF_q\ra\cF_q\cn\lb\mt\sum_{X\sseq A_i(\lb);|X|=k}
   q^{\sum_{x\in X}a_i(\lb\cup X,x)-r_i(\lb,x)}\cdot(\lb\cup X) ,$$
where we let $F_i:=F_i^{(1)}$; note that $F_i^{(0)}=\id_{\cF_q}$.
Let $\cU$ be the $\Z[q,q^{-1}]$-algebra generated by 
the divided power operators $\{F_i^{(k)};i\in\{0\ld p-1\};k\in\N_0\}$,
and let $\baraut$ be the $\Z$-linear ring involution of $\cU$
defined by $\ov{q}:=q^{-1}$ and $\ov{F_i^{(k)}}:=F_i^{(k)}$.
Moreover, let $\cL\leq\cF_q$ be the $\cU$-submodule of
$\cF_q$ generated by the empty partition $[]\in\cP_0$, and
let $\baraut$ be the $\Z$-linear involution on $\cL$ defined by
$\ov{[]\cdot u}:=[]\cdot\ov{u}$, for all $u\in\cU$.
Then, from \cite{Kas}, there exists a unique $\Z[q,q^{-1}]$-basis 
$\{B^\mu\in\cF_q;\mu\in\coprod_{n\in\N_0}\cP_n^{p\text{-reg}}\}$ 
of $\cL$, being called its global crystal basis, such that 
$B^\mu-\mu\in q\Z[q]\otm_\Z\cF\sseq\cF_q$ and $\ov{B^\mu}=B^\mu$.

\abs
Actually we have $B^\mu\in\Z[q,q^{-1}]\otm_\Z\Z[\cP_n]$ whenever
$\mu\in\cP_n^{p\text{-reg}}$. Writing 
$B^\mu=\sum_{\lb\in\cP_n}b_{\lb,\mu}(q)\cdot\lb$ in terms of
the standard basis $\cP_n$, where $b_{\lb,\mu}(q)\in\Z[q,q^{-1}]$ 
such that $b_{\mu,\mu}(q)=1$,
gives rise to the crystallized decomposition matrices $\cD_n^q$.
These again, by \cite{LasLecThi}, share the celebrated unitriangularity
properties, and by \cite{VarVas} fulfill the positivity condition
$b_{\lb,\mu}(q)\in q\N_0[q]$ whenever $\lb\lhd\mu$. In particular,
restricting to the sets $\cP_n(\leq m)$ and $\cP_n^{p\text{-reg}}(\leq m)$,
we also get crystallized decomposition matrices $\cD_n^q(\leq m)$.
Moreover, as was conjectured in \cite{LasLecThi} and
proved in \cite{Ari}, the desired connection to the decomposition numbers 
$[S^\lb\cn D_\zt^\mu]$ is simply given by 
$[S^\lb\cn D_\zt^\mu]=b_{\lb,\mu}(1)$, in other words
the matrices $\cD^\zt_n$ are obtained by evaluating the crystallized 
decomposition matrices $\cD^q_n$ at $q=1$.

\absb
The elements $B^\mu\in\cF_q$, for $\mu\in\cP_n^{p\text{-reg}}$,
are found recursively by the following algorithm \cite{LasLecThi}: 
For the empty partition $[]\in\cP_0$ we have $B^{[]}=[]\in\cF_q$, 
hence we may assume that $n\geq 1$. We consider the
sequence of $p$-ladders, in the sense of \cite[6.3.45]{JamKer},
running from left to right, and starting with the $0$-th ladder passing
through the $(0,0)$-node; then the $j$-th ladder, where $j\in\N_0$,
consists of nodes of residue $i\equiv p-j\pmod{p}$.
Given $\mu\in\cP_n^{p\text{-reg}}$, 
recording the number of nodes of $\mu$ lying on the various $p$-ladders, 
yields a sequence $[k_0\ld k_r]$, where $r\in\N$ and $k_i\in\N_0$
such that $k_r\geq 1$ and $\sum_{i=1}^r k_i=n$. Then let
$$ F^\mu:=F_0^{(k_0)}F_{p-1}^{(k_1)}F_{p-2}^{(k_2)}
          \cdots F_{p-r}^{(k_r)}\in\cU ,$$
be the associated product of divided power operators,
where subscripts are read modulo $p$. Applying it,
from left to right, to the empty partition $[]\in\cP_0$ yields
$A^\mu:=[]\cdot F^\mu\in\cF_q$.
Using the condition $B^\mu-\mu\in q\Z[q]\otm_\Z\cF$, it turns out that
$A^\mu=B^\mu+\sum_{\tau\in\cP_n^{p\text{-reg}};\tau\lhd\mu}
 a_{\tau,\mu}(q)\cdot B^\tau\in\cF_q$,
where $a_{\tau,\mu}(q)\in\Z[q,q^{-1}]$.
Thus we have $A^\mu=B^\mu$ whenever $\mu$ is smallest in
$\cP_n^{p\text{-reg}}$ with respect to the dominance partial order,
and otherwise the $a_{\tau,\mu}(q)$ are uniquely determined recursively
using the elements $B^\tau$ for $\tau\in\cP_n^{p\text{-reg}}$ such that
$\tau\lhd\mu$, and the invariance property $\ov{B^\mu}=B^\mu$.

\absb
In the sequel, we will also need the following alternative description, 
taken from \cite[Sect.6.25]{Mat}, saying that we may also use
induction on $n\in\N_0$ as follows: 
Since $B^{[]}=[]\in\cF_q$ anyway, we may assume that $n\geq 1$. 
Let $X$ be the intersection of the $r$-th $p$-ladder with $\mu$,
that is the rightmost one intersecting $\mu$ non-trivially,
hence we have $\mu\smin X\in\cP_{n-k_r}^{p\text{-reg}}$. By induction 
we may let $A^{\prime\mu}:=B^{\mu\smin X}\cdot F_i^{(k)}\in\cF_q$. 
Then we have 
$A^{\prime\mu}=B^\mu+\sum_{\tau\in\cP_n^{p\text{-reg}};\tau\lhd\mu}
 a'_{\tau,\mu}(q)\cdot B^\tau\in\cF_q$,
where $a'_{\tau,\mu}(q)\in\Z[q,q^{-1}]$, and we may again proceed as above.

\AbsT{The truncated Fock space.}\label{llttrunc}
Let $\cF^{>m}:=\bigoplus_{n\in\N_0}\Z[\cP_n\smin\cP_n(\leq m)]$,
where $m\in\N_0$. Hence we have
$\cF(\leq m):=\cF/\cF^{>m}\cong\bigoplus_{n\in\N_0}\Z[\cP_n(\leq m)]$.
Then it follows directly from the definition of the maps
$F_i^{(k)}\cn\cF_q\ra\cF_q$, that 
$\cF^{>m}_q:=\Z[q,q^{-1}]\otm_\Z\cF^{>m}\leq\cF_q$
is a $\cU$-submodule, thus
$\cF_q(\leq m):=\cF_q/\cF^{>m}_q\cong \Z[q,q^{-1}]\otm_\Z(\cF/\cF^{>m})$
is a $\cU$-module as well.
Moreover, $\cL\cap\cF^{>m}\sseq\cL$ is invariant with respect to $\baraut$,
hence we get an induced $\Z$-linear involution on 
$(\cL+\cF_q^{>m})/\cF^{>m}_q\leq\cF_q(\leq m)$. 
We are going to consider the natural projection 
$\{B^\mu+\cF^{>m}_q\in\cF_q(\leq m); 
   \mu\in\coprod_{n\in\N_0}\cP_n^{p\text{-reg}}\}$
of the global crystal basis to $\cF_q(\leq m)$; 
this is indicated by writing
$B^\mu\equiv\sum_{\lb\in\cP_n(\leq m)}b_{\lb,\mu}(q)\cdot\lb$.

\abs
It follows from the description in \ref{llt},
and the fact that for all $n\in\N_0$ the set $\cP_n(\leq m)$
is an ideal of $\cP_n$ with respect to the reversed dominance partial order,
that instead of running the LLT algorithm in $\cF_q$ and projecting
the global crystal basis to $\cF_q(\leq m)$ afterwards, we may just
run the LLT algorithm in the quotient space $\cF_q(\leq m)$ right
from the beginning, and still end up with the elements
$B^\mu+\cF^{>m}_q\in\cF_q(\leq m)$.
In other words, in terms of the standard basis $\cP_n(\leq m)$,
this directly computes the matrices $\cD_n^q(\leq m)$.

\absb
Let $n\geq 2m+1$ and $\mu\in\cP_n^{p\text{-reg}}(m)$; hence 
we have $\mu_1=n-m>m\geq\mu_2$.
Assume that the rightmost $p$-ladder intersecting $\mu$ has
intersection of cardinality at least $2$ with $\mu$. Then, due
to $p$-regularity, this ladder intersects $\mu$ in its (rightmost)
$(1,n-m)$-node $x$ and in the $(p,n-m-1)$-node.
Hence we have $m\geq (n-m-1)(p-1)$, implying
$n\leq\fl{\frac{pm}{p-1}}+1\leq 2m+1$.
We conclude that, whenever $n\geq 2m+2$,
the rightmost $p$-ladder meeting $\mu$ intersects it just in $\{x\}$.

\abs 
Thus, in this case, since $x$ has residue $r\equiv n-m-1\pmod{p}$,
using the notation of \ref{llt}, we get
$A^{\prime\mu}=B^{\mu\smin\{x\}}\cdot F_r\in\cF_q$.
Moreover, since $\mu\smin\{x\}\in \cP_{n-1}^{p\text{-reg}}(m)$ again,
we conclude that $F^\mu\in\cU$ is ultimately periodic of the form
$F^\mu=\cdots(F_{r-p+1}F_{r-p+2}\cdots F_{r-1}F_r)\cdot
      (F_{r-p+1}F_{r-p+2}\cdots F_{r-1}F_r)$,
where the periodic tail of $F^\mu$ only depends on the 
congruence class of $n$ modulo $p$. 

\abs 
Finally, still assuming $n\geq 2m+1$, the action map 
$F_i\cn\cP_n(m)\ra\cP_{n+1}(m)\dcup\cP_{n+1}(m+1)$,
where $i\in\{0\ld p-1\}$, can be described,
using the identification $\cP_n(m)\ra\cP_m\cn\lb\mt\ov{\lb}$
from \ref{partitions}, by a map
$F_i\cn\cP_m\ra\cP_m\dcup\cP_{m+1}$,
only depending on the residue class of $n$ modulo $p$,
where actually we have $F_i\cn\cP_m\ra\cP_{m+1}$ whenever
$i\not\equiv n-m\pmod{p}$.

\Rem\label{lltrem}
Motivated by the above observations on the dependence of the 
combinatorics just on the residue class of $n$ modulo $p$,
we are wondering whether, for $n$ large enough, the outcome of the
LLT algorithm on the truncated Fock space might be described generically 
in terms of $\coprod_{n\in\N_0}\cP_n(\leq m)$, using the 
identification from \ref{partitions} and treating $n$ as a parameter,
possibly depending on the residue class of $n$ modulo $p$.

\abs
This is indeed true in a strong sense for the 
following special case:
Let $n\geq 2m+1$. Then for $\mu\in\cP_n^{p\text{-reg}}(m)$
we have $B^\mu\equiv\sum_{\lb\in\cP_n(m)}b_{\lb,\mu}(q)\cdot\lb$,
where $b_{\lb,\lb}(q)=1$ and $b_{\lb,\mu}(q)\in q\N_0[q]$ for $\lb\lhd\mu$,
that is modulo $\cF_q(>m)$ all partitions occurring in $B^\mu$ 
belong to $\cP_n(m)$. Since the $(1,n-m+1)$-node $x$ belongs to $A_i(\lb)$, 
for all $\lb\in\cP_n(m)$, where $i\equiv n-m\pmod{p}$, we get
$A^{\prime(\mu\cup\{x\})}:=B^\mu\cdot F_i\equiv
\sum_{\lb\in\cP_n(m)}b_{\lb,\mu}(q)\cdot(\lb\cup\{x\})$.
Thus, from the properties of the $b_{\lb,\mu}(q)$, we get
$B^{\mu\cup\{x\}}\equiv A^{\prime(\mu\cup\{x\})}$.  
Hence we infer that
$B^\mu\equiv\sum_{\ov{\lb}\in\cP_m}b_{\ov{\lb},\ov{\mu}}(q)\cdot\ov{\lb}$,
with coefficients only depending on $\ov{\mu}\in\cP_m^{p\text{-reg}}$, 
but being independent of $n$.
\QED

\abs
But, in view of the results in \ref{crystdec3} and \ref{crystdec2}, 
covering the case $m=4$ for $p=3$ and $p=2$, respectively,
it cannot possibly be expected that, in general, for 
$\mu\in\cP_n^{p\text{-reg}}(\leq m)$ the basis element
$B^\mu\equiv\sum_{\lb\in\cP_n(\leq m)}b_{\lb,\mu}(q)\cdot\lb$
only depends on $\ov{\mu}$. 
Still, from these results it is tempting to expect that $B^\mu$ 
modulo $\cF_q(>m)$ only depends on $\ov{\mu}$ and the congruence 
class of $n$ modulo $p$, that is we are wondering whether for $n\geq 2m+1$
and $\mu=[\mu_1,\ov{\mu}]\in\cP_n^{p\text{-reg}}(\leq m)$
we have
$$ B^{[\mu_1+p,\ov{\mu}]}\equiv\sum_{\lb=[\lb_1,\ov{\lb}]\in\cP_n(\leq m)}
   b_{\lb,\mu}(q)\cdot [\lb_1+p,\ov{\lb}] .$$
If this held true, then, viewing it as a statement on the 
crystallized decomposition matrix $\cD_n^q(\leq m)$, this
would entail a similar statement on the decomposition matrix 
$\cD_n^\zeta(\leq m)$ of the Iwahori-Hecke algebra $\cH_n(u)$,
while the analogous statement for the decomposition matrix
$\cD_n(\leq m)$ of the symmetric group cannot possibly hold, 
as for example the results in \ref{decmat3} and \ref{decmat2} show.
Anyway, to our knowledge this has not yet been examined in the 
literature, and we leave it as an open question to the reader.

\section{Decomposition numbers in characteristic $3$}\label{char3}

\abs
In order to get an overview over the irreducible representations
parameterized by $\mu\in\cP_n^{3\text{-reg}}(\leq 4)$, we determine
the crystallized decomposition matrices $\cD_n^q(\leq 4)$ and the 
decomposition matrices $\cD_n(\leq 4)$ for $p=3$.

\AbsT{Crystallized decomposition matrices.}\label{crystdec3}
We apply the LLT algorithm to the truncated Fock space,
according to the description in \ref{llt} and \ref{llttrunc}.

\absb 
The maps $F_i\cn\cF_q(\leq 4)\ra\cF_q(\leq 4)$,
where $i\in\{0,1,2\}$ runs through the residue classes modulo $3$,
are given in Table \ref{foperators3}, which should be read as follows:
We use the round bracket notation introduced in \ref{partitions}.
Given $\ov{\lb}\in\cP_m$, for some $m\in\N_0$, then the action of
$F_i$ on $(\ov{\lb})\in\cP_n(m)$ only depends on the residue class 
of $n$ modulo $3$ as soon as $n-m\geq\ov{\lb}_1+1$; in particular,
this holds for all $\ov{\lb}\in\cP_m$ whenever $n\geq 2m+1$.
For example, for $\ov{\lb}=(2)$, for all $n\geq 5$ such that
$n\equiv 2\pmod{3}$, we have
$F_1\cn
 [n-2,2]\mt [n-2,3]+q\cdot [n-2,2,1]$.

\begin{table}\caption{Action on truncated Fock space.}\label{foperators3}
$$ \begin{array}{|l||ccc|}
\hline
() & n\equiv 0 & n\equiv 1 & n\equiv 2 \\
\hline
\hline
F_0 & () & & \\
\hline
F_1 & & () & \\
\hline
F_2 & & & () \\
    & q^{-1}(1) & (1) & q(1) \\
\hline
\end{array} 
\quad 
\begin{array}{|l||ccc|}
\hline
(1) & n\equiv 0 & n\equiv 1 & n\equiv 2 \\
\hline
\hline
F_0 & & (1) & \\
    & (2) & q(2) & q^{-1}(2) \\
\hline
F_1 & & & (1) \\
    & q^{-1}(1^2) & (1^2) & q(1^2) \\
\hline
F_2 & (1) & & \\
\hline
\end{array} $$
$$ \begin{array}{|l||ccc|}
\hline
(2) & n\equiv 0 & n\equiv 1 & n\equiv 2 \\
\hline
\hline
F_0 & & & (2) \\
\hline
F_1 & (2) & & \\
    & q(3) & q^{-1}(3) & (3) \\
    & q^2(2,1) & (2,1) & q(2,1) \\
\hline
F_2 & & (2) & \\
\hline
\end{array} 
\quad 
\begin{array}{|l||ccc|}
\hline
(1^2) & n\equiv 0 & n\equiv 1 & n\equiv 2 \\
\hline
\hline
F_0 & & & (1^2) \\
    & q^{-1}(2,1) & (2,1) & q(2,1) \\
    & (1^3) & q(1^3) & q^2(1^3) \\
\hline
F_1 & (1^2) & & \\
\hline
F_2 & & (1^2) & \\
\hline
\end{array} $$
$$ \begin{array}{|l||ccc|}
\hline
(3) & n\equiv 0 & n\equiv 1 & n\equiv 2 \\
\hline
\hline
F_0 & (3) & & \\
\hline
F_1 & & (3) & \\
    & q^{-1}(3,1) & (3,1) & q^{-2}(3,1) \\
\hline
F_2 & & & (3) \\
    & q^{-1}(4) & (4) & q(4) \\
\hline
\end{array}
\quad 
\begin{array}{|l||ccc|}
\hline
(2,1)\! & n\equiv 0 & n\equiv 1 & n\equiv 2 \\
\hline
\hline
F_0 & (2,1) & & \\
    & (2,1^2) & q^{-2}(2,1^2) & q^{-1}(2,1^2) \\
\hline
F_1 & & (2,1) & \\
    & (3,1) & q(3,1) & q^{-1}(3,1) \\
\hline
F_2 & & & (2,1) \\
    & q^{-1}(2^2) & (2^2) & q(2^2) \\
\hline
\end{array} $$
$$ \begin{array}{|l||ccc|}
\hline
(1^3) & n\equiv 0 & n\equiv 1 & n\equiv 2 \\
\hline
\hline
F_0 & (1^3) & & \\
    & q(2,1^2) & q^{-1}(2,1^2) & (2,1^2) \\
\hline
F_1 & & (1^3) & \\
\hline
F_2 & & & (1^3) \\
    & q^{-1}(1^4) & (1^4) & q(1^4) \\
\hline
\end{array}
\quad 
\begin{array}{|l||ccc|}
\hline
\ov{\lb}\in\cP_4 & n\equiv 0 & n\equiv 1 & n\equiv 2 \\
\hline
\hline
F_0 & & \ov{\lb} & \\
\hline
F_1 & & & \ov{\lb} \\
\hline
F_2 & \ov{\lb} & & \\
\hline
\end{array} $$
\abs\hrulefill
\end{table}

\absb 
Next, given $\mu\in\cP_n^{3\text{-reg}}(\leq 4)$,
the elements $F^\mu\in\cU$ are as shown in Table \ref{llt3},
where we again use the round bracket notation, and silently assume that
for $\ov{\mu}\in\cP_m^{3\text{-reg}}$ we have $n-m\geq\ov{\mu}_1$,
or even $n-m\geq\ov{\mu}_1+1$ to ensure $3$-regularity.

\abs
We observe that the elements $B^\mu+\cF^{>4}_q\in\cF_q(\leq 4)$,
where $\ov{\mu}\in\coprod_{m=0}^4\cP_m^{3\text{-reg}}$ 
is fixed, are ultimately periodic, and only depend
on the residue class of $n$ modulo $3$; in the first column of
Table \ref{llt3} we give the bound where periodicity sets in.
But since we do not have an {\it a priori} proof of this fact, we
determine $B^\mu+\cF^{>4}_q$ explicitly, proving periodicity on the fly,
see also \ref{lltrem}.
We proceed through 
$\ov{\mu}\in\coprod_{m=0}^4\cP_m^{3\text{-reg}}$ with decreasing $m$
and in reversed lexicographical ordering:

\absb
For $\ov{\mu}\in\cP_4^{3\text{-reg}}$ we get $B^\mu+\cF^{>4}_q$ as 
indicated in Table \ref{llt3}. This is seen as follows:
For example, for $\ov{\mu}:=(2,1^2)$, applying 
$F_0F_2F_1^{(2)}F_0^{(2)}$ to the empty partition
yields $B^\mu\equiv\ov{\mu}$ for $n=6$. 
Then, successively applying $F_2,F_0,F_1,F_2,F_0,F_1,\ldots$,
and using Table \ref{foperators3}, by induction we get
$B^\mu\equiv\ov{\mu}$ for all $n\geq 6$.
For the other elements of $\cP_4^{3\text{-reg}}$ we argue similarly.

\begin{table}\caption{Products of divided power operators}\label{llt3}
$$ \begin{array}{|l|l|l|l|}
\hline
 & \ov{\mu} & F^\mu & B^\mu\equiv \\
\hline
\hline \vspc
n\geq 6 & (2,1^2)
 & F_0F_2F_1^{(2)}F_0^{(2)}\cdot F_2F_0F_1\cdot F_2F_0F_1\cdots
 & \ov{\mu} \\
n\geq 7 & (2^2)
 & F_0F_2F_1^{(2)}F_0F_2^{(2)}\cdot F_0F_1F_2\cdot F_0F_1F_2\cdots
 & \ov{\mu}+q(1^4) \\
n\geq 7 & (3,1) 
 & F_0F_2F_1^{(2)}F_0F_2F_1\cdot F_0F_1F_2\cdot F_0F_1F_2\cdots
 & \ov{\mu} \\
n\geq 8 & (4)
 & F_0F_2F_1F_0F_2F_1F_0F_2\cdot F_1F_2F_0\cdot F_1F_2F_0\cdots
 & \ov{\mu}+q(2^2) \\
\hline \vspc
n\geq 6 & (2,1)
 & F_0F_2F_1^{(2)}F_0F_2\cdot F_0F_1F_2\cdot F_0F_1F_2\cdots & \\
n\geq 7 & (3)
 & F_0F_2F_1F_0F_2F_1F_0\cdot F_1F_2F_0\cdot F_1F_2F_0\cdots & \\
n\geq 5 & (1^2)
 & F_0F_2F_1^{(2)}F_2\cdot F_0F_1F_2\cdot F_0F_1F_2\cdots & \\
n\geq 6 & (2)
 & F_0F_2F_1F_0F_2F_0\cdot F_1F_2F_0\cdot F_1F_2F_0\cdots & \\
n\geq 7 & (1)
 & F_0F_2F_1F_2F_0F_1F_2\cdot F_0F_1F_2\cdot F_0F_1F_2\cdots & \\
n\geq 7 & ()
 & F_0F_1F_2F_0F_1F_2F_0\cdot F_1F_2F_0\cdot F_1F_2F_0\cdots & \\
\hline
\end{array} $$
\abs\hrulefill
\end{table}

\absb
For $\ov{\mu}\in\cP_n^{3\text{-reg}}(\leq 3)$ we proceed similarly,
where the subsequent computations should be read as follows:
For example, for $\ov{\mu}:=(2,1)$, applying $F_0F_2F_1^{(2)}F_0F_2$ to the 
empty partition, we get $B^\mu+\cF^{>4}_q$ for $n=6$ as indicated below.
Then, proceeding by induction on $n\geq 6$, we successively apply
$F_0,F_1,F_2,F_0,F_1,F_2,\ldots$, where after each application of $F_0$
we additionally have to add a suitable multiple of $B^{(2,1^2)}+\cF^{>4}_q$; 
note that at this stage $B^{(2,1^2)}+\cF^{>4}_q$ has already been obtained
and proved to behave periodically. Closing the circle, the last line
shows that we indeed get $B^\mu+\cF^{>4}_q$ in terms of the residue class 
of $n$ modulo $3$.

\abs
For $\ov{\mu}:=(2,1)$ this for $n\geq 6$ yields:
$$ \begin{array}{clcll}
& B^\mu & \equiv 
& \ov{\mu}+q(1^3)+q(2^2)+q^2(1^4) 
& \text{if }n\equiv 0\pmod{3} \\
\stackrel{F_0}{\lmt} 
&&& \ov{\mu}+q(1^3)+(1+q^2)(2,1^2) \\
\stackrel{B^{(2,1^2)}}{\lmt}
& B^\mu & \equiv 
& \ov{\mu}+q(1^3)+q^2(2,1^2)
& \text{if }n\equiv 1\pmod{3} \\
\stackrel{F_1}{\lmt}
& B^\mu & \equiv 
& \ov{\mu}+q(1^3)+q(3,1)
& \text{if }n\equiv 2\pmod{3} \\
\stackrel{F_2}{\lmt}
&&& \ov{\mu}+q(1^3)+q(2^2)+q^2(1^4) \\
\end{array} $$

\abs
For $\ov{\mu}:=(3)$ this for $n\geq 7$ yields:
$$ \begin{array}{clcll}
& B^\mu & \equiv 
& \ov{\mu}+q(2,1)+q(2,1^2) 
& \text{if }n\equiv 1\pmod{3} \\
\stackrel{F_1}{\lmt} 
&&& \ov{\mu}+q(2,1)+(1+q^2)(3,1) \\
\stackrel{B^{(3,1)}}{\lmt}
& B^\mu & \equiv 
& \ov{\mu}+q(2,1)+q^2(3,1)
& \text{if }n\equiv 2\pmod{3} \\
\stackrel{F_2}{\lmt} 
& B^\mu & \equiv 
& \ov{\mu}+q(2,1)+q(4)+q^2(2^2)
& \text{if }n\equiv 0\pmod{3} \\
\stackrel{F_0}{\lmt} 
&&& \ov{\mu}+q(2,1)+q(2,1^2) \\
\end{array} $$

\abs
For $\ov{\mu}:=(1^2)$ this for $n\geq 5$ yields:
$$ \begin{array}{clcll}
& B^\mu & \equiv 
& \ov{\mu} 
& \text{if }n\equiv 2\pmod{3} \\
\stackrel{F_0}{\lmt} 
& B^\mu & \equiv 
& \ov{\mu}+q(2,1)+q^2(1^3)
& \text{if }n\equiv 0\pmod{3} \\
\stackrel{F_1}{\lmt} 
& B^\mu & \equiv 
& \ov{\mu}+q(3,1)
& \text{if }n\equiv 1\pmod{3} \\
\stackrel{F_2}{\lmt} 
&&& \ov{\mu} \\
\end{array} $$

\abs
For $\ov{\mu}:=(2)$ this for $n\geq 6$ yields:
$$ \begin{array}{clcll}
& B^\mu & \equiv 
& \ov{\mu} 
& \text{if }n\equiv 0\pmod{3} \\
\stackrel{F_1}{\lmt} 
& B^\mu & \equiv 
& \ov{\mu}+q(3)+q^2(2,1)
& \text{if }n\equiv 1\pmod{3} \\
\stackrel{F_2}{\lmt} 
& B^\mu & \equiv 
& \ov{\mu}+q(4)+q^2(2^2)
& \text{if }n\equiv 2\pmod{3} \\
\stackrel{F_0}{\lmt} 
&&& \ov{\mu} \\
\end{array} $$

\abs
For $\ov{\mu}:=(1)$ this for $n\geq 7$ yields:
$$ \begin{array}{clcll}
& B^\mu & \equiv 
& \ov{\mu}+q(2^2)
& \text{if }n\equiv 1\pmod{3} \\
\stackrel{F_0}{\lmt} 
& B^\mu & \equiv 
& \ov{\mu}+q(2)+q(2^2)
& \text{if }n\equiv 2\pmod{3} \\
\stackrel{F_1}{\lmt} 
& B^\mu & \equiv 
& \ov{\mu}+q(1^2)+q(3)+q^2(2,1)+q(2^2)
& \text{if }n\equiv 0\pmod{3} \\
\stackrel{F_2}{\lmt} 
&&& \ov{\mu}+(4)+2q(2^2) \\ 
\stackrel{B^{(4)}}{\lmt}
&&& \ov{\mu}+q(2^2) \\
\end{array} $$

\abs
For $\ov{\mu}:=()$ this for $n\geq 7$ yields:
$$ \begin{array}{clcll}
& B^\mu & \equiv 
& \ov{\mu}+q(2)+q(2,1)+q(2,1^2)
& \text{if }n\equiv 1\pmod{3} \\
\stackrel{F_1}{\lmt} 
&&& \ov{\mu}+(3)+2q(2,1)+q(3,1) \\
\stackrel{B^{(3)}}{\lmt}
& B^\mu & \equiv 
& \ov{\mu}+q(2,1)
& \text{if }n\equiv 2\pmod{3} \\
\stackrel{F_2}{\lmt} 
& B^\mu & \equiv 
& \ov{\mu}+q(1)+q(2,1)+q(2^2)
& \text{if }n\equiv 0\pmod{3} \\
\stackrel{F_0}{\lmt} 
&&& \ov{\mu}+q(2)+q(2,1)+q(2,1^2) \\
\end{array} $$

\absb
In conclusion, this yields the crystallized decomposition matrix
$\cD_n^q(\leq 4)$, for $n\geq 8$, as exhibited in Table \ref{dec3llt}.
Here, for $i\in\{0,1,2\}$ running through the residue classes modulo $3$,
reminiscent of the Kronecker symbol we let
$$ \dt_i=\dt_i(n):=\left\{\begin{array}{ll} 
1, &\text{if }n\equiv i\pmod{3}, \\
0, & \text{otherwise}. \\
\end{array}\right. $$
We remark that for $n\leq 10$ these results are also contained
in the explicit crystallized decomposition matrices given in 
\cite[Sect.10.4]{LasLecThi}.

\begin{table}\caption{Crystallized decomposition matrix for $\cP_n(\leq 4)$.} 
\label{dec3llt}
$$ \begin{array}{|l||r|r|rr|rr|rrrr|}
\hline \rule{0em}{3em}
    & \spc \begin{rotate}{90} $()$ \end{rotate}
    & \spc \begin{rotate}{90} $(1)$ \end{rotate}
    & \spc \begin{rotate}{90} $(2)$ \end{rotate}
    & \spc \begin{rotate}{90} $(1^2)$ \end{rotate}
    & \spc \begin{rotate}{90} $(3)$ \end{rotate}
    & \spc \begin{rotate}{90} $(2,1)$ \end{rotate}
    & \spc \begin{rotate}{90} $(4)$ \end{rotate}
    & \spc \begin{rotate}{90} $(3,1)$ \end{rotate}
    & \spc \begin{rotate}{90} $(2^2)$ \end{rotate}
    & \spc \begin{rotate}{90} $(2,1^2)$ \end{rotate} \\
\hline
\hline
() \vspc &
1 & & & & & & & & & \\
\hline
(1) \vspc &
q\dt_0 & 1 & & & & & & & & \\
\hline
(2) \vspc &
q\dt_1 & q\dt_2 & 1 & & & & & & & \\
(1^2) \vspc &
. & q\dt_0 & . & 1 & & & & & & \\
\hline
(3) \vspc &
. & q\dt_0 & q\dt_1 & . & 1 & & & & & \\
(2,1) \vspc &
q & q^2\dt_0 & q^2\dt_1 & q\dt_0 & q & 1 & & & & \\
(1^3) \vspc &
. & . & . & q^2\dt_0 & . & q & & & & \\
\hline
(4) \vspc &
. & . & q\dt_2 & . & q\dt_0 & . & 1 & & & \\
(3,1) \vspc &
. & . & . & q\dt_1 & q^2\dt_2 & q\dt_2 & . & 1 & & \\
(2^2) \vspc &
q^2\dt_0 & q & q^2\dt_2 & . & q^2\dt_0 & q\dt_0 & q & . & 1 & \\
(2,1^2) \vspc &
q\dt_1 & . & . & . & q\dt_1 & q^2\dt_1 & . & . & . & 1 \\
(1^4) \vspc &
. & . & . & . & . & q^2\dt_0 & . & . & q & . \\
\hline
\end{array} $$
\abs\hrulefill
\end{table}

\AbsT{Decomposition matrices.}\label{decmat3}
We proceed to determine the decomposition matrix $\cD_n(\leq 4)$.
The result is given in Table \ref{dec3}, where we assume that $n\geq 8$.
Here, we again use the Kronecker notation introduced in \ref{crystdec3},
and let
$$ \al=\al(n):=\left\{\begin{array}{ll} 
1, &\text{if }n\equiv 2,3,4\pmod{9}, \\
0, & \text{otherwise}, \\
\end{array}\right. $$
and
$$ \bt=\bt(n):=\left\{\begin{array}{ll} 
1, &\text{if }n\equiv 4,5,6\pmod{9}, \\
0, & \text{otherwise}. \\
\end{array}\right. $$
The decomposition matrices $\cD_n$ for $n\leq 7$ are known,
and for example given in \cite[p.143]{JamLN}, or accessible 
in the databases mentioned in Section \ref{intro}. 

\abs
To determine the decomposition matrix $\cD_n(\leq 4)$, we make use
of the crystallized decomposition matrix $\cD_n^q(\leq 4)$. 
While computing the entries of $\cD_n(\leq 4)$
we also determine the adjustment matrix $\cA_n(\leq 4)$, the result
is given in Table \ref{dec3adj}. In the subsequently given details 
we make use of the following notation:
Given $\lb=[\lb_1,\ov{\lb}]\in\cP_n$ the associated Specht module 
is also written as $S_n(\ov{\lb})$, and similarly for  
$\mu=[\mu_1,\ov{\mu}]\in\cP_n^{3\text{-reg}}$, the associated 
irreducible module is also written as $D_n(\ov{\mu})$.

\begin{table}\caption{Decomposition matrix for $\cP_n(\leq 4)$.}\label{dec3}
$$ \begin{array}{|l||r|r|rr|rr|rrrr|}
\hline \rule{0em}{3em}
    & \spc \begin{rotate}{90} $()$ \end{rotate}
    & \spc \begin{rotate}{90} $(1)$ \end{rotate}
    & \spc \begin{rotate}{90} $(2)$ \end{rotate}
    & \spc \begin{rotate}{90} $(1^2)$ \end{rotate}
    & \spc \begin{rotate}{90} $(3)$ \end{rotate}
    & \spc \begin{rotate}{90} $(2,1)$ \end{rotate}
    & \spc \begin{rotate}{90} $(4)$ \end{rotate}
    & \spc \begin{rotate}{90} $(3,1)$ \end{rotate}
    & \spc \begin{rotate}{90} $(2^2)$ \end{rotate}
    & \spc \begin{rotate}{90} $(2,1^2)$ \end{rotate} \\
\hline
\hline
() \vspc &
1 & & & & & & & & & \\
\hline
(1) \vspc &
\dt_0 & 1 & & & & & & & & \\
\hline
(2) \vspc &
\dt_1 & \dt_2 & 1 & & & & & & & \\
(1^2) \vspc &
. & \dt_0 & . & 1 & & & & & & \\
\hline
(3) \vspc &
\al & \dt_0 & \dt_1 & . & 1 & & & & & \\
(2,1) \vspc &
1\!+\!\al & \dt_0 & \dt_1 & \dt_0 & 1 & 1 & & & & \\
(1^3) \vspc &
. & . & . & \dt_0 & . & 1 & & & & \\
\hline
(4) \vspc &
\al\dt_0 & \bt & \dt_2 & . & \dt_0 & . & 1 & & & \\
(3,1) \vspc &
(\al\!+\!\bt)\dt_2 & . & . & \dt_1 & \dt_2 & \dt_2 & . & 1 & & \\
(2^2) \vspc &
(1\!+\!\al)\dt_0 & 1\!+\!\bt & \dt_2 & . & \dt_0 & \dt_0 & 1 & . & 1 & \\
(2,1^2) \vspc &
(1\!+\!\al)\dt_1 & . & . & . & \dt_1 & \dt_1 & . & . & . & 1 \\
(1^4) \vspc &
. & . & . & . & . & \dt_0 & . & . & 1 & . \\
\hline
\end{array} $$
\abs\hrulefill
\end{table}

\begin{table}\caption{Adjustment matrix for $\cP_n(\leq 4)$.}\label{dec3adj}
$$ \begin{array}{|l||r|r|rr|rr|rrrr|}
\hline \rule{0em}{3em}
    & \spc \begin{rotate}{90} $()$ \end{rotate}
    & \spc \begin{rotate}{90} $(1)$ \end{rotate}
    & \spc \begin{rotate}{90} $(2)$ \end{rotate}
    & \spc \begin{rotate}{90} $(1^2)$ \end{rotate}
    & \spc \begin{rotate}{90} $(3)$ \end{rotate}
    & \spc \begin{rotate}{90} $(2,1)$ \end{rotate}
    & \spc \begin{rotate}{90} $(4)$ \end{rotate}
    & \spc \begin{rotate}{90} $(3,1)$ \end{rotate}
    & \spc \begin{rotate}{90} $(2^2)$ \end{rotate}
    & \spc \begin{rotate}{90} $(2,1^2)$ \end{rotate} \\
\hline
\hline
() \vspc      &
1 & & & & & & & & & \\
\hline
(1) \vspc     &
. & 1 & & & & & & & & \\
\hline
(2) \vspc     &
. & . & 1 & & & & & & & \\
(1^2) \vspc   &
. & . & . & 1 & & & & & & \\
\hline
(3) \vspc     &
\al & . & . & . & 1 & & & & & \\
(2,1) \vspc   &
. & . & . & . & . & 1 & & & & \\
\hline
(4) \vspc     &
. & \bt & . & . & . & . & 1 & & & \\
(3,1) \vspc   &
\bt\dt_2 & . & . & . & . & . & . & 1 & & \\
(2^2) \vspc   &
. & . & . & . & . & . & . & . & 1 & \\
(2,1^2) \vspc &
. & . & . & . & . & . & . & . & . & 1 \\
\hline
\end{array} $$
\abs\hrulefill
\end{table}

\absb
The decomposition numbers of the Specht modules $S^{[n-m,m]}=S_n(m)$, where 
$m\in\{0\ld\fl{\frac{n}{2}}\}$, that is those belonging to partitions
having at most two parts, are known by \cite{JamI},
see also \cite[Thm.24.15]{JamLN}.

\absb
The decomposition numbers of the Specht modules $S^{[n-m,1^m]}=S_n(1^m)$, 
where $m\in\{0\ld n-1\}$, that is those belonging to hook partitions, 
are known by \cite{Pee}, see also \cite[Thm.24.1]{JamLN}.
More precisely, if $n\not\equiv 0\pmod{3}$, then the Specht module 
$S_n(1^m)$ is reducible, hence by \cite[Thm.6.3.50]{JamKer} is 
isomorphic to $D^{[n-m,1^m]^R}$.
If $n\equiv 0\pmod{3}$, then for $m\in\{1\ld n-2\}$ the Specht module 
$S_n(1^m)$ is uniserial with two distinct constituents, 
such that $\Hom_{\F_3\cS_n}(S_n(1^m),S_n(1^{m+1}))\neq\{0\}$ 
for $m\in\{0\ld n-2\}$; hence by induction on $m\in\{0\ld n-1\}$,
again using \cite[Thm.6.3.50]{JamKer}, it follows that the head 
constituent of $S_n(1^m)$ is isomorphic to $D^{[n-m,1^m]^R}$.
 
\absb
Whenever $\lb\in\cP_n(m)$ and $\mu\in\cP_n^{3\text{-reg}}(m)$ for 
some $m\in\{0\ld n-1\}$, that is $\lb_1=\mu_1=n-m$, the principle
of first row removal, see \cite{JamIII}, yields 
$[S_n(\ov{\lb})\cn D_n(\ov{\mu})]=[S^{\ov{\lb}}\cn D^{\ov{\mu}}]$,
where the latter are decomposition numbers of $\cS_m$, which 
since $m\leq 4$ are easily determined or can be 
looked up in \cite[p.143]{JamLN}.

\absb
Whenever $\lb\in\cP_n(3)\dcup\cP_n(4)$ and 
$\mu\in\cP_n^{3\text{-reg}}(2)\dcup\cP_n^{3\text{-reg}}(3)$ 
such that $l(\lb)=l(\mu)=3$, the principle of first column removal, 
see \cite{JamIII}, yields 
$[S_n(\ov{\lb})\cn D_n(\ov{\mu})]
=[S_{n-3}(\lb_2-1,\lb_3-1)\cn D_{n-3}(\mu_2-1,\mu_3-1)]$.
This settles the cases
$$ \begin{array}{lcl}
\mbox{}[S_n(2,1)\cn D_n(1^2)]&=&[S_{n-3}(1)\cn D_{n-3}(())], \\
\mbox{}[S_n(3,1)\cn D_n(2,1)]&=&[S_{n-3}(2)\cn D_{n-3}(1)],
\rule{0em}{1.2em} \\
\mbox{}[S_n(3,1)\cn D_n(1^2)]&=&[S_{n-3}(2)\cn D_{n-3}(())], 
\rule{0em}{1.2em} \\
\mbox{}[S_n(2^2)\cn D_n(2,1)]&=&[S_{n-3}(1^2)\cn D_{n-3}(1)], 
\rule{0em}{1.2em} \\
\mbox{}[S_n(2^2)\cn D_n(1^2)]&=&[S_{n-3}(1^2)\cn D_{n-3}(())], 
\rule{0em}{1.2em} \\
\end{array} $$
where the decomposition numbers of $\cS_{n-3}$ involved are two-part 
or hook partition cases which we have already dealt with above.
  
\absb
In order to complete $\cD_n(\leq 3)$ 
it remains to determine the row belonging to the partition $(2,1)$.
But from $\cD_n(\leq 3)=\cD^\zt_n(\leq 3)\cdot\cA_n(\leq 3)$,
using the entries of $\cD_n(\leq 3)$ already known, we conclude that 
$\cA_n(\leq 3)$ actually coincides with the upper left-hand 
$(6\tm 6)$-sub-matrix of the matrix given in Table \ref{dec3adj}.
This in turn determines $\cD_n(\leq 3)$.

\absb
We turn our attention to the rows of $\cD_n(\leq 4)$
belonging to partitions in $\cP_n(4)$: From
$\cD_n(\leq 4)=\cD^\zt_n(\leq 4)\cdot\cA_n(\leq 4)$,
a comparison of the known entries of $\cD_n(\leq 4)$ with 
the matrix $\cD^\zt_n(\leq 4)$ shows that the 
rows of $\cA_n(\leq 4)$ belonging to
the partitions $(4)$ and $(2^2)$ are as given in Table \ref{dec3adj}.
This in turn completes the row of $\cD_n(\leq 4)$ belonging to
the partition $(2^2)$.

\absb
Considering $3$-contents, reflecting the distribution of modules into
$3$-blocks, we conclude that
$[S_n(2,1^2)\cn D_n(\ov{\mu})]=0$ for $\ov{\mu}\in\{(1^2),(1)\}$, and
$[S_n(3,1)\cn D_n(\ov{\mu})]=0$ for $\ov{\mu}\in\{(2),(1)\}$.
This leaves the six entries
$[S_n(2,1^2)\cn D_n(\ov{\mu})]=0$ for $\ov{\mu}\in\{(2,1),(3),(2),()\}$,
and $[S_n(3,1)\cn D_n(\ov{\mu})]=0$ for $\ov{\mu}\in\{(3),()\}$
to be determined, which we deal with in turn:
First note that for the cases mentioned a consideration of $3$-contents
shows that $[S_n(2,1^2)\cn D_n(\ov{\mu})]>0$ only if $n\equiv 1\pmod{3}$,
and that $[S_n(3,1)\cn D_n(\ov{\mu})]>0$ only if $n\equiv 2\pmod{3}$.

\absb
We consider $[S_n(2,1^2)\cn D_n(2,1)]$ for $n\equiv 1\pmod{3}$: 
Applying \ref{jamestrick} with $0$-removal leads to the inequality
$$ \begin{array}{cl}
& [S_n(2,1^2)\cn D_n(2,1)]\cdot [D_n(2,1)\da_0\cn D_{n-1}(1^2)] \\
\leq & [S_{n-1}(2,1)\cn D_{n-1}(1^2)]+[S_{n-1}(1^3)\cn D_{n-1}(1^2)]=2.
\rule{0em}{1.2em} \\ \end{array} $$
Since by \ref{modbranch} we have
$[D_n(2,1)\da_0\cn D_{n-1}(1^2)]=2$, using Table \ref{dec3llt} implies
$1=[S_n(2,1^2)\cn D^\zt_n(2,1)]\leq [S_n(2,1^2)\cn D_n(2,1)]\leq 1$.

\absb
We consider $[S_n(2,1^2)\cn D_n(3)]$ for $n\equiv 1\pmod{3}$: 
Again applying \ref{jamestrick} with $0$-removal leads to the inequality
$$ \begin{array}{cl}
& [S_n(2,1^2)\cn D_n(3)]\cdot [D_n(3)\da_0\cn D_{n-1}(3)] \\
\leq & [S_{n-1}(2,1)\cn D_{n-1}(3)]+[S_{n-1}(1^3)\cn D_{n-1}(3)]=1,
\rule{0em}{1.2em} \\ \end{array} $$
which using Table \ref{dec3llt} implies
$1=[S_n(2,1^2)\cn D^\zt_n(3)]\leq [S_n(2,1^2)\cn D_n(3)]\leq 1$.

\absb
We consider $[S_n(2,1^2)\cn D_n(2)]$ for $n\equiv 1\pmod{3}$: 
Applying \ref{jamestrick} with $0$-removal, since by \ref{modbranch}
we have $[D_n(2)\da_0\cn D_{n-1}(1)]=1$, we get the inequality
$$ \begin{array}{cl}
& [S_n(2,1^2)\cn D_n(2)]
+[S_n(2,1^2)\cn D_n(2,1)]\cdot[D_n(2,1)\da_0\cn D_{n-1}(1)] \\
\leq & [S_{n-1}(2,1)\cn D_{n-1}(1)]+[S_{n-1}(1^3)\cn D_{n-1}(1)]=1.
\rule{0em}{1.2em} \\ \end{array} $$
We have $[S_n(2,1^2)\cn D_n(2,1)]=1$, and from 
$[D_n(2,1)\da_0]=[S_n(1^3)\da_0]=[S_{n-1}(1^3)]+[S_{n-1}(1^2)]$ 
we get 
$$ [D_n(2,1)\da_0\cn D_{n-1}(1)]
=[S_{n-1}(1^3)\cn D_{n-1}(1)]+[S_{n-1}(1^2)\cn D_{n-1}(1)]=1, $$
hence we infer $[S_n(2,1^2)\cn D_n(2)]=0$.

\absb
We consider $[S_n(2,1^2)\cn D_n(())]$ for $n\equiv 1\pmod{3}$: 
Applying \ref{jamestrick} with $0$-removal, we get the inequality
$$ \begin{array}{lcl}
[S_n(2,1^2)\cn D_n(())]
& \leq & [S_{n-1}(2,1)\cn D_{n-1}(())]+[S_{n-1}(1^3)\cn D_{n-1}(())] \\
& = & \left\{\begin{array}{ll}
2, & \text{if }n\equiv 4\pmod{9}, \\
1, & \text{if }n\equiv 1,7\pmod{9}. \\
\end{array}\right. \rule{0em}{2em} \\ \end{array} $$
On the other hand, by what we already know about $\cA_n(\leq 4)$,
we conclude that
$$ [S_n(2,1^2)\cn D^\zt_n(())]+\al\cdot [S_n(2,1^2)\cn D^\zt_n(3)]
=(1+\al)\dt_1\leq [S_n(2,1^2)\cn D_n(())] ,$$ 
implying equality.
This completes the row of $\cD_n(\leq 4)$ belonging to the partition
$(2,1^2)$, and shows that the row of $\cA_n(\leq 4)$ belonging
to $(2,1^2)$ is as in Table \ref{dec3adj}.

\absb
We consider $[S_n(3,1)\cn D_n(3)]$ for $n\equiv 2\pmod{3}$: 
Applying \ref{jamestrick} with $1$-removal, we get the inequality
$$ \begin{array}{cl}
& [S_n(3,1)\cn D_n(3)]\cdot [D_n(3)\da_1\cn D_{n-1}(2)] \\
\leq & [S_{n-1}(3)\cn D_{n-1}(2)]+[S_{n-1}(2,1)\cn D_{n-1}(2)]=2. 
\rule{0em}{1.2em} \\ \end{array} $$
Since from \ref{modbranch} we obtain
$[D_n(3)\da_1\cn D_{n-1}(2)]=2$, using Table \ref{dec3llt} thus implies
$1=[S_n(3,1)\cn D^\zt_n(3)]\leq [S_n(3,1)\cn D_n(3)]\leq 1$.

\absb
We consider $[S_n(3,1)\cn D_n(())]$ for $n\equiv 2\pmod{3}$: 
Applying \ref{jamestrick} with $1$-removal, 
since we have $[S_n(3,1)\cn D_n(3)]=1$, and
$[D_n(2,1)\da_1\cn D_{n-1}(())]=0$ and $[D_n(3,1)\da_1\cn D_{n-1}(())]=0$
by \ref{modbranch}, we this time even get the strict equality
$$ \begin{array}{cl}
& [S_n(3,1)\cn D_n(())]+[D_n(3)\da_1\cn D_{n-1}(())] \\
= & [S_{n-1}(3)\cn D_{n-1}(())]+[S_{n-1}(2,1)\cn D_{n-1}(())] 
\rule{0em}{1.2em} \\
= & \left\{\begin{array}{ll}
3, & \text{if }n\equiv 5\pmod{9}, \\
1, & \text{if }n\equiv 2,8\pmod{9}. \\
\end{array}\right. \rule{0em}{2em} \\ \end{array} $$
If $n\equiv 5,8\pmod{9}$ then 
$[D_n(3)\da_1]=[S_n(3)\da_1]=[S_{n-1}(2)]+[S_{n-1}(3)]$, hence 
$$ \begin{array}{cl}
& [D_n(3)\da_1\cn D_{n-1}(())] \\
= & [S_{n-1}(2)\cn D_{n-1}(())]+[S_{n-1}(3)\cn D_{n-1}(())] 
\rule{0em}{1.2em} \\
= & \left\{\begin{array}{ll}
2, & \text{if }n\equiv 5\pmod{9}, \\
1, & \text{if }n\equiv 8\pmod{9}. \\
\end{array}\right. \rule{0em}{2em} \\ \end{array} $$
Similarly, if $n\equiv 2\pmod{9}$ then 
$[D_n(3)\da_1]=[S_n(3)\da_1]-[S_n(())\da_1]
=[S_{n-1}(2)]+[S_{n-1}(3)]-[S_{n-1}(())]$, hence 
$$ \begin{array}{cl}
& [D_n(3)\da_1\cn D_{n-1}(())] \\
= & [S_{n-1}(2)\cn D_{n-1}(())]+[S_{n-1}(3)\cn D_{n-1}(())]
-[S_{n-1}(())\cn D_{n-1}(())]=0.
\rule{0em}{1.2em} \\ \end{array} $$
In conclusion we thus have
$$ [S_n(3,1)\cn D_n(())]
=\left\{\begin{array}{ll}
1, & \text{if }n\equiv 2,5\pmod{9}, \\
0, & \text{if }n\equiv 8\pmod{9}, \\
\end{array}\right. $$
yielding $[S_n(3,1)\cn D_n(())]=(\al+\bt)\dt_2$.
This completes $\cD_n(\leq 4)$, and shows that the row of 
$\cA_n(\leq 4)$ belonging to the partition $(3,1)$ 
is as given in Table \ref{dec3adj}.

\AbsT{Degree formulae.}\label{deg3formulae}
We are now prepared to obtain degree formulae for the irreducible
modular representations parameterized by $\cP_n^{3\text{-reg}}(\leq 4)$.
They are given in Tables \ref{deg3leq3} and \ref{deg3eq4}, where 
for $\mu\in\cP_n^{3\text{-reg}}$ we let $d^\mu:=\dim_{\F_3}(D^\mu)$:

\abs
In view of the decomposition matrix $\cD_n(\leq 4)$ for $n\geq 8$ 
given in Table \ref{dec3}, and the known decomposition numbers 
for $n\leq 7$, see the comment at the beginning of \ref{decmat3}, 
these follow straightforwardly from the hook length formula 
for the dimension of Specht modules, 
see \cite{FraRobThr} and also \cite[Thm.20.1]{JamLN}.

\begin{table}\caption{Degree formulae for $\cP_n^{3\text{-reg}}(\leq 3)$.} 
             \label{deg3leq3}
$$ \begin{array}{|l|l|lr|}
\hline
\mu & d^\mu & \multicolumn{2}{l|}{\text{condition}} \\
\hline
\hline
[n] & 1 & & \vspc \\
\hline
[n-1,1]   & n-2 & n\equiv 0   & \pmod{3} \vspc \\ 
(n\geq 2) & n-1 & n\equiv 1,2 & \pmod{3} \vspc \\
\hline
[n-2,2]   & \frac{1}{2}(n^2-5n+2) & n\equiv 2 & \pmod{3} \vspc \\
(n\geq 4) & \frac{1}{2}(n^2-3n-2) & n\equiv 1 & \pmod{3} \vspc \\
          & \frac{1}{2}(n^2-3n)   & n\equiv 0 & \pmod{3} \vspc \\
\hline
[n-2,1^2] & \frac{1}{2}(n^2-5n+6) & n\equiv 0 & \pmod{3} \vspc \\
(n\geq 4) & \frac{1}{2}(n^2-3n+2) & n\equiv 1,2 & \pmod{3} \vspc \\
\hline
[n-3,3]   & \frac{1}{6}(n^3-9n^2+14n)   & n\equiv 4 & \pmod{9} \vspc \\
(n\geq 6) & \frac{1}{6}(n^3-9n^2+14n+6) & n\equiv 1,7 & \pmod{9} \vspc \\
          & \frac{1}{6}(n^3-6n^2-n+6)   & n\equiv 3 & \pmod{9} \vspc \\
          & \frac{1}{6}(n^3-6n^2-n+12)  & n\equiv 0,6 & \pmod{9} \vspc \\
          & \frac{1}{6}(n^3-6n^2+5n-6)  & n\equiv 2 & \pmod{9} \vspc \\
          & \frac{1}{6}(n^3-6n^2+5n)    & n\equiv 5,8 & \pmod{9} \vspc \\
\hline
[n-3,2,1] & \frac{1}{6}(n^3-9n^2+26n-24) & n\equiv 0 & \pmod{3} \vspc \\
(n\geq 5) & \frac{1}{6}(n^3-6n^2+11n-6)  & n\equiv 1,2 & \pmod{3} \vspc \\
\hline
\end{array} $$
\abs\hrulefill
\end{table}

\begin{table}\caption{Degree formulae for $\cP_n^{3\text{-reg}}(4)$.} 
             \label{deg3eq4}
$$ \begin{array}{|l|l|lr|}
\hline
\mu & d^\mu & \multicolumn{2}{l|}{\text{condition}} \\
\hline
\hline
[n-4,4] & \frac{1}{24}(n^4-14n^3+47n^2-34n) & n\equiv 6 & \pmod{9} \vspc \\
(n\geq 8) 
      & \frac{1}{24}(n^4-14n^3+47n^2-10n-48) & n\equiv 0,3 & \pmod{9} \vspc \\
      & \frac{1}{24}(n^4-10n^3+11n^2+22n) & n\equiv 5 & \pmod{9} \vspc \\
      & \frac{1}{24}(n^4-10n^3+11n^2+46n-24) & n\equiv 2,8 & \pmod{9} \vspc \\
      & \frac{1}{24}(n^4-10n^3+23n^2-38n+24) & n\equiv 5 & \pmod{9} \vspc \\
      & \frac{1}{24}(n^4-10n^3+23n^2-14n) & n\equiv 1,7 & \pmod{9} \vspc \\
\hline
[n-4,3,1] & \frac{1}{24}(3n^4-38n^3+129n^2-118n) & n\equiv 5 & \pmod{9} \vspc \\
(n\geq 7) & \frac{1}{24}(3n^4-38n^3+129n^2-118n+24) 
          & n\equiv 2,8 & \pmod{9} \vspc \\
        & \frac{1}{24}(3n^4-30n^3+69n^2-18n-24) & n\equiv 1 & \pmod{3} \vspc \\
        & \frac{1}{24}(3n^4-30n^3+81n^2-54n)    & n\equiv 0 & \pmod{3} \vspc \\
\hline
[n-4,2^2] & \frac{1}{24}(n^4-14n^3+71n^2-154n+120) 
          & n\equiv 0 &\pmod{3} \vspc \\
(n\geq 7) & \frac{1}{24}(n^4-10n^3+35n^2-50n+24) 
          & n\equiv 1,2 &\pmod{3} \vspc \\
\hline
[n-4,2,1^2] & \frac{1}{24}(3n^4-38n^3+153n^2-190n-24) 
            & n\equiv 1 & \pmod{3} \vspc \\
(n\geq 6)   & \frac{1}{24}(3n^4-30n^3+93n^2-90n) 
            & n\equiv 0,2 & \pmod{3} \vspc \\
\hline
\end{array} $$
\abs\hrulefill
\end{table}

\section{Decomposition numbers in characteristic $2$}\label{char2}

\abs
In order to get an overview over the irreducible representations
parameterized by $\mu\in\cP_n^{2\text{-reg}}(\leq 4)$, we determine
the crystallized decomposition matrices $\cD_n^q(\leq 4)$ and the 
decomposition matrices $\cD_n(\leq 4)$ for $p=2$.

\AbsT{Crystallized decomposition matrices.}\label{crystdec2}
We apply the LLT algorithm to the truncated Fock space,
according to the description in \ref{llt} and \ref{llttrunc}.
We proceed entirely similar to \ref{crystdec3}, and thus only
record the results:

\absb
The maps $F_i\cn\cF_q(\leq 4)\ra\cF_q(\leq 4)$,
where $i\in\{0,1\}$ runs through the residue classes modulo $2$,
are given in Table \ref{foperators2}.

\begin{table}\caption{Action on truncated Fock space.}\label{foperators2}
$$ \begin{array}{|l||cc|}
\hline
() & n\equiv 0 & n\equiv 1 \\
\hline
\hline
F_0 & () & \\
\hline
F_1 & & () \\
    & q^{-1}(1) & q(1) \\
\hline
\end{array} 
\quad
\begin{array}{|l||cc|}
\hline
(1) & n\equiv 0 & n\equiv 1 \\
\hline
\hline
F_0 & & (1) \\
    & q^{-1}(2) & q(2) \\
    & (1^2) & q^2(1^2) \\
\hline
F_1 & (1) & \\
\hline
\end{array} $$
$$ \begin{array}{|l||cc|}
\hline
(2) & n\equiv 0 & n\equiv 1 \\
\hline
\hline
F_0 & (2) & \\
    & (2,1) & q^{-2}(2,1) \\
\hline
F_1 & & (2) \\
    & q^{-1}(3) & q(3) \\
\hline
\end{array} 
\quad
\begin{array}{|l||cc|}
\hline
(1^2) & n\equiv 0 & n\equiv 1 \\
\hline
\hline
F_0 & (1^2) & \\
    & q(2,1) & q^{-1}(2,1) \\
\hline
F_1 & & (1^2) \\
    & q^{-1}(1^3) & q(1^3) \\
\hline
\end{array} $$
$$ \begin{array}{|l||cc|}
\hline
(3) & n\equiv 0 & n\equiv 1 \\
\hline
\hline
F_0 & & (3) \\
    & q^{-1}(4) & q(4) \\
    & (3,1) & q^2(3,1) \\
\hline
F_1 & (3)  & \\
\hline
\end{array} 
\quad
\begin{array}{|l||cc|}
\hline
(2,1) & n\equiv 0 & n\equiv 1 \\
\hline
\hline
F_0 & & (2,1) \\
\hline
F_1 & (2,1) & \\
    & q(3,1) & q^{-1}(3,1) \\
    & q^2(2^2) & (2^2) \\
    & q^3(2,1^2) & q(2,1^2) \\
\hline
\end{array} $$
$$ \begin{array}{|l||cc|}
\hline
(1^3) & n\equiv 0 & n\equiv 1 \\
\hline
\hline
F_0 & & (1^3) \\
    & q^{-1}(2,1^2) & q(2,1^2) \\
    & (1^4) & q^2(1^4) \\
\hline
F_1 & (1^3) & \\
\hline
\end{array} 
\quad
\begin{array}{|l||cc|}
\hline
\ov{\lb}\in\cP_4 & n\equiv 0 & n\equiv 1 \\
\hline
\hline
F_0 & \ov{\lb} & \\
\hline
F_1 & & \ov{\lb} \\
\hline
\end{array} $$
\abs\hrulefill
\end{table}

\absb 
Given $\mu\in\cP_n^{2\text{-reg}}(\leq 4)$,
the elements $F^\mu\in\cU$ are as shown in Table \ref{llt2}.
Again we observe that the elements $B^\mu+\cF^{>4}_q\in\cF_q(\leq 4)$,
where $\ov{\mu}\in\coprod_{m=0}^4\cP_m^{2\text{-reg}}$ 
is fixed, are ultimately periodic, and only depend
on the residue class of $n$ modulo $2$; in the first column of
Table \ref{llt2} we give the bound where periodicity sets in.

\absb
For $\ov{\mu}\in\cP_4^{2\text{-reg}}$ we get $B^\mu+\cF^{>4}_q$ 
as follows:
$$ \begin{array}{lcl}
B^{(3,1)} & \equiv & (3,1) + q(2^2) + q^2(2,1^2) \\
B^{(4)}   & \equiv & (4) + q(3,1) + q(2,1^2) + q^2(1^4) \\
\end{array} $$

\begin{table}\caption{Products of divided power operators}\label{llt2}
$$ \begin{array}{|l|l|l|}
\hline
 & \ov{\mu} & F^\mu \\
\hline
\hline \vspc
n\geq 8 & (3,1) 
 & F_0F_1^{(2)}F_0^{(3)}F_1^{(2)}\cdot F_0F_1\cdot F_0F_1\cdots \\
%
n\geq 9 & (4)
 & F_0F_1^{(2)}F_0^{(2)}F_1^{(2)}F_0^{(2)}\cdot F_1F_0\cdot F_1F_0\cdots \\
\hline \vspc
n\geq 6 & (2,1)
 & F_0F_1^{(2)}F_0^{(3)}\cdot F_1F_0\cdot F_1F_0\cdots \\
n\geq 7 & (3)
 & F_0F_1^{(2)}F_0^{(2)}F_1^{(2)}\cdot F_0F_1\cdot F_0F_1\cdots \\
n\geq 7 & (2)
 & F_0F_1^{(2)}F_0^{(2)}F_1F_0\cdot F_1F_0\cdot F_1F_0\cdots \\
n\geq 5 & (1)
 & F_0F_1^{(2)}F_0F_1\cdot F_0F_1\cdot F_0F_1\cdots \\
n\geq 5 & ()
 & F_0F_1F_0F_1F_0\cdot F_1F_0\cdot F_1F_0\cdots \\
\hline
\end{array} $$
\abs\hrulefill
\end{table}

\absb
For $\ov{\mu}\in\cP_n^{2\text{-reg}}(\leq 3)$ we get the following: 

\abs
For $\ov{\mu}:=(2,1)$ this for $n\geq 6$ yields:
$$ \begin{array}{clcll}
& B^\mu & \equiv 
& \ov{\mu}
& \text{if }n\equiv 0\pmod{2} \\
\stackrel{F_1}{\lmt} 
& B^\mu & \equiv 
& \ov{\mu}+q(3,1)+q^2(2^2)+q^3(2,1^2) 
& \text{if }n\equiv 1\pmod{2} \\
\stackrel{F_0}{\lmt}
&&& \ov{\mu} \\
\end{array} $$

\abs
For $\ov{\mu}:=(3)$ this for $n\geq 7$ yields:
$$ \begin{array}{clcll}
& B^\mu & \equiv
& \ov{\mu}+q(1^3)
& \text{if }n\equiv 1\pmod{2} \\
\stackrel{F_0}{\lmt} 
& B^\mu & \equiv
& \ov{\mu}+q(1^3)+q(4)+q^2(3,1) \\
&&& \rule{0.6em}{0em} +q^2(2,1^2)+q^3(1^4)
& \text{if }n\equiv 0\pmod{2} \\
\stackrel{F_1}{\lmt}
&&& \ov{\mu}+q(1^3) \\
\end{array} $$

\abs
For $\ov{\mu}:=(2)$ this for $n\geq 7$ yields:
$$ \begin{array}{clcll}
& \hspace*{-1.5em} B^\mu & \equiv 
& \ov{\mu}+q(1^2)+q^2(2,1)+q(2,1^2)+q^2(1^4) 
& \hspace*{-1em}\text{if }n\equiv 1\pmod{2} \\
\stackrel{F_1}{\lmt} 
& \hspace*{-1.5em} B^\mu & \equiv 
& \ov{\mu}+q(1^2)+q(3)+q^2(1^3)+q(3,1) \\
&&& \rule{0.6em}{0em} +q^2(2^2)+(q+q^3)(2,1^2)+q^2(1^4)
& \hspace*{-1em}\text{if }n\equiv 0\pmod{2} \\
\stackrel{F_0}{\lmt}
&&& \ov{\mu}+q(1^2)+(1+q^2)(2,1)+(4)+2q(3,1) \\
&&& \rule{0.6em}{0em} +q^2(2^2)+(2q+q^3)(2,1^2)+2q^2(1^4) \\
\stackrel{B^{(2,1)},B^{(4)}}{\lmt}
&&& \ov{\mu}+q(1^2)+q^2(2,1)+q(2,1^2)+q^2(1^4) \\ 
\end{array} $$

\abs
For $\ov{\mu}:=(1)$ this for $n\geq 5$ yields:
$$ \begin{array}{clcll}
& B^\mu & \equiv 
& \ov{\mu}+q(1^3)
& \text{if }n\equiv 1\pmod{2} \\
\stackrel{F_0}{\lmt} 
& B^\mu & \equiv 
& \ov{\mu}+q(2)+q^2(1^2)+q(1^3) \\
&&& \rule{0.6em}{0em} +q^2(2,1^2)+q^3(1^4)  
& \text{if }n\equiv 0\pmod{2} \\
\stackrel{F_1}{\lmt}
&&& \ov{\mu}+(3)+2q(1^3) \\
\stackrel{B^{(3)}}{\lmt}
&&& \ov{\mu}+q(1^3) \\ 
\end{array} $$

\abs
For $\ov{\mu}:=()$ this for $n\geq 5$ yields:
$$ \begin{array}{clcll}
& B^\mu & \equiv 
& \ov{\mu}+q(1^2)+q^2(1^4)
& \text{if }n\equiv 1\pmod{2} \\
\stackrel{F_1}{\lmt} 
& B^\mu & \equiv 
& \ov{\mu}+q(1)+q(1^2)+q^2(1^3)+q^2(1^4) 
& \text{if }n\equiv 0\pmod{2} \\
\stackrel{F_0}{\lmt}
&&& \ov{\mu}+(2)+2q(1^2)+q^2(2,1) \\
&&& \rule{0.6em}{0em} +q(2,1^2)+2q^2(1^4) \\
\stackrel{B^{(2)}}{\lmt}
&&& \ov{\mu}+q(1^2)+q^2(1^4) \\
\end{array} $$

\absb
In conclusion, this yields the crystallized decomposition matrix
$\cD_n^q(\leq 4)$, for $n\geq 8$, as exhibited in Table \ref{dec2llt}.
Again, for $i\in\{0,1\}$ running through the residue classes modulo $2$,
we use the Kronecker type notation
$$ \dt=\dt_i(n):=\left\{\begin{array}{ll} 
1, &\text{if }n\equiv i\pmod{2}, \\
0, & \text{otherwise}. \\
\end{array}\right. $$
We remark that for $n\leq 13$ these results are also contained
in the explicit crystallized decomposition matrices given in 
\cite[Sect.10.3]{LasLecThi}.

\begin{table}\caption{Crystallized decomposition matrix for $\cP_n(\leq 4)$.} 
\label{dec2llt}
$$ \begin{array}{|l||r|r|r|rr|rr|}
\hline \rule{0em}{3em}
    & \spc \begin{rotate}{90} $()$ \end{rotate}
    & \spc \begin{rotate}{90} $(1)$ \end{rotate}
    & \spc \begin{rotate}{90} $(2)$ \end{rotate}
    & \spc \begin{rotate}{90} $(3)$ \end{rotate}
    & \spc \begin{rotate}{90} $(2,1)$ \end{rotate}
    & \spc \begin{rotate}{90} $(4)$ \end{rotate}
    & \spc \begin{rotate}{90} $(3,1)$ \end{rotate} \\
\hline
\hline
() \vspc &
1 & & & & & & \\
\hline
(1) \vspc &
q\dt_0 & 1 & & & & & \\
\hline
(2) \vspc &
. & q\dt_0 & 1 & & & & \\
(1^2) \vspc &
q & q^2\dt_0 & q & & & & \\
\hline
(3) \vspc &
. & . & q\dt_0 & 1 & & & \\
(2,1) \vspc &
. & . & q^2\dt_1 & . & 1 & & \\
(1^3) \vspc &
q^2\dt_0 &  q & q^2\dt_0 & q & . & & \\
\hline
(4) \vspc &
. & . & . & q\dt_0 & . & 1 & \\
(3,1) \vspc &
. & . & q\dt_0 & q^2\dt_0 & q\dt_1 & q & 1 \\
(2^2) \vspc &
. & . & q^2\dt_0 & . & q^2\dt_1 & . & q \\
(2,1^2) \vspc &
. & q^2\dt_0 & q + q^3\dt_0 & q^2\dt_0 & q^3\dt_1 & q & q^2 \\
(1^4) \vspc &
q^2 & q^3\dt_0 & q^2 & q^3\dt_0 & . & q^2 & . \\
\hline
\end{array} $$
\abs\hrulefill
\end{table}

\AbsT{Decomposition matrices.}\label{decmat2}
We proceed to determine the decomposition matrix $\cD_n(\leq 4)$, and
at the same time the adjustment matrix $\cA_n(\leq 4)$.
The results are given in Tables \ref{dec2} and \ref{dec2adj},
respectively, where we assume that $n\geq 8$.
Here, we again use the Kronecker type notation 
introduced in \ref{crystdec2}, and let
$$ \al=\al(n):=\left\{\begin{array}{ll} 
1, &\text{if }n\equiv 1,2\pmod{4}, \\
0, & \text{otherwise}, \\
\end{array}\right. $$
and
$$ \bt=\bt(n):=\left\{\begin{array}{ll} 
1, &\text{if }n\equiv 3,4,5,6\pmod{8}, \\
0, & \text{otherwise}. \\
\end{array}\right. $$

\abs
The decomposition matrices $\cD_n$ for $n\leq 7$ are known,
and for example given in \cite[p.137]{JamLN}, or accessible 
in the databases mentioned in Section \ref{intro}.
We again denote the Specht module associated with 
$\lb=[\lb_1,\ov{\lb}]\in\cP_n$ by $S_n(\ov{\lb})$, and 
the irreducible module associated with 
$\mu=[\mu_1,\ov{\mu}]\in\cP_n^{2\text{-reg}}$ by $D_n(\ov{\mu})$.


\begin{table}\caption{Decomposition matrix for $\cP_n(\leq 4)$.}\label{dec2}
$$ \begin{array}{|l||r|r|r|rr|rr|}
\hline \rule{0em}{3em}
    & \spc \begin{rotate}{90} $()$ \end{rotate}
    & \spc \begin{rotate}{90} $(1)$ \end{rotate}
    & \spc \begin{rotate}{90} $(2)$ \end{rotate}
    & \spc \begin{rotate}{90} $(3)$ \end{rotate}
    & \spc \begin{rotate}{90} $(2,1)$ \end{rotate}
    & \spc \begin{rotate}{90} $(4)$ \end{rotate}
    & \spc \begin{rotate}{90} $(3,1)$ \end{rotate} \\
\hline
\hline
() \vspc &
1 & & & & & & \\
\hline
(1) \vspc &
\dt_0 & 1 & & & & & \\
\hline
(2) \vspc &
\al & \dt_0 & 1 & & & & \\
(1^2) \vspc &
1+\al & \dt_0 & 1 & & & & \\
\hline
(3) \vspc &
\al\dt_0 & 1-\al & \dt_0 & 1 & & & \\
(2,1) \vspc &
\dt_1 & . & \dt_1 & . & 1 & & \\
(1^3) \vspc &
(1+\al)\dt_0 & 2-\al & \dt_0 & 1 & . & & \\
\hline
(4) \vspc &
\bt & (1\!-\!\al)\dt_0 & \al & \dt_0 & . & 1 & \\
(3,1) \vspc &
1\!+\!(1\!-\!\al)\dt_0\!+\!\al\dt_1\!+\!\bt & (1\!-\!\al)\dt_0 & \dt_0+\al 
                                            & \dt_0 & \dt_1 & 1 & 1 \\
(2^2) \vspc &
1+(1-\al)\dt_0+\al\dt_1 & . & \dt_0 & . & \dt_1 & . & 1 \\
(2,1^2) \vspc &
1+\dt_0+2\al\dt_1+\bt & (2\!-\!\al)\dt_0 & 1+\dt_0+\al 
                       & \dt_0 & \dt_1 & 1 & 1 \\
(1^4) \vspc &
1+\al+\bt & (2\!-\!\al)\dt_0 & 1+\al & \dt_0 & . & 1 & . \\
\hline
\end{array} $$
\abs\hrulefill
\end{table}

\begin{table}\caption{Adjustment matrix for $\cP_n(\leq 4)$.}\label{dec2adj}
$$ \begin{array}{|l||r|r|r|rr|rr|}
\hline \rule{0em}{3em}
    & \spc \begin{rotate}{90} $()$ \end{rotate}
    & \spc \begin{rotate}{90} $(1)$ \end{rotate}
    & \spc \begin{rotate}{90} $(2)$ \end{rotate}
    & \spc \begin{rotate}{90} $(3)$ \end{rotate}
    & \spc \begin{rotate}{90} $(2,1)$ \end{rotate}
    & \spc \begin{rotate}{90} $(4)$ \end{rotate}
    & \spc \begin{rotate}{90} $(3,1)$ \end{rotate} \\
\hline
\hline
() \vspc      &
1 & & & & & & \\
\hline
(1) \vspc     &
. & 1 & & & & & \\
\hline
(2) \vspc     &
\al & . & 1 & & & & \\
\hline
(3) \vspc     &
. & 1-\al & . & 1 & & & \\
(2,1) \vspc   &
(1-\al)\dt_1 & . & . & . & 1 & & \\
\hline
(4) \vspc     &
\bt & . & \al & . & . & 1 & \\
(3,1) \vspc   &
(2-2\al)\dt_0+2\al\dt_1 & . & . & . & . & . & 1 \\
\hline
\end{array} $$
\abs\hrulefill
\end{table}

\absb
The decomposition numbers of the Specht modules $S^{[n-m,m]}=S_n(m)$, 
where $m\in\{0\ld\fl{\frac{n}{2}}\}$,
are known by \cite{Jam2}, see also \cite[Thm.24.15]{JamLN}.

\absb
Whenever $\lb\in\cP_n(m)$ and $\mu\in\cP_n^{2\text{-reg}}(m)$ for 
some $m\in\{0\ld n-1\}$, 
the principle of first row removal, see \cite{JamIII}, yields 
$[S_n(\ov{\lb})\cn D_n(\ov{\mu})]=[S^{\ov{\lb}}\cn D^{\ov{\mu}}]$,
where the latter are decomposition numbers of $\cS_m$, which
since $m\leq 4$ are easily determined or can be
looked up in \cite[p.137]{JamLN}.

\absb
In order to complete $\cD_n(\leq 2)$ 
it remains to determine the row belonging to the partition $(1^2)$.
But from $\cD_n(\leq 2)=\cD^\zt_n(\leq 2)\cdot\cA_n(\leq 2)$,
using the entries of $\cD_n(\leq 2)$ already known, we conclude that 
$\cA_n(\leq 2)$ actually coincides with the upper left-hand 
$(3\tm 3)$-sub-matrix of the matrix given in Table \ref{dec2adj}.
This in turn determines $\cD_n(\leq 2)$.

\absb
Similarly, from 
$\cD_n(\leq 3)=\cD^\zt_n(\leq 3)\cdot\cA_n(\leq 3)$,
using the entries of $\cD_n(\leq 3)$ already known, we conclude that 
the row of $\cA_n(\leq 3)$ belonging to
the partition $(3)$ is as given in Table \ref{dec2adj}.
This in turn completes the row of $\cD_n(\leq 3)$ belonging to $(1^3)$.
In order to complete $\cD_n(\leq 3)$ 
it remains to determine the row belonging to the partition $(2,1)$.
Considering $2$-contents, reflecting the distribution of modules into
$2$-blocks, we conclude that $[S_n(2,1)\cn D_n(1)]=0$, and
$[S_n(2,1)\cn D_n(\ov{\mu})]>0$, where $\ov{\mu}\in\{(2),()\}$,
only if $n\equiv 1\pmod{2}$.

\absb
We consider $[S_n(2,1)\cn D_n(2)]$ for $n\equiv 1\pmod{2}$:
Applying \ref{jamestrick} with $0$-removal yields the inequality
$$ \begin{array}{cl}
& [S_n(2,1)\cn D_n(2)]\cdot [D_n(2)\da_0\cn D_{n-1}(1)] \\
\leq & [S_{n-1}(2)\cn D_{n-1}(1)]+[S_{n-1}(1^2)\cn D_{n-1}(1)]=2. 
\rule{0em}{1.2em} \\ \end{array} $$
Since by \ref{modbranch} we have $[D_n(2)\da_0\cn D_{n-1}(1)]=2$, 
using Table \ref{dec2llt} implies
$1=[S_n(2,1)\cn D^\zt_n(2)]\leq [S_n(2,1)\cn D_n(2)]\leq 1$.

\absb
We consider $[S_n(2,1)\cn D_n(())]$ for $n\equiv 1\pmod{2}$:
Applying \ref{jamestrick} with $0$-removal, using $[S_n(2,1)\cn D_n(2)]=1$ 
yields the inequality
$$ \begin{array}{cl}
& [S_n(2,1)\cn D_n(())]
+[D_n(2)\da_0\cn D_{n-1}(())] \\
\leq & [S_{n-1}(2)\cn D_{n-1}(())]+[S_{n-1}(1^2)\cn D_{n-1}(())] 
\rule{0em}{1.2em} \\
=&1+2\cdot\al(n-1) \rule{0em}{1.2em} \\
=& \left\{\begin{array}{ll}
1, & \text{if }n\equiv 1\pmod{4}, \\
3, & \text{if }n\equiv 3\pmod{4}. \\
\end{array}\right. \rule{0em}{2em} \\ \end{array} $$
Moreover, if $n\equiv 3\pmod{4}$ then we have 
$[D_n(2)\da_0]=[S_n(2)\da_0]=[S_{n-1}(2)]+[S_{n-1}(1)]$, thus
$$ [D_n(2)\da_0\cn D_{n-1}(())]
=[S_{n-1}(2)\cn D_{n-1}(())]+[S_{n-1}(1)\cn D_{n-1}(())]=2 .$$
Hence in any case we have $[S_n(2,1)\cn D_n(())]\leq 1$.
On the other hand, the hook length formula \cite{FraRobThr},
see also \cite[Thm.20.1]{JamLN}, readily yields
$\dim_{\Q}(S_n(2,1))=\frac{1}{3}n(n-2)(n-4)$, which is odd.
Since all irreducible $2$-modular representations of 
$\cS_n$ are self-contragredient, these all have even degree,
apart from the trivial module; see also \cite[Thm.11.8]{JamLN}.
Thus $[S_n(2,1)\cn D_n(())]$ is odd, 
hence we infer $[S_n(2,1)\cn D_n(())]=1$, completing $\cD_n(\leq 3)$.
We turn our attention to $\cD_n(\leq 4)$:

\absb
Whenever $\lb\in\cP_n(4)$ and $\mu\in\cP_n^{2\text{-reg}}(3)$
such that $l(\lb)=l(\mu)=3$, the principle of first column removal, 
see \cite{JamIII}, this time yields 
$[S_n(\ov{\lb})\cn D_n(\ov{\mu})]
=[S_{n-3}(\lb_2-1,\lb_3-1)\cn D_{n-3}(\mu_2-1,\mu_3-1)]$.
This settles the cases
$$ \begin{array}{lcl}
\mbox{}[S_n(3,1)\cn D_n(2,1)]&=&[S_{n-3}(2)\cn D_{n-3}(1)], \\
\mbox{}[S_n(2^2)\cn D_n(2,1)]&=&[S_{n-3}(1^2)\cn D_{n-3}(1)], 
\rule{0em}{1.2em} \\ \end{array} $$
where the latter decomposition numbers belong to $\cD_{n-3}(\leq 3)$ 
which we have already dealt with above.
Next we consider the row of $\cD_n(\leq 4)$ 
belonging to the partition $(2^2)$:

\absb
We consider $[S_n(2^2)\cn D_n(3)]$: A consideration of $2$-contents shows
that we have $[S_n(2^2)\cn D_n(3)]>0$ only if $n\equiv 0\pmod{2}$.
In this case, applying \ref{jamestrick} with $0$-removal, 
$S_n(2^2)\da_0=\{0\}$ implies
$[S_n(2^2)\cn D_n(3)]\cdot [D_n(3)\da_0\cn D_{n-1}(3)]=0$.
Since by \ref{modbranch} we have $[D_n(3)\da_0\cn D_{n-1}(3)]=1$, 
we get $[S_n(2^2)\cn D_n(3)]=0$.

\absb
We consider $[S_n(2^2)\cn D_n(1)]$: Similarly,
a consideration of $2$-contents shows
that $[S_n(2^2)\cn D_n(1)]>0$ only if $n\equiv 0\pmod{2}$.
In this case, applying \ref{jamestrick} with $0$-removal, 
$S_n(2^2)\da_0=\{0\}$ implies
$[S_n(2^2)\cn D_n(1)]\cdot [D_n(1)\da_0\cn D_{n-1}(1)]=0$.
Since by \ref{modbranch} we have $[D_n(1)\da_0\cn D_{n-1}(1)]=1$, 
we get $[S_n(2^2)\cn D_n(1)]=0$.

\absb
We consider $[S_n(2^2)\cn D_n(2)]$: 
If $n\equiv 1\pmod{2}$,
then applying \ref{jamestrick} with $0$-removal yields the inequality
$$ [S_n(2^2)\cn D_n(2)]\cdot [D_n(2)\da_0\cn D_{n-1}(1)] 
\leq [S_{n-1}(2^2)\cn D_{n-1}(1)] .$$
Since by \ref{modbranch} we have $[D_n(2)\da_0\cn D_{n-1}(1)]=2$,
and we have already seen that $[S_{n-1}(2^2)\cn D_{n-1}(1)]=0$, 
we conclude that $[S_n(2^2)\cn D_n(2)]=0$ in this case.

\abs
Similarly, if $n\equiv 0\pmod{2}$, 
then applying \ref{jamestrick} with $1$-removal yields 
$$ \begin{array}{cl}
& [S_n(2^2)\cn D_n(2)]\cdot [D_n(2)\da_1\cn D_{n-1}(2)] \\
\leq & [S_{n-1}(2^2)\cn D_{n-1}(2)]+[S_{n-1}(2,1)\cn D_{n-1}(2)].
\rule{0em}{1.2em} \\ \end{array} $$
Since by \ref{modbranch} we have $[D_n(2)\da_1\cn D_{n-1}(2)]=1$,
and we have already seen that $[S_{n-1}(2^2)\cn D_{n-1}(2)]=0$, 
using Table \ref{dec2llt} we get that 
$1\leq [S_n(2^2)\cn D_n(2)]\leq 1$, thus
$[S_n(2^2)\cn D_n(2)]=1$ in this case.

\absb
We consider $[S_n(2^2)\cn D_n(())]$: Letting $x_n:=[S_n(2^2)\cn D_n(())]$,
by what we already know we have
$$[S_n(2^2)]=x_n[D_n(())]+\dt_0[D_n(2)]+\dt_1[D_n(2,1)]+[D_n(3,1)] .$$
Thus, if $n\equiv 1\pmod{2}$, then by \ref{modbranch} we have 
$[D_n(2,1)\da_0]=[D_{n-1}(2)]$ and $[D_n(3,1)\da_0]=[D_{n-1}(3,1)]$,
hence by $0$-removal we get 
$$ [S_{n-1}(2^2)]=[S_n(2^2)\da_0]
=x_n[D_{n-1}(())]+[D_{n-1}(2)]+[D_n(3,1)] ,$$
implying that $x_{n-1}=x_n$.
Next, it follows from $\cD_n(\leq 4)=\cD^\zt_n(\leq 4)\cdot\cA_n(\leq 4)$,
where $\cD^\zt_n(\leq 4)$ is obtained from the crystallized decomposition
matrix $\cD^q_n(\leq 4)$ given in Table \ref{dec2llt} by evaluating 
at $q=1$, that $[S_n(2^2)]=[S_n(3,1)]-[S_n(4)]$.
Hence, if $n\equiv 1\pmod{2}$, then \cite[Thm.7.1.(i),(ii)]{Jam2} says that
$$ x_n=\left\{\begin{array}{ll} 
2, &\text{if }n\equiv 1\pmod{4}, \\
1, &\text{if }n\equiv 3\pmod{4}. \\
\end{array}\right. $$

\absb
Thus we have completed the row of $\cD_n(\leq 4)$ 
belonging to the partition $(2^2)$.
To conclude we just observe that the results obtained so far show that 
the adjustment matrix $\cA_n(\leq 4)$ is as given in Table \ref{dec2adj},
which in turn completes $\cD_n(\leq 4)=\cD^\zt_n(\leq 4)\cdot\cA_n(\leq 4)$.

\AbsT{Degree formulae.}\label{deg2formulae}
We are now prepared to obtain degree formulae for the irreducible 
modular representations parameterized by $\cP_n^{2\text{-reg}}(\leq 4)$.
They are given in Tables \ref{deg2leq3} and \ref{deg2eq4}, where 
for $\mu\in\cP_n^{2\text{-reg}}$ we again let $d^\mu:=\dim_{\F_2}(D^\mu)$:

\abs
In view of the decomposition matrix $\cD_n(\leq 4)$ for $n\geq 8$ 
given in Table \ref{dec2}, and the known decomposition numbers 
for $n\leq 7$, see the comment at the beginning of \ref{decmat2},
these follow straightforwardly from the hook length formula 
for the dimension of Specht modules, 
see \cite{FraRobThr} and also \cite[Thm.20.1]{JamLN}.

\begin{table}\caption{Degree formulae for $\cP_n^{2\text{-reg}}(\leq 3)$.} 
             \label{deg2leq3}
$$ \begin{array}{|l|l|lr|}
\hline
\mu & d^\mu & \multicolumn{2}{l|}{\text{condition}} \\
\hline
\hline
[n] & 1 & & \vspc \\
\hline
[n-1,1]   & n-2 & n\equiv 0 & \pmod{2} \vspc \\ 
(n\geq 3) & n-1 & n\equiv 1 & \pmod{2} \vspc \\
\hline
[n-2,2]   & \frac{1}{2}(n^2-5n+2) & n\equiv 2 & \pmod{4} \vspc \\
(n\geq 4) & \frac{1}{2}(n^2-5n+4) & n\equiv 0 & \pmod{4} \vspc \\
          & \frac{1}{2}(n^2-3n-2) & n\equiv 1 & \pmod{4} \vspc \\
          & \frac{1}{2}(n^2-3n)   & n\equiv 3 & \pmod{4} \vspc \\
\hline
[n-3,3]   & \frac{1}{6}(n^3-9n^2+14n)    & n\equiv 0 & \pmod{4} \vspc \\
(n\geq 7) & \frac{1}{6}(n^3-9n^2+20n-12) & n\equiv 2 & \pmod{4} \vspc \\
          & \frac{1}{6}(n^3-6n^2-n+6)    & n\equiv 3 & \pmod{4} \vspc \\
          & \frac{1}{6}(n^3-6n^2+5n)     & n\equiv 1 & \pmod{4} \vspc \\
\hline
[n-3,2,1] & \frac{1}{6}(2n^3-15n^2+25n-6) & n\equiv 3 & \pmod{4} \vspc \\
(n\geq 6) & \frac{1}{6}(2n^3-15n^2+25n)   & n\equiv 1 & \pmod{4} \vspc \\
          & \frac{1}{6}(2n^3-12n^2+16n)   & n\equiv 0 & \pmod{2} \vspc \\
\hline
\end{array} $$
\abs\hrulefill
\end{table}

\begin{table}\caption{Degree formulae for $\cP_n^{2\text{-reg}}(4)$.} 
             \label{deg2eq4}
$$ \begin{array}{|l|l|lr|}
\hline
\mu & d^\mu & \multicolumn{2}{l|}{\text{condition}} \\
\hline
\hline
[n-4,4]    & \frac{1}{24}(n^4-14n^3+47n^2-34n) 
           & n\equiv 6 & \pmod{8} \vspc \\
(n \geq 9) & \frac{1}{24}(n^4-14n^3+47n^2-34n+24) 
           & n\equiv 2 & \pmod{8} \vspc \\
           & \frac{1}{24}(n^4-14n^3+59n^2-94n+24) 
           & n\equiv 4 & \pmod{8} \vspc \\
           & \frac{1}{24}(n^4-14n^3+59n^2-94n+48) 
           & n\equiv 0 & \pmod{8} \vspc \\
           & \frac{1}{24}(n^4-10n^3+11n^2+22n) 
           & n\equiv 5 & \pmod{8} \vspc \\
           & \frac{1}{24}(n^4-10n^3+11n^2+22n+24) 
           & n\equiv 1 & \pmod{8} \vspc \\
           & \frac{1}{24}(n^4-10n^3+23n^2-14n-24) 
           & n\equiv 3 & \pmod{8} \vspc \\
           & \frac{1}{24}(n^4-10n^3+23n^2-14n) 
           & n\equiv 7 & \pmod{8} \vspc \\
\hline
[n-4,3,1] & \frac{1}{24}(2n^4-28n^3+118n^2-140n-48) 
          & n\equiv 1 & \pmod{4} \vspc \\
(n\geq 8) & \frac{1}{24}(2n^4-28n^3+118n^2-140n) 
          & n\equiv 3 & \pmod{4} \vspc \\
          & \frac{1}{24}(2n^4-20n^3+46n^2+20n-96) 
          & n\equiv 0 & \pmod{4} \vspc \\
          & \frac{1}{24}(2n^4-20n^3+46n^2+20n-48) 
          & n\equiv 2 & \pmod{4} \vspc \\
\hline
\end{array} $$
\abs\hrulefill
\end{table}

\section{James's Theorem revisited}\label{james}

\abs
Let $p$ be a rational prime. For $n\in\N_0$ and 
$\mu\in\cP_n^{p\text{-reg}}$ we let $d^\mu:=\dim_{\F_p}(D^\mu)$.
In \cite{JamMinimal}, a description of the growth behavior of $d^\mu$,
for $\mu\in\cP_n^{p\text{-reg}}(m)$, is given, where 
$m\in\N_0$ is fixed and $n\in\N_0$ is allowed to tend to infinity.
We are going to present two improvements of \cite[La.4]{JamMinimal},
the first one for arbitrary $p\geq 2$, while the second one takes
care of particular phenomena occurring in the case $p=2$.

\abs
Later on, we will need the following statement, 
which is an immediate consequence of {\cite[La.3]{JamMinimal} and its proof,
and is also contained in \cite[Thm.5.1]{JanSei} for $p\geq 3$.
We present another proof in the spirit of this work:

\PropT{\cite[La.3]{JamMinimal}}\label{twostepres}\mbox{}\\
Let $n\geq 2$ and $\mu\in\cP_n^{p\text{-reg}}$ such that the 
restriction $D^\mu\da_{\cS_{n-2}}$ of $D^\mu$ to $\cS_{n-2}$
is irreducible. Then we have $\mu=[n]$ or $\mu=[1^n]^R$.

\Pf

Let $b:=\mu_1$ be the largest part of $\mu$, 
where we may assume that $b\geq 2$.
By \ref{modbranch}, $\mu$ has a single normal node, $x$ say, which is
necessarily good. Since the $(a,b)$-node,
where $a\in\{1\ld p-1\}$ is the multiplicity of $b$ in $\mu$,
is the rightmost removable one and hence is normal, 
we infer that $x$ is the $(a,b)$-node.
Let $\mu^\ast:=\mu\smin\{x\}\in\cP_{n-1}^{p\text{-reg}}$ for short.
Then, by assumption, $\mu^\ast$ also has a single normal node, 
$x^\ast$ say, which also is necessarily good.
We distinguish two cases:

\absb
Let $a\geq 2$; note that this only occurs if $p\geq 3$.
Then $\mu^\ast$ has largest part $b$, occurring with
multiplicity $a-1$, and second largest part $b-1$, occurring with 
multiplicity $a'\in\{1\ld p-1\}$. By the same reasoning as above 
we infer that $x^\ast$ is the $(a-1,b)$-node. 
Hence the rightmost removable node $y$ of $\mu^\ast$ left of $x^\ast$,
the $(a+a'-1,b-1)$-node, is not normal. 
As $y$ and the addable node $x$ have distinct $p$-residues, we conclude 
that $y$ has the same $p$-residue as the addable $(1,b+1)$-node.
Thus we get $a+a'\equiv 0\pmod{p}$. Since $3\leq a+a'\leq 2p-2$ we 
conclude $a+a'=p$ and hence $y$ is the $(p-1,b-1)$-node.

\abs
Assume there is a further removable node $z$ of $\mu^\ast$ left of $y$,
and pick the rightmost one amongst those. Since $z$ is not normal, 
it has a $p$-residue different from that of the addable $(1,b+1)$-node 
and the removable node $y$. Thus it has the same
$p$-residue as the addable node $x$. But $z$ also is a removable
node of $\mu$, contradicting the fact that $\mu$ has a unique normal node. 
Hence $\mu$ only has the parts $b$ and $b-1$,
with multiplicities $a$ and $p-a-1$, respectively, that is $\mu=[1^n]^R$.

\absb
Let $a=1$. Then $\mu^\ast$ has largest part $b-1$, occurring with
multiplicity $a'\in\{1\ld p-1\}$, thus $x^\ast$ is the $(a',b-1)$-node.
We again distinguish two cases:

\absc
Let $a'\geq 2$; note again that this only occurs if $p\geq 3$.
Then $x^\ast$ also is a removable node of $\mu$, and since it is
not normal, it has the same $p$-residue as the addable $(1,b+1)$-node.
Thus we get $a'\equiv -1\pmod{p}$, thus $a'=p-1$.

\abs
Assume there is a further removable node $z$ of $\mu^\ast$ 
left of $x^\ast$, and pick the rightmost one amongst those. 
Since $z$ is not normal, it has a $p$-residue different from 
that of the addable $(1,b+1)$-node and the removable node $x^\ast$.
Thus it has the same $p$-residue as the addable node $x$.
But $z$ also is a removable node of $\mu$, contradicting 
the fact that $\mu$ has a unique normal node.
Hence $\mu$ only has the parts $b$ and $b-1$,
with multiplicities $1$ and $p-2$, respectively, that is $\mu=[1^n]^R$.

\absc
Let $a'=1$. Assume there is a further removable node $z$ of $\mu^\ast$ 
left of $x^\ast$, and pick the rightmost one amongst those. Since $z$ is not
normal, it has the same $p$-residue as the addable node $x$.
But $z$ also is a removable node of $\mu$, contradicting the
fact that $\mu$ has a unique normal node.
Hence $\mu$ only has the part $b$, with multiplicity $1$, that is $\mu=[n]$.
\QED

\Thm\label{jamesthm}
Let $m\in\N_0$, let $n_0\in\N_0$ such that $n_0\geq 2m+5$ if $p\geq 3$, 
and let $f\cn\{n_0,n_0+1,\ldots\}\ra\R$ fulfilling
$2\cdot f(n)\geq\max\{f(n+1),f(n+2)\}$ for all $n\geq n_0$.
Moreover, for $n\in\N_0$ let 
$$ \cM_n(\geq m+1):=\{\mu\in\cP_n^{p\text{-reg}};
                 \mu,\mu^M\not\in\cP_n^{p\text{-reg}}(\leq m)\} .$$

\abs\it
Then we have $d^\mu\geq f(n)$ for all $n\geq n_0$ and $\mu\in\cM_n(\geq m+1)$,
provided this inequality is known to hold for all $\mu\in\cM_n(\geq m+1)$
such that
$$ n\in\{n_0,n_0+1\} \quad\text{or}\quad
\{\mu,\mu^M\}\cap\cP_n^{p\text{-reg}}(m+1)\neq\emp .$$

\Pf
We proceed by induction on $n\geq n_0$, where by the assumptions made
we may assume that $n\geq n_0+2$, and
$\mu\in\cP_n^{p\text{-reg}}(k)$ and $\mu^M\in\cP_n^{p\text{-reg}}(k')$,
where $k,k'\geq m+2$. We consider the restriction $D^\mu\da_{\cS_{n-1}}$
of $D^\mu$ to $\cS_{n-1}$; 
recall that $D^{\mu^M}\cong D^\mu\otm D^{[1^n]^R}$, where 
$D^{[1^n]^R}$ is the sign representation, and hence
$D^{\mu^M}\da_{\cS_{n-1}}\cong(D^\mu\da_{\cS_{n-1}})\otm D^{[1^{n-1}]^R}$.
We distinguish two cases:

\abs\bfa
Assume that $D^\mu$ restricts reducibly to $\cS_{n-1}$.
By \ref{modbranch} we conclude that $\mu$ has at least two normal nodes.
We again distinguish two cases:

\abs\bfi
There is an $i$-good node $x\in N_i(\mu)$, for some $i\in\{0\ld p-1\}$,
such that $r_i(\mu,x)>a_i(\mu,x)$. Hence by \ref{modbranch} we have 
$[D^\mu\da_i\cn D^{\mu\smin\{x\}}]\geq 2$, where,
since $x$ cannot possibly belong to the first row of $\mu$, we have 
$$ \mu\smin\{x\}\in\cP_{n-1}^{p\text{-reg}}(k-1) .$$
Moreover, by \cite[Alg.4.8]{KleIII}, there is an $i'$-good node 
$x'\in N_{i'}(\mu^M)$, where $i'\equiv -i\pmod{p}$ and 
$r_{i'}(\mu^M,x')-a_{i'}(\mu^M,x')=r_i(\mu,x)-a_i(\mu,x)$,
entailing that
$$ (\mu\smin\{x\})^M=\mu^M\smin\{x'\}\in\cP_{n-1}^{p\text{-reg}}(k'-1) .$$
Thus by induction we have 
$d^{\mu}\geq 2\cdot d^{\mu\smin\{x\}}\geq 2\cdot f(n-1)\geq f(n)$. 

\abs\bfii
There are an $i$-good node $x\in N_i(\mu)$ and a $j$-good node 
$y\in N_j(\mu)$, for some $i\neq j\in\{0\ld p-1\}$, such that 
$r_i(\mu,x)=a_i(\mu,x)$ and $r_j(\mu,y)=a_j(\mu,y)$.
Hence by \ref{modbranch} we have
$[D^\mu\da_i\cn D^{\mu\smin\{x\}}]=1$ and
$[D^\mu\da_j\cn D^{\mu\smin\{y\}}]=1$, where 
$$ \mu\smin\{x\},\,\mu\smin\{y\}\in
\cP_{n-1}^{p\text{-reg}}(k-1)\dcup\cP_{n-1}^{p\text{-reg}}(k) .$$
Moreover, there are an $i'$-good node $x'\in N_{i'}(\mu^M)$ and 
a $j'$-good node $y'\in N_{j'}(\mu^M)$,
where $i'\equiv -i\pmod{p}$ and $j'\equiv -j\pmod{p}$, and 
$r_{i'}(\mu^M,x')=a_{i'}(\mu^M,x')$ and $r_{j'}(\mu^M,y')=a_{j'}(\mu^M,y')$,
entailing that both
$$ (\mu\smin\{x\})^M=\mu^M\smin\{x'\}\in
\cP_{n-1}^{p\text{-reg}}(k'-1)\dcup\cP_{n-1}^{p\text{-reg}}(k') $$ 
and
$$ (\mu\smin\{y\})^M=\mu^M\smin\{y'\}\in
\cP_{n-1}^{p\text{-reg}}(k'-1)\dcup\cP_{n-1}^{p\text{-reg}}(k') .$$ 
Thus by induction we have 
$d^{\mu}\geq d^{\mu\smin\{x\}}+d^{\mu\smin\{y\}}\geq 2\cdot f(n-1)\geq f(n)$. 

\abs\bfb
Assume that $D^\mu$ restricts irreducibly to $\cS_{n-1}$.
By \ref{modbranch} we conclude that $\mu$ has a single normal node,
$x$ say, which hence is an $i$-good node for some $i\in\{0\ld p-1\}$,
and we have $r_i(\mu,x)=a_i(\mu,x)$.
Letting $\mu^\ast:=\mu\smin\{x\}$ for short, we have 
$D^\mu\da_{\cS_{n-1}}\cong D^{\mu^\ast}$, where for some $\eps\in\{0,1\}$
we have 
$$ \mu^\ast\in\cP_{n-1}^{p\text{-reg}}(k-\eps) .$$ 
Moreover, $\mu^M$ also has a single normal node, $x'$ say, 
which, by \cite[Alg.4.8]{KleIII}, is an $i'$-good node
where $i'\equiv -i\pmod{p}$ and $r_{i'}(\mu^M,x')=a_{i'}(\mu^M,x')$.
Letting $(\mu^M)^\ast:=\mu^M\smin\{x'\}$, for some $\eps'\in\{0,1\}$ we get
$$ (\mu^\ast)^M=(\mu^M)^\ast\in\cP_{n-1}^{p\text{-reg}}(k'-\eps') .$$

\abs
Assume that $D^{\mu^\ast}$ restricts irreducibly to $\cS_{n-2}$.
Then from \ref{twostepres} we conclude that $\mu=[n]$ or $\mu=[1^n]^R$.
In the first case we have $\mu\in\cP_n^{p\text{-reg}}(0)$,
while in the latter case 
we have $\mu=[1^n]^R=[n]^M$, hence $\mu^M\in\cP_n^{p\text{-reg}}(0)$.
Thus in both cases we arrive at a contradiction. Hence we may assume
that $D^{\mu^\ast}$ restricts reducibly to $\cS_{n-2}$.
We again distinguish two cases:

\abs\bfi
We have $k-\eps\geq m+2$ and $k'-\eps'\geq m+2$. 
Then, by what we have already seen in (a), 
$D^{\mu^\ast}\da_{\cS_{n-2}}$ has constituents
$D^{\mu^\ast\smin\{y\}}$ and $D^{\mu^\ast\smin\{z\}}$,
where $y$ and $z$ are good nodes of $\mu^\ast$,
which are possibly identical, and 
$$ \mu^\ast\smin\{y\},\,\mu^\ast\smin\{z\}\in
\cP_{n-2}^{p\text{-reg}}(k-\eps-1)\dcup\cP_{n-2}^{p\text{-reg}}(k-\eps) .$$
Moreover, there are good nodes $y'$ and $z'$ of $(\mu^\ast)^M$, 
which are possibly identical, such that
$$ (\mu^\ast\smin\{y\})^M=(\mu^\ast)^M\smin\{y'\}\in
\cP_{n-2}^{p\text{-reg}}(k'-\eps'-1)\dcup\cP_{n-2}^{p\text{-reg}}(k'-\eps') $$ 
and 
$$ (\mu^\ast\smin\{z\})^M=(\mu^\ast)^M\smin\{z'\}\in
\cP_{n-2}^{p\text{-reg}}(k'-\eps'-1)\dcup\cP_{n-2}^{p\text{-reg}}(k'-\eps') .$$
Thus by induction we have 
$d^{\mu}=d^{\mu^\ast}\geq 
 d^{\mu^\ast\smin\{y\}}+d^{\mu^\ast\smin\{z\}}\geq 2\cdot f(n-2)\geq f(n)$. 

\abs\bfii
We have $k-\eps=m+1$ or $k'-\eps'=m+1$. Hence we may assume that 
$$ \mu\in\cP_n^{p\text{-reg}}(m+2) \quad\text{and}\quad 
\mu^\ast=\mu\smin\{x\}\in\cP_{n-1}^{p\text{-reg}}(m+1) .$$
Thus $x$, being the rightmost removable node of $\mu$,
does not belong to the first row of $\mu$. Hence the rightmost
node of the first row is not removable, that is for the first two parts
of $\mu$ we have $\mu_1=\mu_2$. This cannot happen if $p=2$, thus
we have $p\geq 3$. Moreover, we have $n-m-2=\mu_1=\mu_2\leq m+2$,
or equivalently $n\leq 2m+4$. Hence, since we are assuming that 
$n\geq 2m+5$, this case cannot happen either.
\QED

\Rem\label{jamesrem}
The statement of \ref{jamesthm}, and the strategy of proof, are 
reminiscent of \cite[La.4]{JamMinimal}. But while the proof there 
employs the ordinary branching rule, see for example 
\cite[Thm.9.2]{JamLN}, 
here we make use of the modular branching rule \ref{modbranch}, which
of course had not been available at the time of writing of \cite{JamMinimal}.

\abs
The statement of \cite[La.4]{JamMinimal} is similar to the
one given here, but our condition 
`$\{\mu,\mu^M\}\cap\cP_n^{p\text{-reg}}(m+1)\neq\emp$' 
is replaced there by the stronger one
$$\text{`}\{\mu,\mu^M\}\cap
  \left(\cP_n^{p\text{-reg}}(m+1)\dcup\cP_n^{p\text{-reg}}(m+2)\right)
  \neq\emp\text{'} .$$ 
As it turns out in practice, only having to check a weaker condition 
saves quite a bit of explicit computation, and lends itself to a 
treatment independent of $p$, see \ref{degbnd}. In particular,
for the case $p=3$ and $m=3$ this leads to a shorter proof of
\cite[Prop.3.1]{Dan}, together with a smaller lower bound $n_0$.

\Rem\label{seitzrem}
The assumptions of case (b) of the proof of
\ref{jamesthm} say that the restriction $D^\mu\da_{\cS_{n-1}}$ of 
$D^\mu$ to $\cS_{n-1}$ is irreducible, hence in theses cases $\mu$ is a 
Jantzen-Seitz partition \cite{JanSei}, see \cite{For,KleI,Kle0}.

\abs
Case (b)(ii) of the proof of \ref{jamesthm} shows that for $p\geq 3$
we can weaken the conditions slightly more as follows:
Note that $\cM_n(\geq m+1)=\emp$ for $n\leq m+1$ anyway. Then, 
leaving out the condition `$n_0\geq 2m+5$', we assume
instead that the desired inequality holds for all 
$n\in\{m+3\ld 2m+4\}$ and all Jantzen-Seitz partitions $\mu\in\cM_n(\geq m+2)$
such that $\mu\in\cP_n^{p\text{-reg}}(m+2)$ and $\mu_1=\mu_2$. 

\abs
For fixed $m\in\N_0$, the latter condition leads to a finite
set of cases to be checked additionally. Actually, for small $m$
this set turns out to be very small: For any $m\in\{0\ld 6\}$ there is 
just the single case $n=2m+4$ and $\mu=[\frac{n}{2},\frac{n}{2}]=[m+2,m+2]$
if $p\geq 5$, and no case at all for $p=3$.

\abs
Unfortunately, this pattern does not continue: 
For $m=7$ and $p=5$, next to $\mu=[9^2]\in\cP_{18}^{5\text{-reg}}(9)$,
we get $\mu=[8^2,1]\in\cP_{17}^{5\text{-reg}}(9)$ and
$\mu=[7^2,1^2]\in\cP_{16}^{5\text{-reg}}(9)$, 
and there are many more examples for larger $m$.
But we have not been able to find any example for $p=3$, 
so that we are wondering whether in this case there are any at all.
We are tempted to ask for a classification of the 
Jantzen-Seitz partitions $\mu\in\cP_n^{p\text{-reg}}$ such that
$\mu_1=\mu_2$. To our knowledge this is not available in the literature, 
and we leave it as an open question to the reader.
\QED

\abs\abs
Although the above theorem also holds for the case $p=2$,
it would not be strong enough for our purposes, the reason being
the existence of `very small' representations escaping the desired
growth behavior in low degrees. These we consider next, in order to
proceed to prove an improved theorem. 

\AbsT{Basic spin representations.}\label{basicspin}
Let $p:=2$. 

\absb
For $n\in\N$ let $\mu_{\text{bs}}(n)\in\cP_n^{2\text{-reg}}$ 
be defined by
$$ \mu_{\text{bs}}(n):=\left\{\begin{array}{ll}
\mbox{}[\frac{n+2}{2},\frac{n-2}{2}], &\text{if }n\equiv 0\pmod{2}, \\
\mbox{}[\frac{n+1}{2},\frac{n-1}{2}], &\text{if }n\equiv 1\pmod{2}. \\
\end{array}\right. $$
The associated irreducible $\F_2\cS_n$-module 
$D_n(\text{bs})$, 
being called the basic spin module, by \cite[Cor.5.7]{Jam2}, 
see also \cite[La.5.3]{Ben} or \cite[Tbl.III]{Wal}, has dimension 
$$ \dim_{\F_2}(D_n(\text{bs}))=2^{\fl{\frac{n-1}{2}}}
=\left\{\begin{array}{ll}
2^{\frac{n-2}{2}}, &\text{if }n\equiv 0\pmod{2}, \\
2^{\frac{n-1}{2}}, &\text{if }n\equiv 1\pmod{2}. \\
\end{array}\right. $$
For later use, in Table \ref{bspin} we record the results obtained 
from removing $i$-good and adding $i$-cogood nodes to 
$\mu_{\text{bs}}(n)$, for $n\geq 3$. 
Here, the rows are indexed by the congruence classes of $n$ modulo $4$,
and missing entries indicate the non-existence of $i$-normal and 
$i$-conormal nodes, respectively.
Using this, the degree formulae and \ref{modbranch},
for $n\geq 2$ we infer
$$ \begin{array}{ccccrl}
D_n(\text{bs})\da_{\cS_{n-1}}
& \cong & D_n(\text{bs})\da_0
& \cong & D_{n-1}(\text{bs}), & \text{if }n\equiv 0\pmod{4}, \\
\mbox{}[D_n(\text{bs})\da_{\cS_{n-1}}]
& = & [D_n(\text{bs})\da_0]
& = & 2\cdot [D_{n-1}(\text{bs})], & \text{if }n\equiv 1\pmod{4}, \\
D_n(\text{bs})\da_{\cS_{n-1}}
& \cong & D_n(\text{bs})\da_1
& \cong & D_{n-1}(\text{bs}), & \text{if }n\equiv 2\pmod{4}, \\
\mbox{}[D_n(\text{bs})\da_{\cS_{n-1}}]
& = & [D_n(\text{bs})\da_1]
& = & 2\cdot [D_{n-1}(\text{bs})], & \text{if }n\equiv 3\pmod{4}. \\
\end{array} $$

\absb
For $n\geq 6$ let $\mu_{\text{bbs}}(n)\in\cP_n^{2\text{-reg}}$ 
be defined by
$$ \mu_{\text{bbs}}(n):=\left\{\begin{array}{ll}
\mbox{}[\frac{n}{2},\frac{n-2}{2},1], &\text{if }n\equiv 0\pmod{2}, \\
\mbox{}[\frac{n+1}{2},\frac{n-3}{2},1], &\text{if }n\equiv 1\pmod{2}. \\
\end{array}\right. $$
The associated irreducible $\F_2\cS_n$-module 
$D_n(\text{bbs})$, 
being called the second basic spin module, by
\cite[Thm.1.2]{Ben} and \cite[Tbl.IV]{Wal}, 
see also \cite[Thm.7.1]{Jam2} or \cite[Cor.4.2]{Dan},
has dimension 
$$ \dim_{\F_2}(D_n(\text{bbs}))
=\left\{\begin{array}{ll}
(n-3)\cdot 2^{\frac{n-2}{2}}, &\text{if }n\equiv 0\pmod{4}, \\
(n-4)\cdot 2^{\frac{n-3}{2}}, &\text{if }n\equiv 1\pmod{4}, \\
(n-2)\cdot 2^{\frac{n-2}{2}}, &\text{if }n\equiv 2\pmod{4}, \\
(n-2)\cdot 2^{\frac{n-3}{2}}, &\text{if }n\equiv 3\pmod{4}. \\
\end{array}\right. $$
Again, in Table \ref{bspin} we record the results obtained 
from removing $i$-good and adding $i$-cogood nodes to 
$\mu_{\text{bbs}}(n)$, for $n\geq 6$.
Using this, the degree formulae and \ref{modbranch},
for $n\geq 7$ we infer
$$ \begin{array}{ccrl}
\mbox{}[D_n(\text{bbs})\da_1]
& = & 2\cdot [D_{n-1}(\text{bbs})], & \text{if }n\equiv 0\pmod{4}. \\
D_n(\text{bbs})\da_0 
& \cong & D_{n-1}(\text{bbs}), & \text{if }n\equiv 1\pmod{4}, \\
\mbox{}[D_n(\text{bbs})\da_0]
& = & 3\cdot [D_{n-1}(\text{bs})]+2\cdot [D_{n-1}(\text{bbs})],
& \text{if }n\equiv 2\pmod{4}, \\
D_n(\text{bbs})\da_{\cS_{n-1}}
& \cong & D_{n-1}(\text{bs})\oplus D_{n-1}(\text{bbs}),
& \text{if }n\equiv 3\pmod{4}, \\
\end{array} $$
where we only record non-zero $i$-restrictions, and for
$n\equiv 3\pmod{4}$ the direct summands refer to $0$- and
$1$-restriction, respectively; for $n=6$ we get
$[D_6(\text{bbs})\da_0]= 3\cdot [D_5(\text{bs})]+4\cdot [D([5])]$.

\absb
We conclude from Table \ref{bspin} that adding an $i$-cogood node to
$\mu_{\text{bs}}(n)$, we end up with partitions in 
$\{\mu_{\text{bs}}(n+1),\mu_{\text{bbs}}(n+1)\}$, and repeating this step, 
possibly with a different $i$, we end up with partitions in
$\{\mu_{\text{bs}}(n+2),\mu_{\text{bbs}}(n+2)\}$.

\begin{table}\caption{Removal and addition of nodes.}\label{bspin}
$$ \begin{array}{|r||l|l||l|l|}
\hline
\multicolumn{5}{|c|}{\mu_{\text{bs}}(n)} \\
\hline
n & 0\text{-good} & 1\text{-good} & 0\text{-cogood} & 1\text{-cogood} \\
\hline
\hline
0 & \mu_{\text{bs}}(n-1) & & \mu_{\text{bs}}(n+1) & \\
1 & \mu_{\text{bs}}(n-1) & & \mu_{\text{bbs}}(n+1) & \mu_{\text{bs}}(n+1) \\
2 & & \mu_{\text{bs}}(n-1) & \mu_{\text{bbs}}(n+1) & \mu_{\text{bs}}(n+1) \\
3 & & \mu_{\text{bs}}(n-1) & \mu_{\text{bs}}(n+1) & \\
\hline
\end{array} $$ 
$$ \begin{array}{|r||l|l||l|l|}
\hline
\multicolumn{5}{|c|}{\mu_{\text{bbs}}(n)} \\
\hline
n & 0\text{-good} & 1\text{-good} & 0\text{-cogood} & 1\text{-cogood} \\
\hline
\hline
0 & & \mu_{\text{bbs}}(n-1) 
  & \mu_{\text{bbs}}(n+1) & [\frac{n}{2},\frac{n-2}{2},2] \\
1 & \mu_{\text{bbs}}(n-1) & & & [\frac{n+1}{2},\frac{n-3}{2},2] \\
2 & \mu_{\text{bs}}(n-1) & & & \mu_{\text{bbs}}(n+1) \\
3 & \mu_{\text{bs}}(n-1) & \mu_{\text{bbs}}(n-1) 
  & & \mu_{\text{bbs}}(n+1) \\
\hline
\end{array} $$

\abs\hrulefill
\end{table}

\Thm\label{james2thm}
Let $p:=2$, let $m\in\N_0$, let $n_0\geq 5$, let
$f\cn\{n_0,n_0+1,\ldots\}\ra\R$ fulfilling
$2\cdot f(n)\geq\max\{f(n+1),f(n+2)\}$ for all $n\geq n_0$,
and for $n\in\N$ let
$$ \cM'_n(\geq m+1):=\cP_n^{p\text{-reg}}\smin
  (\cP_n^{p\text{-reg}}(\leq m)\cup\{\mu_{\text{bs}}(n)\}) .$$

\abs\it
Then we have $d^\mu\geq f(n)$ for all $n\geq n_0$ and
$\mu\in\cM'_n(\geq m+1)$, provided this inequality is known to hold
for all $\mu\in\cM'_n(\geq m+1)$ such that
$$ n\in\{n_0,n_0+1\} \quad\text{or}\quad \mu\in\cP_n^{p\text{-reg}}(m+1) .$$

\Pf
We again proceed by induction on $n\geq n_0$, 
where by the assumptions made we may assume that $n\geq n_0+2\geq 7$, 
and $\mu\in\cP_n^{2\text{-reg}}(k)$ where $k\geq m+2$, and are again 
going to consider the restriction of $D^\mu$ to $\cS_{n-1}$. 
We first assume $\mu\neq\mu_{\text{bbs}}(n)$, and distinguish two cases:

\abs\bfa
Assume that $D^\mu$ restricts reducibly to $\cS_{n-1}$.
By \ref{modbranch} we conclude that $\mu$ has at least two normal nodes.
We again distinguish two cases:

\abs\bfi
There is an $i$-good node $x\in N_i(\mu)$, for some $i\in\{0,1\}$,
such that $r_i(\mu,x)>a_i(\mu,x)$. Hence by \ref{modbranch} we have 
$[D^\mu\da_i\cn D^{\mu\smin\{x\}}]\geq 2$, where,
since $x$ cannot possibly belong to the first row of $\mu$, we have 
$$ \mu\smin\{x\}\in\cP_{n-1}^{2\text{-reg}}(k-1) .$$
Since from Table \ref{bspin} we conclude that 
$\mu\smin\{x\}\neq\mu_{\text{bs}}(n-1)$, by induction we have 
$d^{\mu}\geq 2\cdot d^{\mu\smin\{x\}}\geq 2\cdot f(n-1)\geq f(n)$. 

\abs\bfii
There are an $i$-good node $x\in N_i(\mu)$ and a $j$-good node 
$y\in N_j(\mu)$, for some $i\neq j\in\{0,1\}$, such that 
$r_i(\mu,x)=a_i(\mu,x)$ and $r_j(\mu,y)=a_j(\mu,y)$.
Hence by \ref{modbranch} we have
$[D^\mu\da_i\cn D^{\mu\smin\{x\}}]=1$ and
$[D^\mu\da_j\cn D^{\mu\smin\{y\}}]=1$, where both
$$ \mu\smin\{x\},\,\mu\smin\{y\}\in
\cP_{n-1}^{2\text{-reg}}(k-1)\dcup\cP_{n-1}^{2\text{-reg}}(k) .$$
Since from Table \ref{bspin} we conclude that
$\mu\smin\{x\}\neq\mu_{\text{bs}}(n-1)\neq\mu\smin\{y\}$, 
by induction we have 
$d^{\mu}\geq d^{\mu\smin\{x\}}+d^{\mu\smin\{y\}}\geq 2\cdot f(n-1)\geq f(n)$. 

\abs\bfb
Assume that $D^\mu$ restricts irreducibly to $\cS_{n-1}$.
By \ref{modbranch} we conclude that $\mu$ has a single normal node,
$x\in N_i(\mu)$ for some $i\in\{0,1\}$, which hence is $i$-good 
and we have $r_i(\mu,x)=a_i(\mu,x)$.
Since $\mu$ is $2$-regular, for the first two parts of $\mu$ 
we infer $\mu_1\neq\mu_2$, implying that the rightmost node
of the first row is removable, and hence coincides with the normal node $x$.
Letting $\mu^\ast:=\mu\smin\{x\}$ for short,
we have $D^\mu\da_{\cS_{n-1}}\cong D^{\mu^\ast}$, where
$\mu^\ast\in\cP_{n-1}^{2\text{-reg}}(k)$.

\abs
Next, in view of Table \ref{bspin}, assume that
$n\equiv 1\pmod{4}$ and $\mu=[\frac{n-1}{2},\frac{n-3}{2},2]$. 
Then $\mu$ has three $1$-normal nodes, a contradiction.
Similarly, assume that $n\equiv 2\pmod{4}$ and
$\mu=[\frac{n}{2},\frac{n-4}{2},2]$.
Then $\mu$ has both a $0$-normal and a $1$-normal node,
a contradiction as well.
Thus, since
$\mu\not\in\{\mu_{\text{bs}}(n),\mu_{\text{bbs}}(n)\}$,
from Table \ref{bspin} we conclude that
$\mu^\ast\not\in\{\mu_{\text{bs}}(n-1),\mu_{\text{bbs}}(n-1)\}$.

\abs
Assume that $D^{\mu^\ast}$ restricts irreducibly to $\cS_{n-2}$. 
Then from \ref{twostepres} we conclude that 
$\mu=[n]\in\cP_n^{2\text{-reg}}(0)$,
a contradiction. Hence we may assume 
that $D^{\mu^\ast}$ restricts reducibly to $\cS_{n-2}$.
Thus, by what we have already seen in (a), 
$D^{\mu^\ast}\da_{\cS_{n-2}}$ has constituents
$D^{\mu^\ast\smin\{y\}}$ and $D^{\mu^\ast\smin\{z\}}$,
where $y$ and $z$ are good nodes of $\mu^\ast$,
which are possibly identical, and 
$$ \mu^\ast\smin\{y\},\,\mu^\ast\smin\{z\}\in
\cP_{n-2}^{2\text{-reg}}(k-1)\dcup\cP_{n-2}^{2\text{-reg}}(k) .$$
Thus by induction we have 
$d^{\mu}=d^{\mu^\ast}\geq 
 d^{\mu^\ast\smin\{y\}}+d^{\mu^\ast\smin\{z\}}\geq 2\cdot f(n-2)\geq f(n)$. 

\absb
Let now $\mu=\mu_{\text{bbs}}(n)$. We again distinguish two cases:

\abs\bfi
Let $n\equiv 0\pmod{2}$, hence $k=\frac{n}{2}$. 
Then, by \ref{basicspin},
$\mu$ has a normal node $x$, belonging to the second row of $\mu$, such that 
$$ \mu\smin\{x\}=\mu_{\text{bbs}}(n-1)\in\cP_{n-1}^{2\text{-reg}}(k-1) $$
and $[D_n(\text{bbs})\da_{\cS_{n-1}}\cn D_{n-1}(\text{bbs})]=2$.
Hence by induction we from this infer that
$d^{\mu}\geq 2\cdot d^{\mu\smin\{x\}}\geq 2\cdot f(n-1)\geq f(n)$. 

\abs\bfii
Let $n\equiv 1\pmod{2}$, such that $n\geq 9$, hence $k=\frac{n-1}{2}$.
Then, by \ref{basicspin}, $\mu$ has a good node $x$, 
belonging to the first row of $\mu$, such that 
$$ \mu\smin\{x\}=\mu_{\text{bbs}}(n-1)\in\cP_{n-1}^{2\text{-reg}}(k) $$
and $[D_n(\text{bbs})\da_{\cS_{n-1}}\cn D_{n-1}(\text{bbs})]=1$.
By what we have already seen in (i), we have
$[D_{n-1}(\text{bbs})\da_{\cS_{n-2}}\cn D_{n-2}(\text{bbs})]=2$,
where $\mu_{\text{bbs}}(n-2)\in\cP_{n-1}^{2\text{-reg}}(k-1)$.
Thus by induction we have 
$d^{\mu}\geq d^{\mu\smin\{x\}}\geq 2\cdot f(n-2)\geq f(n)$. 

\abs
Finally, for $n=7$, which occurs only if $n_0=5$, we have
$\mu=[4,2,1]$, hence $k=3$, implying $m\leq 1$.
We have $d^{\mu}=20$, as well as $d^{[6,1]}=6$ and $d^{[5,2]}=14$,
by \cite{Jam2}, see also \cite[Thm.24.15]{JamLN}, or by looking up 
the decomposition matrix $\cD_7$ in \cite[p.137]{JamLN}
or in the databases mentioned in Section \ref{intro}.
This implies $d^{\mu}\geq f(7)$. 
\QED

\Rem\label{danzrem}
Note that using the notation of \ref{jamesthm} we have
$\cM'_n(\geq m+1):=\cM_n(\geq m+1)\smin\{\mu_{\text{bs}}(n)\}$,
and that $\cM'_n(\geq m+1)=\emp$ for $n\leq 4$.

\abs
The statement of \ref{james2thm}, and the strategy of proof, are 
reminiscent of \cite[Thm.4.3]{Dan}, but there stronger conditions
on the growth behavior of $f$ are needed, and our condition 
`$\mu\in\cP_n^{p\text{-reg}}(m+1)$' is replaced there by the stronger one  
$$\text{`}\mu\in\cP_n^{p\text{-reg}}(m+1)\dcup\cP_n^{p\text{-reg}}(m+2)
  \text{'} .$$
Again, in practice, only having to check weaker conditions 
saves quite a bit of explicit computation, see \ref{degbnd}. In particular,
for the case $m=3$ this leads to a shorter proof of 
\cite[Prop.4.4]{Dan}, together with a smaller lower bound $n_0$.

\section{Explicit results}\label{results}

\abs
We are now prepared to tackle the problem of classifying all
irreducible representations $D^\mu$ of $\cS_n$, where 
$\mu\in\cP_n^{p\text{-reg}}$, satisfying $d^\mu:=\dim_{\F_p}(D^\mu)\leq n^3$.
The general strategy to obtain this classification is described in 
\ref{strategy}, and subsequently it is applied to compile explicit 
lists for $p\leq 7$. To start with, as an immediate corollary 
of the information collected above, we have the following:

\AbsT{Degree formulae.}\label{degformulae} 
The formulae for $d^\mu$, where $\mu\in\cP_n^{p\text{-reg}}(\leq 4)$,
given in \ref{deg3formulae} and \ref{deg2formulae} for $p=3$ and $p=2$, 
respectively, together with the results in \cite[La.1.21]{BruKle}
giving similar formulae for $p\geq 5$,
show that for all $\mu\in\cP_n^{p\text{-reg}}(m)$, where 
$m\in\{0\ld 4\}$ and $n\geq 2m$, independent of $p\geq 2$ we get
$$ d^\mu=\dim_{\F_p}(D^\mu)\geq\left\{\begin{array}{cl}
1   & \text{if }m=0, \\
n-2 & \text{if }m=1, \rule{0em}{1.2em} \\
\frac{n^2-5n+2}{2} & \text{if }m=2, \rule{0em}{1.2em} \\
\frac{n^3-9n^2+14n}{6} & \text{if }m=3, \rule{0em}{1.2em} \\
\frac{n^4-14n^3+47n^2-34n}{24} & \text{if }m=4. \rule{0em}{1.2em} \\
\end{array}\right. $$
These lower bounds are best possible in the sense that,
for each $m\in\{0\ld 4\}$, equality holds for infinitely many
values of $n\geq 2m$. We remark that the cases $m\leq 1$ are already
noted in \cite[Thm.6]{JamMinimal}, and that the cases $m\in\{2,3\}$
for $p\geq 5$ also appear in \cite[La.1.18, La.1.20]{BruKle};
the remaining lower bounds also follow from a close inspection
of the results given without proof in \cite[App.]{JamMinimal}.

\AbsT{Degree bounds.}\label{degbnd}
We are going to apply Theorems \ref{jamesthm} and \ref{james2thm}, 
for $p\geq 3$ and $p=2$, respectively, in the case $m=3$. 
Hence a natural choice for the function $f(n)$ should be closely
related to the lower bound function given in \ref{degformulae} 
for the case $m+1=4$; 
see also \cite[Prop.3.1]{Dan} and \cite[Prop.4.4]{Dan} for the
cases $p=3$ and $p=2$, respectively.
Thus we let
$$ f_4\cn\Z\ra\Z\cn n\mt\frac{n^4-14n^3+47n^2-34n}{24} .$$
It is easily seen that $f_4(n)$ has integral values,
which are positive for $n\geq 10$, that it is strictly 
increasing for $n\geq 8$, and fulfills $2\cdot f_4(n)\geq f_4(n+2)$
for $n\geq 16$. 

\absb
In order to obtain a strong bound $n_0$ in Theorems \ref{jamesthm} and 
\ref{james2thm}, we use a modification of the above function: 
Let $g\cn\{11,12,\ldots\}\ra\N$ be defined by
$g(11):=55$ and $g(12):=89$, and
$$ g(2k-1):=2^{k-6}\cdot g(11)\quad\text{and}\quad
   g(2k):=2^{k-6}\cdot g(12) ,$$
for $k\geq 6$; the choice of $g(11)$ and $g(12)$ will become clear below.
Then from $g(11)<g(12)<2\cdot g(11)$ it follows that $g(n)$ is
strictly increasing, and by construction we have $2\cdot g(n)=g(n+2)$,
for all $n\geq 11$.
Having this, we let 
$$ f\cn\{11,12,\ldots\}\ra\N\cn n\mt\min\{f_4(n),g(n)\} ,$$
where is it is easily seen that $g(n)\leq f_4(n)$ if and only if
$n\in\{11\ld 25\}$. Hence it follows that $f(n)$ is strictly 
increasing, and fulfills $2\cdot f(n)\geq f(n+2)$, for all $n\geq 11$.
For convenience, here are a few values:
$$ \begin{array}{|c||cccccccccccc|} \hline
n   & 10 & 11 & 12 & 13 & 14 & 15 & 16 & 17 & 18 & \ldots & 25 & 26 \\
\hline
f_4 & 15 & 55 & 121 & 221 & 364 & 560 & 820 & 1156 & 1581 & \ldots 
    & 8350 & 10075 \\
\hline 
g   &    & 55 & 89 & 110 & 178 & 220 & 356 & 440 & 712 & \ldots 
    & 7040 & 11392 \\
\hline 
\end{array} $$

\absb
By \ref{degformulae} we have $d^\mu\geq f_4(n)\geq f(n)$,
for all $n\geq 11$ and $\mu\in\cP_n^{p\text{-reg}}(4)$.
The decomposition matrices of $\cS_n$ for $n\in\{11,12\}$ 
are known, see the relevant comments in \ref{strategy}, and 
are accessible in the databases mentioned in Section \ref{intro}.
An explicit check, using the computer algebra system 
\GAP{} \cite{GAP}, yields the following:

\abs
If $p\geq 3$,
then for $\mu\in\cM_{11}(\geq 4)$ we have $d^\mu\geq 55$,
where equality is assumed precisely 
for $p=5$ and $\{\mu,\mu^M\}=\{[7,4],[5,3^2]\}$,
and for $\mu\in\cM_{12}(\geq 4)$ we have $d^\mu\geq 89$,
where equality is assumed precisely 
for $p=5$ and $\{\mu,\mu^M\}=\{[6^2],[4^3]\}$.
If $p=2$, then for $\mu\in\cM'_{11}(\geq 4)$ we have $d^\mu\geq 144$,
and for $\mu\in\cM'_{12}(\geq 4)$ we have $d^\mu\geq 164$,
while $d^{\mu_{\text{bs}}(11)}=d^{\mu_{\text{bs}}(12)}=32$.

\abs
Hence, letting $m:=3$ and $n_0:=11$, 
by Theorem \ref{jamesthm} for all $p\geq 3$, independent of $p$, 
we have $d^\mu\geq f(n)$ for all $n\geq 11$ and $\mu\in\cM_n(\geq 4)$.
Similarly, by Theorem \ref{james2thm} for $p=2$
we have $d^\mu\geq f(n)$ for all $n\geq 11$ and $\mu\in\cM'_n(\geq 4)$.

\AbsT{Strategy.}\label{strategy}
We now apply the following strategy:

\absb
Using the function $f\cn\{11,12,\ldots\}\ra\N$ given in \ref{degbnd},
it is easily seen that for $n\geq 11$ we have $f(n)\leq n^3$
if and only if $n\leq 36$. Hence for $n\geq 37$ we conclude that
$d^\mu\leq n^3$ implies $\mu\in\cP_n^{p\text{-reg}}(\leq 3)$.
Moreover, due to the optimal choice of $f_4(n)$ in the sense of \ref{degbnd},
it is to be expected that this bound is rather sharp;
by the results in \ref{low2}--\ref{low7}
it is sharp for $p\in\{2,3,5\}$, and almost for $p=7$.

\abs
Conversely, for all $n\in\N_0$ and $\mu\in\cP_n^{p\text{-reg}}(\leq 3)$, 
by \cite[Cor.2]{JamMinimal} we indeed have $d^\mu\leq n^3$; 
note that this also follows from the degree formulae given in
\ref{deg2formulae} and \ref{deg3formulae}, and in \cite[La.1.21]{BruKle}
for $p=2$, $p=3$ and $p\geq 5$, respectively.

\absb
Hence it remains to consider the cases $n\in\{1\ld 36\}$.
Here, we make use of the decomposition matrices
of $\cS_n$, as far as they are explicitly known.
We briefly report on the current state of the art:
 
\abs
Whenever $n\leq p<2n$, then the Sylow $p$-subgroups of $\cS_n$
are cyclic of order $p$. Hence the decomposition numbers are
described in terms of Brauer trees, see for example \cite[Sect.VII]{Fei}.
Moreover, since $\cS_n$ only has rational-valued irreducible ordinary 
characters, see \cite[Thm.2.1.3]{JamKer}, all the Brauer trees are 
actually stems. This together with unitriangularity properties determines 
these decomposition matrices completely. In turn, this yields the degrees 
$d^\mu$, for all $\mu\in\cP_n^{p\text{-reg}}$.

\abs
For $p\geq 5$ and $2p\leq n<5p$, it was proven in a number of steps 
in \cite{Fay2,Fay1,Ric,Sco} that the adjustment matrix
$\cA_n$, see \ref{hecke}, is the identity matrix.
(We remark that it follows from \cite{GeckTrees} that for $p\geq 2$ 
and $n\leq p<2n$ the adjustment matrix $\cA_n$ also is the identity
matrix, but we will not need this fact.)
Hence the decomposition matrix of $\cS_n$ coincides with
the $\zt_p$-modular decomposition matrix of the generic
Iwahori-Hecke algebra $\cH_n$, where the latter decomposition matrix 
can be computed explicitly by the LLT algorithm, see \ref{llt},
and for all $\mu\in\cP_n^{p\text{-reg}}$ we have
$d^\mu=\dim_{\Q_{\zt_p}}(D_{\zt_p}^\mu)$.

\abs
The $2$-modular and $3$-modular decomposition matrices of $\cS_n$ 
for $n\leq 13$ are given in \cite[p.137--142]{JamLN} and in 
\cite[p.143--152]{JamLN}, respectively. Moreover, those for $p=2$
and $n\in\{14,15\}$ have been determined in \cite{BenDec},
up to an ambiguity for $n=15$, which was solved, together with
the cases $n\in\{16,17\}$ in \cite{JM}.
Moreover, by unpublished work of the author, the $2$-modular 
and $3$-modular decomposition matrices of $\cS_n$ are also known 
for $n\in\{18,19\}$ and $n\in\{14\ld 18\}$, respectively, 
and are accessible in the databases mentioned in Section \ref{intro}.
(In the latter work, a similar technique as in \cite{JM} was employed.
It consists of a mixture of character theoretic and module theoretic
computations, in particular encompassing so-called condensation methods.
A recent description of this approach is given in the thesis \cite{Maa},
which has been written under the author's supervision.)

\absb
By the above, the cases $p\geq 11$ are settled, thus 
we are left with the cases $p\in\{2,3,5,7\}$ and
$n\in\{n_p\ld 36\}$, where $n_2:=20$, $n_3:=19$, $n_5:=25$ and $n_7:=35$. 
Of course, for $\mu\in\cP_n^{p\text{-reg}}(\leq 4)$ we make use of the 
degree formulae in \ref{deg2formulae}, \ref{deg3formulae} and in 
\cite[La.1.21]{BruKle} for $p=2$, $p=3$ and $p\geq 5$, respectively,
as well as of those for $p=2$ and 
$\mu\in\{\mu_{\text{bs}}(n),\mu_{\text{bbs}}(n)\}$ given in \ref{basicspin}. 
Moreover, since the decomposition numbers of the Specht modules
$S^{[n-m,m]}$, where $m\in\{0\ld\fl{\frac{n}{2}}\}$, are
known by \cite{JamI,Jam2}, see also \cite[Thm.24.15]{JamLN},
this also yields the degree $d^\mu$ 
whenever $\mu=[n-m,m]\in\cP_n^{p\text{-reg}}$.
Apart from these cases, given $\mu\in\cP_n^{p\text{-reg}}$, the task is 
to find upper and lower bounds for $d^\mu$, which are good 
enough to ensure that $d^\mu\leq n^3$ or $d^\mu>n^3$, respectively.
 
\abs
To find upper bounds, by \ref{hecke}, in general we have
$d^\mu\leq\dim_{\Q_{\zt_p}}(D_{\zt_p}^\mu)$, where the latter 
dimension can be determined from the $\zt_p$-modular decomposition 
matrix of the generic Iwahori-Hecke algebra $\cH_n$.
As we are dealing with explicit cases, we use the implementation
of the LLT algorithm available in the {\sf SPECHT} package \cite{SPECHT},
available in the computer algebra system {\sf CHEVIE} \cite{CHEVIE},
in order to compute the relevant $\zt_p$-modular decomposition matrices.

\abs
To find lower bounds, we proceed recursively, starting from $n=n_p-1$,
and apply the modular branching rules \ref{modbranch}. 
As it turns out, these estimates are good enough to exclude all cases  
where actually $d^\mu>n^3$. Moreover, as soon as the upper and lower 
bounds found coincide, we have determined $d^\mu$ precisely.
This is particularly strong, in view of the recursive nature of the
procedure, as soon as 
$D^\mu\da_{\cS_{n-1}}\cong\bigoplus_{i=0}^{p-1}D^\mu\da_i$ 
is semi-simple, that is all the $D^\mu\da_i$ are irreducible. 
In particular, as it turns out, these techniques are sufficient to 
determine $d^\mu$ precisely, for all the cases where actually $d^\mu\leq n^3$.

\absb
We have implemented the strategy described above in the computer 
algebra system \GAP{} \cite{GAP}, so that it can be carried out 
explicitly for any fixed prime $p\geq 2$. Below, we are showing 
the results thus obtained for $p\in\{2,3,5,7\}$. (These results are
of course available electronically on request from the author.)

\AbsT{Low-degree representations for $p=2$.}\label{low2}
We obtain the following:

\absb
Generic cases $\mu\in\cP_n^{2\text{-reg}}(\leq 3)$, for $n\in\N_0$,
with degree formulae in \ref{deg2formulae}.
 
\absb
Basic spin cases $\mu=\mu_{\text{bs}}(n)$ and 
second basic spin cases $\mu=\mu_{\text{bbs}}(n)$, 
with degree formulae given in \ref{basicspin},
as well as `almost generic' cases 
$\mu\in\cP_n^{2\text{-reg}}(4)\smin
 \{\mu_{\text{bs}}(n),\mu_{\text{bbs}}(n)\}$,
with degree formulae given in \ref{deg2formulae}:
$$ \begin{array}{|l|l|} \hline
\mu & n \\ \hline \hline
\mu_{\text{bs}}(n) & \{9\ld 30\}\dcup\{32\} \\ 
\mu_{\text{bbs}}(n) & \{8\ld 19\}\dcup\{21\} \\
\hline
\end{array} 
\quad \quad \quad
\begin{array}{|l|l|} \hline
\ov{\mu} & n \\ \hline \hline
(4) & \{11\ld 34\}\dcup\{36\} \\ 
(3,1) & \{10\ld 21\}\dcup\{23\} \\
\hline
\end{array} $$
Note that $\mu_{\text{bs}}(n)$ for $n\leq 8$,
and $\mu_{\text{bbs}}(n)$ for $n\in\{6,7\}$ yield generic cases, and
that $\mu\in\cP_n^{2\text{-reg}}(4)$ for small $n$ 
yields basic spin or second basic spin cases:
$[5,4]=\mu_{\text{bs}}(9)$ and $[6,4]=\mu_{\text{bs}}(10)$, and
$[4,3,1]=\mu_{\text{bbs}}(8)$ and $[5,3,1]=\mu_{\text{bbs}}(9)$.

\absb
`Exceptional' cases 
$\mu\in\cP_n^{2\text{-reg}}\smin(\cP_n^{2\text{-reg}}(\leq 4)
 \cup\{\mu_{\text{bs}}(n),\mu_{\text{bbs}}(n)\})$:
$$ \begin{array}{ccc}
\begin{array}[t]{|l|l|r|} \hline
n & \mu & d^\mu \\ \hline \hline
9 & [4,3,2] & 160 \\ 
\hline
10 & [ 5, 3, 2 ] & 200 \\
   & [ 4, 3, 2, 1 ] & 768 \\
\hline
11 & [ 6, 3, 2 ] & 848 \\
   & [ 5, 4, 2 ] & 416 \\
   & [ 5, 3, 2, 1 ] & 1168 \\
\hline
12 & [ 7, 4, 1 ] & 1408 \\
   & [ 7, 3, 2 ] & 1046 \\
   & [ 6, 4, 2 ] & 416 \\
\hline
\end{array} \rule{1em}{0em} &
\begin{array}[t]{|l|l|r|} \hline
n & \mu & d^\mu \\ \hline \hline
13 & [ 8, 5 ] & 560 \\
   & [ 8, 4, 1 ] & 1572 \\
   & [ 6, 5, 2 ] & 1728 \\
\hline
14 & [ 9, 5 ] & 560 \\
   & [ 8, 4, 2 ] & 2510 \\
   & [ 7, 5, 2 ] & 2016 \\
\hline
15 & [ 10, 5 ] & 1288 \\
   & [ 9, 6 ] & 1912 \\
\hline
\end{array} \rule{1em}{0em} &
\begin{array}[t]{|l|l|r|} \hline
n & \mu & d^\mu \\ \hline \hline
16 & [ 11, 5 ] & 1288 \\
   & [ 10, 6 ] & 1912 \\
   & [ 8, 6, 2 ] & 4096 \\
\hline
17 & [ 12, 5 ] & 3808 \\
   & [ 11, 6 ] & 4488 \\
\hline
18 & [ 13, 5 ] & 3808 \\
   & [ 12, 6 ] & 4488 \\
\hline
20 & [ 15, 5 ] & 6972 \\
\hline
\end{array} \\
\end{array} $$

\AbsT{Low-degree representations for $p=3$.}\label{low3}
We obtain the following:

\absb
Generic cases $\mu\in\cP_n^{3\text{-reg}}(\leq 3)$, for $n\in\N_0$,
with degree formulae in \ref{deg3formulae}.
 
\absb
`Almost generic' cases $\mu\in\cP_n^{3\text{-reg}}(4)$, where
$\mu^M\not\in\cP_n^{3\text{-reg}}(\leq 3)$,
with degree formulae given in \ref{deg3formulae};
we also note the cases where $\mu^M\in\cP_n^{3\text{-reg}}(\leq 4)$,
which only occurs for $n\leq 9$:
$$ \begin{array}{|l|l|} \hline
\ov{\mu} & n \\ \hline \hline
(4) & \{10\ld 33\}\dcup\{36\} \\ 
(3,1) & \{9\ld 17\} \\
(2^2) & \{10\ld 33\}\dcup\{36\} \\ 
(2,1^2) & \{8\ld 17\} \\ \hline
\end{array} 
\quad \quad \quad 
\begin{array}{|l|llll|} \hline
\ov{\mu} & 6 & 7 & 8 & 9 \\ \hline \hline
(4) & - & - & () & () \\
(3,1) & - & (1) & (1) & \\
(2^2) & - & (1^2) & (2,1) & (2,1) \\
(2,1^2) & (2) & (2) & (2,1^2) & \\ \hline
\end{array} $$

\absb
`Exceptional' cases
$\mu,\mu^M\in\cP_n^{3\text{-reg}}\smin\cP_n^{3\text{-reg}}(\leq 4)$,
where we record the lexicographically largest of $\mu$ and $\mu^M$:
$$ \begin{array}{ccc}
\begin{array}[t]{|l|l|r|} \hline
n & \mu & d^\mu \\ \hline \hline
10 & [5,3,1^2] & 567 \\ \hline
11 & [6,4,1] & 693 \\
   & [6,3,2] & 252 \\
   & [6,3,1^2] & 791 \\
   & [6,2^2,1] & 714 \\ \hline
12 & [7,5] & 297 \\
   & [7,4,1] & 1013 \\
   & [7,3,2] & 252 \\
   & [7,3,1^2] & 1431 \\
   & [7,2^2,1] & 1728 \\
   & [6,3,2,1] & 1428 \\ \hline
\end{array} \rule{1em}{0em} &
\begin{array}[t]{|l|l|r|} 
\hline
n & \mu & d^\mu \\ \hline \hline
13 & [8,5] & 428 \\
   & [8,4,1] & 1275 \\
   & [8,3,2] & 792 \\
   & [8,2^2,1] & 1938 \\
   & [7,3^2] & 924 \\
   & [7,3,2,1] & 1428 \\ \hline
14 & [9,5] & 428 \\
   & [9,3,2] & 1287 \\
   & [8,6] & 1000 \\
   & [8,3^2] & 1716 \\ \hline
15 & [10,5] & 1548 \\
   & [10,3,2] & 1287 \\ 
   & [9,6] & 1428 \\ 
   & [9,3^2] & 1716 \\ \hline 
\end{array} \rule{1em}{0em} &
\begin{array}[t]{|l|l|r|} \hline
n & \mu & d^\mu \\ \hline \hline
16 & [11,5] & 2108 \\ 
   & [11,3,2] & 3003 \\
   & [10,6] & 1428 \\
   & [9,7] & 3417 \\ \hline
17 & [12,5] & 2108 \\
   & [12,3,2] & 4368 \\
   & [10,7] & 4845 \\ \hline
18 & [13,5] & 5508 \\
   & [13,3,2] & 4368 \\
   & [11,7] & 4845 \\ \hline
20 & [15,5] & 7105 \\ \hline
\end{array} \\
\end{array} $$

\AbsT{Low-degree representations for $p=5$.}\label{low5}
We obtain the following:

\absb
Generic cases $\mu\in\cP_n^{5\text{-reg}}(\leq 3)$, for $n\in\N_0$,
with degree formulae in \cite[La.1.21]{BruKle}.
 
\absb
`Almost generic' cases $\mu\in\cP_n^{5\text{-reg}}(4)$, where
$\mu^M\not\in\cP_n^{5\text{-reg}}(\leq 3)$,
with degree formulae given in \cite[La.1.21]{BruKle};
we also note the cases where $\mu^M\in\cP_n^{5\text{-reg}}(\leq 4)$,
which only occurs for $n\leq 10$:
$$ \begin{array}{|l|l|} \hline
\ov{\mu} & n \\ \hline \hline
(4) & \{8\ld 33\}\dcup\{36\} \\ 
(3,1) & \{8\ld 17\} \\
(2^2) & \{8\ld 20\}\dcup\{23\} \\ 
(2,1^2) & \{8\ld 16\} \\
(1^4) & \{9\ld 32\} \dcup\{35\} \\ \hline
\end{array} 
\quad \quad \quad 
\begin{array}{|l|lllll|} \hline
\ov{\mu} & 6 & 7 & 8 & 9 & 10 \\ \hline \hline
(4) & - & - & (2^2) & & \\
(3,1) & - & (2) & & & \\
(2^2) & (3) & (3) & (4) & & \\
(2,1^2) & () & (2,1) & (2,1^2) & & \\ 
(1^4) & (1) & (1^2) & (1^3) & (1^4) & (1^4) \\ \hline
\end{array} $$

\absb
`Exceptional' cases
$\mu,\mu^M\in\cP_n^{5\text{-reg}}\smin\cP_n^{5\text{-reg}}(\leq 4)$,
where we record the lexicographically largest of $\mu$ and $\mu^M$:
$$ \begin{array}{ccc}
\begin{array}[t]{|l|l|r|} \hline
n & \mu & d^\mu \\ \hline \hline
9 & [ 4^2, 1 ] & 83 \\ \hline
10 & [ 5^2 ] & 34 \\
 & [ 5, 4, 1 ] & 217 \\
 & [ 5, 3, 2 ] & 450 \\
 & [ 5, 3, 1^2 ] & 266 \\
 & [ 5, 2^2, 1 ] & 525 \\
 & [ 4^2, 1^2 ] & 300 \\ \hline
11 & [ 6, 5 ] & 89 \\ 
 & [ 6, 4, 1 ] & 372 \\ 
 & [ 6, 3, 2 ] & 605 \\ 
 & [ 6, 3, 1^2 ] & 266 \\
 & [ 6, 2^2, 1 ] & 1100 \\ 
 & [ 6, 2, 1^3 ] & 252 \\
 & [ 5^2, 1 ] & 285 \\ 
 & [ 5, 4, 2 ] & 980 \\ 
 & [ 5, 4, 1^2 ] & 1035 \\
 & [ 5, 3, 2, 1 ] & 1330 \\
 & [ 5, 2^3 ] & 825 \\ \hline
\end{array} \rule{0.4em}{0em} &
\begin{array}[t]{|l|l|r|} 
\hline
n & \mu & d^\mu \\ \hline \hline
12 &  [ 7, 5 ] & 144 \\
 & [ 7, 4, 1 ] & 372 \\
 & [ 7, 3, 1^2 ] & 1266 \\
 & [ 7, 2^2, 1 ] & 1506 \\
 & [ 7, 2, 1^3 ] & 462 \\
 & [ 6^2 ] & 89 \\
 & [ 6, 5, 1 ] & 835 \\
 & [ 6, 3^2 ] & 1650 \\
 & [ 6, 3, 2, 1 ] & 1596 \\
 & [ 6, 2, 1^4 ] & 2100 \\ 
 & [ 5^2, 2 ] & 1265 \\ 
 & [ 5^2, 1^2 ] & 1320 \\ \hline
13 & [ 8, 5 ] & 144 \\
 & [ 8, 4, 1 ] & 2001 \\
 & [ 8, 2, 1^3 ] & 792 \\
 & [ 7, 6 ] & 233 \\
 & [ 7, 5, 1 ] & 1495 \\
 & [ 7, 2^2, 1^2 ] & 924 \\ 
 & [ 6^2, 1 ] & 924 \\ \hline
14 & [ 9, 5 ] & 1001 \\
 & [ 9, 2, 1^3 ] & 1287 \\
 & [ 8, 6 ] & 377 \\
 & [ 8, 5, 1 ] & 1639 \\
 & [ 8, 2^2, 1^2 ] & 1716 \\
 & [ 7^2 ] & 233 \\
 & [ 7, 6, 1 ] & 2652 \\ \hline
\end{array} \rule{0.4em}{0em} &
\begin{array}[t]{|l|l|r|} 
\hline
n & \mu & d^\mu \\ \hline \hline
15 & [ 10, 5 ] & 1625 \\
 & [ 10, 2, 1^3 ] & 1287 \\
 & [ 9, 6 ] & 377 \\
 & [ 9, 2^2, 1^2 ] & 1716 \\
 & [ 8, 7 ] & 610 \\
 & [ 7^2, 1 ] & 2885 \\ \hline
16 & [ 11, 5 ] & 2445 \\
 & [ 11, 2, 1^3 ] & 3003 \\
 & [ 10, 6 ] & 3640 \\
 & [ 9, 7 ] & 987 \\
 & [ 8^2 ] & 610 \\ \hline
17 & [ 12, 5 ] & 3265 \\
 & [ 12, 2, 1^3 ] & 4368 \\
 & [ 10, 7 ] & 987 \\
 & [ 9, 8 ] & 1597 \\ \hline
18 & [ 13, 5 ] & 3265 \\
 & [ 10, 8 ] & 2584 \\
 & [ 9, 9 ] & 1597 \\ \hline
19 & [ 11, 8 ] & 2584 \\
 & [ 10, 9 ] & 4181 \\ \hline
20 & [ 11, 9 ] & 6765 \\
 & [ 10, 10 ] & 4181 \\ \hline
21 & [ 12, 9 ] & 6765 \\ \hline
\end{array} \rule{1em}{0em} \\
\end{array} $$

\AbsT{Low-degree representations for $p=7$.}\label{low7}
We obtain the following:

\absb
Generic cases $\mu\in\cP_n^{7\text{-reg}}(\leq 3)$, for $n\in\N_0$,
with degree formulae in \cite[La.1.21]{BruKle}.
 
\absb
`Almost generic' cases $\mu\in\cP_n^{7\text{-reg}}(4)$, where
$\mu^M\not\in\cP_n^{7\text{-reg}}(\leq 3)$,
with degree formulae given in \cite[La.1.21]{BruKle};
we also note the cases where $\mu^M\in\cP_n^{7\text{-reg}}(\leq 4)$,
which only occurs for $n\leq 9$:
$$ \begin{array}{|l|l|} \hline
\ov{\mu} & n \\ \hline \hline
(4) & \{8\ld 34\} \\ 
(3,1) & \{7\ld 16\} \\
(2^2) & \{7\ld 20\} \\ 
(2,1^2) & \{9\ld 16\} \\
(1^4) & \{9\ld 32\} \dcup\{35\} \\ \hline
\end{array} 
\quad \quad \quad 
\begin{array}{|l|lllll|} \hline
\ov{\mu} & 5 & 6 & 7 & 8 & 9 \\ \hline \hline
(4) & - & - & - & & \\
(3,1) & - & - & (2^2) & & \\
(2^2) & - & (3) & (3,1) & & \\
(2,1^2) & - & (2) & (2,1) & (2,1) & \\ 
(1^4) & () & (1) & (1) & (1^3) & (1^4) \\ \hline
\end{array} $$

\absb
`Exceptional' cases
$\mu,\mu^M\in\cP_n^{7\text{-reg}}\smin\cP_n^{7\text{-reg}}(\leq 4)$,
where we record the lexicographically largest of $\mu$ and $\mu^M$:
$$ \begin{array}{ccc}
\begin{array}[t]{|l|l|r|} \hline
n & \mu & d^\mu \\ \hline \hline
8 & [ 3^2, 2 ] & 42 \\ \hline
9 & [ 4^2, 1 ] & 84 \\
 & [ 4, 3, 2 ] & 168 \\
 & [ 3^3 ] & 42 \\ \hline 
10 & [ 5, 5 ] & 42 \\
 & [ 5, 4, 1 ] & 199 \\
 & [ 5, 3, 2 ] & 384 \\
 & [ 5, 3, 1^2 ] & 567 \\
 & [ 5, 2^2, 1 ] & 525 \\
 & [ 5, 2, 1^3 ] & 448 \\
 & [ 4^2, 2 ] & 252 \\
 & [ 4, 3^2 ] & 210 \\ \hline
11 & [ 6, 5 ] & 131 \\
 & [ 6, 4, 1 ] & 693 \\
 & [ 6, 3, 2 ] & 485 \\
 & [ 6, 3, 1^2 ] & 1232 \\
 & [ 6, 2^2, 1 ] & 626 \\
 & [ 6, 2, 1^3 ] & 924 \\
 & [ 6, 1^5 ] & 252 \\
 & [ 5^2, 1 ] & 199 \\
 & [ 5, 4, 2 ] & 835 \\
 & [ 5, 4, 1^2 ] & 1155 \\
 & [ 5, 3^2 ] & 594 \\
 & [ 4^2, 3 ] & 462 \\ \hline
\end{array} \rule{1em}{0em} &
\begin{array}[t]{|l|l|r|} 
\hline
n & \mu & d^\mu \\ \hline \hline
12 & [ 7, 5 ] & 286 \\
 & [ 7, 4, 1 ] & 1353 \\
 & [ 7, 2, 1^3 ] & 1398 \\
 & [ 7, 1^5 ] & 462 \\
 & [ 6^2 ] & 131 \\
 & [ 6, 5, 1 ] & 1155 \\
 & [ 6, 4, 2 ] & 1320 \\
 & [ 6, 3^2 ] & 1079 \\
 & [ 5^2, 2 ] & 1034 \\
 & [ 5^2, 1^2 ] & 1354 \\
 & [ 4^3 ] & 462 \\ \hline
13 & [ 8, 5 ] & 507 \\
 & [ 8, 4, 1 ] & 2145 \\
 & [ 8, 1^5 ] & 792 \\
 & [ 7, 6 ] & 417 \\
 & [ 7, 3^2 ] & 1079 \\
 & [ 7, 1^6 ] & 924 \\
 & [ 6^2, 1 ] & 1286 \\ \hline
\end{array} \rule{1em}{0em} &
\begin{array}[t]{|l|l|r|} 
\hline
n & \mu & d^\mu \\ \hline \hline
14 & [ 9, 5 ] & 728 \\
 & [ 9, 4, 1 ] & 2366 \\
 & [ 9, 1^5 ] & 792 \\
 & [ 8, 6 ] & 924 \\
 & [ 8, 1^6 ] & 924 \\ 
 & [ 7^2 ] & 417 \\
 & [ 6^2, 2 ] & 2354 \\ \hline
15 & [ 10, 5 ] & 728 \\
 & [ 10, 1^5 ] & 2002 \\
 & [ 9, 6 ] & 1652 \\
 & [ 9, 1^6 ] & 3003 \\
 & [ 8, 7 ] & 1341 \\ \hline
16 & [ 11, 5 ] & 2548 \\
 & [ 11, 1^5 ] & 3003 \\
 & [ 10, 6 ] & 2380 \\
 & [ 9, 7 ] & 2993 \\
 & [ 8^2 ] & 1341 \\ \hline
17
 & [ 12, 5 ] & 3808 \\
 & [ 12, 1^5 ] & 4368 \\
 & [ 11, 6 ] & 2380 \\
 & [ 9, 8 ] & 4334 \\ \hline
18 & [ 13, 5 ] & 5507 \\
 & [ 9^2 ] & 4334  \\ \hline
\end{array}
\end{array} $$

\Rem
For the sake of completeness, we remark that to classify 
the cases $d^\mu\leq n$ and $d^\mu\leq n^2$ it turns out that we may 
just use the functions 
$$ f_2\cn\Z\ra\Z\cn n\mt\frac{n^2-5n+2}{2}\quad\text{and}\quad
   f_3\cn\Z\ra\Z\cn n\mt\frac{n^3-9n^2+14n}{6} $$
from \ref{degformulae}, respectively.
Proceeding as in \ref{degbnd} we find that the conditions
of \ref{jamesthm} and \ref{james2thm} are fulfilled for $n_0=8$ and
$n_0=13$, respectively. Since $f_2(n)>n$ for all $n\geq 8$, and
$f_3(n)>n^2$ for all $n\geq 15$, it only remains to consider the cases
$n\leq 14$ explicitly. In this range, all decomposition matrices of $\cS_n$
are known, so that the strategy in \ref{strategy} simplifies considerably.
Here are the results, where for $p\geq 3$ we only record the 
lexicographically larger of $\mu$ and $\mu^M$, and for all 
representations mentioned degree formulae are given in
\ref{basicspin} and \ref{deg2formulae}, in \ref{deg3formulae},
and in \cite[La.1.21]{BruKle} for $p=2$, $p=3$ and $p\geq 5$, respectively,
or are accessible in the databases mentioned in Section \ref{intro}.

\absb
As for $d^\mu\leq n$, apart from the generic cases 
$\mu\in\cP_n^{p\text{-reg}}(\leq 1)$ for $n\in\N_0$,
we have the following `exceptional' cases
$\mu,\mu^M\in\cP_n^{p\text{-reg}}\smin\cP_n^{p\text{-reg}}(\leq 1)$:
For $p=2$ we get $\{[3,2],[4,2],[5,3]\}$, which are 
basic spin representations;
for $p=3$ we get $\{[4,1^2]\}$; and
for $p\geq 5$ we get $\{[2^2],[3,2],[3,3]\}$.

\absb
As for $d^\mu\leq n^2$, apart from the generic cases 
$\mu\in\cP_n^{p\text{-reg}}(\leq 2)$ for $n\in\N_0$,
we have the following cases where 
$\mu,\mu^M\in\cP_n^{p\text{-reg}}\smin\cP_n^{p\text{-reg}}(\leq 2)$:

\absc
For $p=2$, we get the basic spin and second basic spin cases 
$\mu_{\text{bs}}(n)$ and $\mu_{\text{bbs}}(n)$,
and the `almost generic' cases 
$\mu\in\cP_n^{2\text{-reg}}(3)\smin 
\{\mu_{\text{bs}}(n),\mu_{\text{bbs}}(n)\}$
as follows:
$$ \begin{array}{|l|l|} \hline
\mu & n \\ \hline \hline
\mu_{\text{bs}}(n) & \{7\ld 18\} \\ 
\mu_{\text{bbs}}(n) & \{6\ld 9\} \\
\hline
\end{array}
\quad \quad \quad
\begin{array}{|l|l|} \hline
\ov{\mu} & n \\ \hline \hline
(3) & \{9,10,11,12,\} \\
(2,1) & \{8,9\} \\ 
\hline
\end{array} $$
Note that $\mu_{\text{bs}}(n)$ for $n\leq 6$ yields generic cases, and 
that $\mu\in\cP_n^{2\text{-reg}}(3)$ for small $n$ 
yields basic spin or second basic spin cases:
$[4,3]=\mu_{\text{bs}}(7)$ and $[5,3]=\mu_{\text{bs}}(8)$, and
$[3,2,1]=\mu_{\text{bbs}}(6)$ and $[4,2,1]=\mu_{\text{bbs}}(7)$.
There are no `exceptional' cases.

\absc
For $p=3$, the `almost generic' cases
$\mu\in\cP_n^{3\text{-reg}}(3)$ such that
$\mu^M\not\in\cP_n^{3\text{-reg}}(\leq 2)$, 
and the `exceptional' cases $\mu,\mu^M\not\in\cP_n^{3\text{-reg}}(\leq 3)$
are, respectively:
$$ \begin{array}{|l|l|} \hline
\ov{\mu} & n \\ \hline \hline
(3) & \{8\ld 13\} \\
(2,1) & \{7\ld 12\} \\ 
\hline
\end{array} 
\quad \quad \quad
\begin{array}{|l|l|r|} 
\hline
n & \mu & d^\mu \\ \hline \hline
10 & [6,4] & 90 \\ \hline
12 & [8,4] & 131 \\ \hline
\end{array} $$

\absc
For $p=5$, similarly, the `almost generic' and `exceptional' cases are:
$$ \begin{array}{ccc} 
\begin{array}[t]{|l|l|} \hline
\ov{\mu} & n \\ \hline \hline
(3) & \{6\ld 11\}\dcup\{14\} \\
(2,1) & \{7,8\} \\ 
(1^3) & \{7\ld 11\} \\ 
\hline
\end{array} \rule{1em}{0em} &
\begin{array}[t]{|l|l|r|} 
\hline
n & \mu & d^\mu \\ \hline \hline
8 & [4^2] & 13 \\ \hline
9 & [5,4] & 34 \\ 
 & [5,1^4] & 70 \\ \hline
10 & [6,4] & 55 \\ 
 & [6,1^4] & 70 \\
 & [5^2] & 34 \\ \hline
\end{array} \rule{1em}{0em} &
\begin{array}[t]{|l|l|r|} 
\hline
n & \mu & d^\mu \\ \hline \hline
11 & [7,4] & 55 \\
 & [6,5] & 89 \\ \hline
12 & [7,5] & 144 \\
 & [6^2] & 89 \\ \hline
13 & [8,5] & 144 \\ \hline
\end{array} 
\end{array} $$

\absc
For $p=7$, similarly, the `almost generic' and `exceptional' cases are:
$$ \begin{array}{ccc}  
\begin{array}[t]{|l|l|} \hline
\ov{\mu} & n \\ \hline \hline
(3) & \{6\ld 11\} \\
(2,1) & \{6,7,8\} \\ 
(1^3) & \{8\ld 11\} \\ 
\hline
\end{array} \rule{1em}{0em} &
\begin{array}[t]{|l|l|r|} 
\hline
n & \mu & d^\mu \\ \hline \hline
7 & [3^2,1] & 21 \\ \hline
8 & [4^2] & 14 \\
 & [4,2^2] & 56 \\
 & [3^2,2] & 42 \\ \hline
\end{array} \rule{1em}{0em} &
\begin{array}[t]{|l|l|r|} 
\hline
n & \mu & d^\mu \\ \hline \hline
9 & [5,4] & 42 \\
 & [5,1^4] & 70 \\
 & [3^3] & 42 \\ \hline
10 & [6,4] & 89 \\
 & [5^2] & 42 \\ \hline
12 & [6^2] & 131 \\ \hline
\end{array} 
\end{array} $$


\abs
{\sc 
Institut f\"ur Mathematik \\
Friedrich-Schiller-Universit\"at Jena \\
Ernst-Abbe-Platz 2, 07737 Jena, Germany} \\
{\sf juergen.manfred.mueller@uni-jena.de}


\begin{thebibliography}{References}
\setlength{\parskip}{-0.2em}

\bibitem{Ari} {\sc S. Ariki}:
On the decomposition numbers of the Hecke algebra of $G(m,1,n)$, 
J. Math. Kyoto Univ. 36 (4), 1996, 789--808.                              

\bibitem{Ben} {\sc D. Benson}:
Spin modules for symmetric groups, 
J. London Math. Soc. (2) 38 (2), 1988, 250--262.

\bibitem{BenDec} {\sc D. Benson}: 
Some remarks on the decomposition numbers for the symmetric groups, 
Proc. of Symposia in Pure Math. 47, 1987, 381--394.

\bibitem{BesOls} {\sc C. Bessenrodt, J. Olsson}:
On residue symbols and the Mullineux conjecture,
J. Algebraic Combin. 7 (3), 1998, 227--251.

\bibitem{CTblLib} {\sc T. Breuer}:
{\sf GAP}-package {\sf CTblLib}---The {\sf GAP} character table library,
Version 1.2.2, 2013,
{\sf http://www.gap-system.org/Packages/ctbllib.html}.

\bibitem{BruKle2} {\sc J. Brundan, A. Kleshchev}:
Representation theory of symmetric groups and their double covers,
in: Groups, combinatorics \& geometry, Durham, 2001, 31--53.

\bibitem{BruKle} {\sc J. Brundan, A. Kleshchev}:
Representations of the symmetric group which are irreducible over subgroups,
J. Reine Angew. Math. 530, 2001, 145--190.                                           
\bibitem{Dan} {\sc S. Danz}:
Vertices of low-dimensional simple modules for symmetric groups, 
Comm. Algebra 36 (12), 2008, 4521--4539. 

\bibitem{DipJam} {\sc R. Dipper, G. James}: 
Representations of Hecke algebras of general linear groups, 
Proc. London Math. Soc. 52 (3), 1986, 20--52.

\bibitem{FawBreSax} {\sc J. Fawcett, E. O'Brien, J. Saxl}:
Regular orbits of symmetric and alternating groups,
to appear in J. Algebra.

\bibitem{Fay2} {\sc M. Fayers}:
James's Conjecture holds for weight four blocks of Iwahori-Hecke
algebras, J. Algebra 317 (2), 2007, 593--633.

\bibitem{Fay1} {\sc M. Fayers}:
Decomposition numbers for weight three blocks of symmetric groups
and Iwahori-Hecke algebras, 
Trans. Amer. Math. Soc. 360 (3), 2008, 1341--1376.

\bibitem{Fei} {\sc W. Feit}:
The representation theory of finite groups, North-Holland, 1982.

\bibitem{For} {\sc B. Ford}:
Irreducible restrictions of representations of the symmetric groups,
Bull. London Math. Soc. 27 (5), 1995, 453--459.

\bibitem{ForKle} {\sc B. Ford, A. Kleshchev}:
A proof of the Mullineux conjecture, Math. Z. 226 (2), 1997, 267--308.
                          
\bibitem{FraRobThr} {\sc J. Frame, G. Robinson, R. Thrall}:
The hook graphs of the symmetric groups, Canadian J. Math. 6, 1954, 316--324.

\bibitem{GAP}
{\sc The {\sf GAP} Group}: {\sf GAP}---Groups, Algorithms,
Programming---a system for computational discrete algebra,
Version 4.7.7, 2015, {\sf http://www.gap-system.org}.

\bibitem{GeckTrees} {\sc M. Geck}:
Brauer trees of Hecke algebras, Comm. Algebra 20 (10), 1992, 2937--2973.


                           
\bibitem{HisMal} {\sc G. Hiss, G. Malle}:
Low-dimensional representations of quasi-simple groups,
LMS J. Comput. Math. 4, 2001, 22–-63; corrigenda,
LMS J. Comput. Math. 5, 2002, 95-–126. 

\bibitem{JamMinimal} {\sc G. James}:
On the minimal dimensions of irreducible representations of symmetric 
groups, Math. Proc. Cambridge Philos. Soc. 94 (3), 1983, 417--424.

\bibitem{JamIII} {\sc G. James}:
On the decomposition matrices of the symmetric groups III,
J. Algebra 71 (1), 1981, 115-–122.

\bibitem{JamLN} {\sc G. James}:
The representation theory of the symmetric groups,
Lecture Notes in Mathematics 682, Springer, 1978.

\bibitem{JamI} {\sc G. James}:
On the decomposition matrices of the symmetric groups I,                   
J. Algebra 43 (1), 1976, 42–-44.                                          

\bibitem{Jam2} {\sc G. James}:
Representations of the symmetric groups over the field of order $2$,
J. Algebra 38 (2), 1976, 280-–308.                                        

\bibitem{JamKer} {\sc G. James, A. Kerber}:
The representation theory of the symmetric group,
Encyclopedia of Mathematics and its Applications 16,
Addison-Wesley, 1981.

\bibitem{JamWil} {\sc G. James, A. Williams}:
Decomposition numbers of symmetric groups by induction, 
J. Algebra 228 (1), 2000, 119--142.

\bibitem{ModularAtlas}
{\sc C. Jansen, K. Lux, R. Parker, R. Wilson}:
An Atlas of Brauer Characters, Clarendon Press, Oxford, 1995.

\bibitem{JanSei} {\sc J. Jantzen, G. Seitz}:
On the representation theory of the symmetric groups,
Proc. London Math. Soc. (3) 65 (3), 1992, 475--504.

\bibitem{Kas} {\sc M. Kashiwara}:
On crystal bases of the $q$-analogue of universal enveloping algebras,
Duke Math. J. 63 (2), 1991, 465--516.                                     

\bibitem{KleBook} {\sc A. Kleshchev}:
Linear and projective representations of symmetric groups,
Cambridge Tracts in Mathematics 163, Cambridge University Press, 2005.

\bibitem{Kle} {\sc A. Kleshchev}:
On decomposition numbers and branching coefficients for symmetric and
special linear groups, Proc. London Math. Soc. (3) 75 (3), 1997, 497--558.  

\bibitem{KleIII} {\sc A. Kleshchev}:
Branching rules for modular representations of symmetric groups III:
some corollaries and a problem of Mullineux,
J. London Math. Soc. (2) 54 (1), 1996, 25--38.

\bibitem{KleII} {\sc A. Kleshchev}:
Branching rules for modular representations of symmetric groups II,
J. Reine Angew. Math. 459, 1995, 163-–212.

\bibitem{KleI} {\sc A. Kleshchev}:
Branching rules for modular representations of symmetric groups I,
J. Algebra 178 (2), 1995, 493--511.

\bibitem{Kle0} {\sc A. Kleshchev}:
On restrictions of irreducible modular representations of semisimple
algebraic groups and symmetric groups to some natural subgroups I,
Proc. London Math. Soc. (3) 69 (3), 1994, 515--540.

\bibitem{LasLecThi} {\sc A. Lascoux, B. Leclerc, J. Thibon}:
Hecke algebras at roots of unity and crystal bases of quantum affine
algebras, Comm. Math. Phys. 181 (1), 1996, 205--263.

\bibitem{Lub} {\sc F. L\"ubeck}:
Small degree representations of finite Chevalley groups in 
defining characteristic, LMS J. Comput. Math. 4, 2001, 135--169.

\bibitem{Maa} {\sc L. Maas}:
Modular spin characters of symmetric groups,
Ph.D. Thesis, University of Essen, 2011.
 
\bibitem{SPECHT} {\sc A. Mathas}:
{\sf SPECHT}---decomposition matrices for the Hecke algebras of type $A$,
University of Sydney, 1997, 
{\sf http://www.maths.usyd.edu.au/u/ mathas/specht}.

\bibitem{Mat} {\sc A. Mathas}:
Iwahori-Hecke algebras and Schur algebras of the symmetric group,
University Lecture Series 15, American Mathematical Society, 1999. 

\bibitem{CHEVIE} {\sc J. Michel}:
{\sf CHEVIE}---development version of the {\sf GAP} part of {\sf CHEVIE},
{\sf http://webusers.imj-prg.fr/$\sim$jean.michel/chevie.html}.

\bibitem{JM} {\sc J. M\"uller}:
The $2$-modular decomposition matrices of the symmetric groups 
$S_{15}$, $S_{16}$, and $S_{17}$, Comm. Algebra 28, 2000, 4997--5005.

\bibitem{Mul} {\sc G. Mullineux}: 
Bijections of $p$-regular partitions and $p$-modular 
irreducibles of the symmetric groups,
J. London Math. Soc. (2) 20 (1), 1979, 60--66.

\bibitem{Pee} {\sc M. Peel}:
Hook representations of the symmetric groups,
Glasgow Math. J. 12, 1971, 136-–149.                                         
  
\bibitem{Ric} {\sc M. Richards}:
Some decomposition numbers for Hecke algebras of general linear groups,
Math. Proc. Cambridge Philos. Soc. 119 (3), 1996, 383--402.               

\bibitem{Sco} {\sc J. Scopes}:
Symmetric group blocks of defect two,
Quart. J. Math. Oxford Ser. (2) 46 (182), 1995, 201--234.                 

\bibitem{VarVas} {\sc M. Varagnolo, E. Vasserot}:
On the decomposition matrices of the quantized Schur algebra,
Duke Math. J. 100 (2), 1999, 267--297.                                    

\bibitem{Wal} {\sc D. Wales}:
Some projective representations of $S_{n}$, 
J. Algebra 61 (1), 1979, 37--57.


\bibitem{ModularAtlasProject}
{\sc R. Wilson, J. Thackray, R. Parker, F. Noeske, J. M{\"u}ller, \!
F. L{\"u}beck, C. Jansen, G. Hiss, T. Breuer}:
The Modular Atlas Project,
{\sf http://www.math.rwth-aachen.de/homes/MOC}.


\end{thebibliography}
\end{document}